\documentclass[12pt]{article}
\usepackage{epsfig}
\usepackage{diagbox}
\usepackage{multirow}
\usepackage{array}

\usepackage{amsfonts}
\usepackage{verbatim}
\usepackage{latexsym}
\usepackage{cite}
\usepackage{graphicx,amsmath}
\usepackage{amsmath,amssymb,color}
\usepackage{amsmath,amsfonts,amssymb,epsfig,subfigure}
\usepackage{mathrsfs}
\usepackage{stmaryrd}
\usepackage{graphicx,amsmath}
\usepackage{latexsym}
\usepackage{verbatim}
\usepackage{listings}
\usepackage[framed,numbered,autolinebreaks,useliterate]{mcode}
\usepackage{listings}
\def\nn{\nonumber}
\def\g{\gamma}

\makeatletter \@addtoreset{equation}{section}
\renewcommand{\thesection}{\arabic{section}}
\renewcommand{\theequation}{\thesection.\arabic{equation}}

\begin{document}

\newcommand{\E}{\mathbb{E}}
\newcommand{\PP}{\mathbb{P}}
\newcommand{\RR}{\mathbb{R}}
\newcommand{\I}{\mathbb{I}}
\newcommand{\cd}{(\cdot)}
\newcommand{\LL}{\mathbb{L}}
\newcommand{\Se}{\mathbb{S}}

\newtheorem{theorem}{Theorem}[section]
\newtheorem{lemma}{Lemma}[section]
\newtheorem{coro}{Corollary}[section]
\newtheorem{defn}{Definition}[section]
\newtheorem{assp}{Assumption}
\newtheorem{expl}{Example}[section]
\newtheorem{remark}{Remark}[section]

\newtheorem{thm}{Lemma}[subsection]
\renewcommand{\thethm}{A.\arabic{thm}}

\newtheorem{corollary}{Corollary}[subsection]
\renewcommand{\thecorollary}{D.\arabic{corollary}}

\newcommand\tq{{\scriptstyle{3\over 4 }\scriptstyle}}
\newcommand\qua{{\scriptstyle{1\over 4 }\scriptstyle}}
\newcommand\hf{{\textstyle{1\over 2 }\displaystyle}}
\newcommand\hhf{{\scriptstyle{1\over 2 }\scriptstyle}}

\newcommand{\eproof}{\indent\vrule height6pt width4pt depth1pt\hfil\par\medbreak}

\newcommand{\disp}{\displaystyle}
\newcommand{\ad}{&\!\!\!\disp}
\newcommand{\aad}{&\disp}
\newcommand{\barray}{\begin{array}{ll}}
\newcommand{\earray}{\end{array}}

\def\t{\triangle}  \def\ts{\triangle^*}
\def\ra{\rightarrow}

\def\a{\alpha} \def\b{\beta} \def\d{\delta}\def\g{\gamma}
\def\e{\varepsilon} \def\z{\zeta} \def\y{\eta} \def\o{\theta}
\def\vo{\vartheta} \def\k{\kappa} \def\l{\lambda} \def\m{\mu} \def\n{\nu}
\def\x{\xi}  \def\r{\rho} \def\s{\sigma}
\def\p{\phi} \def\f{\varphi}   \def\w{\omega}
\def\q{\surd} \def\i{\bot} \def\h{\forall} \def\j{\emptyset}
 \def\si{s_{ij}}
\def\be{\beta} \def\de{\delta} \def\up{\upsilon} \def\eq{\equiv}
\def\ve{\vee} \def\we{\wedge}
\def\SS{{\cal S}}
\def\F{{\cal F}} \def\Le{{\cal L}}
\def\Lem{{\cal L}^{\varepsilon, M}}
\def\T{\tau} \def\G{\Gamma}  \def\D{\Delta} \def\O{\Theta} \def\L{\Lambda}
\def\X{\Xi} \def\S{\Sigma} \def\W{\Omega}
\def\M{\partial} \def\N{\nabla} \def\Ex{\exists} \def\K{\times}
\def\V{\bigvee} \def\U{\bigwedge}
\def\X{\Phi}
\def\lf{\left} \def\rt{\right}

\def\tl{\tilde}
\def\trace{\hbox{\rm trace}}
\def\diag{\hbox{\rm diag}}
\def\for{\quad\hbox{for }}
\def\refer{\hangindent=0.3in\hangafter=1}

\newcommand\wD{\widehat{\D}}
\makeatletter \@addtoreset{equation}{section}
\renewcommand{\thesection}{\arabic{section}}
\renewcommand{\theequation}{\thesection.\arabic{equation}}
\allowdisplaybreaks

\newcommand{\dis}{\displaystyle}
\def\nn{\nonumber}

\def\bb{\begin}
\def\bc{\begin{center}}       \def\ec{\end{center}}
\def\ba{\begin{array}}        \def\ea{\end{array}}
\def\be{\begin{equation}}     \def\ee{\end{equation}}
\def\bea{\begin{eqnarray}}    \def\eea{\end{eqnarray}}
\def\beaa{\begin{eqnarray*}}  \def\eeaa{\end{eqnarray*}}
\def\la{\label}

\title{
 Explicit numerical approximation  for logistic models with regime switching in finite and infinite horizons
}
\author{   Xiaoyue Li\thanks{School of Mathematics and Statistics,
Northeast Normal University, Changchun, Jilin, 130024, China.
 Email: lixy209@nenu.edu.cn.},
~~~~Hongfu Yang\thanks{School of Mathematics and Statistics,
Northeast Normal University, Changchun, Jilin, 130024, China. Email: yanghf783@nenu.edu.cn.}
}
\date{}
\maketitle
\setcounter{page}{1} \pagestyle{plain}
\vspace{-1em}

\begin{abstract}
 The stochastic logistic model with regime switching is an important model in the ecosystem. While analytic solution  to this model is positive, current numerical methods are unable to preserve such boundaries in the approximation. So, proposing appropriate numerical method for solving this model which  preserves positivity and  dynamical behaviors  of the model's solution is very important.
 In this paper, we present a positivity preserving truncated Euler-Maruyama scheme  for this model, which taking advantages of being explicit and easily implementable.
 Without additional restriction conditions, strong convergence of the numerical algorithm  is studied, and 1/2  order convergence rate  is obtained.
 In the particular case of this model without switching the  first order strong convergence rate is obtained.
  Furthermore,   the approximation of long-time dynamical properties is realized, including the stochastic permanence, extinctive and stability in distribution.
  Some simulations and examples are provided to confirm the theoretical results and demonstrate the validity of the approach.

 \vspace{3mm} \noindent {\bf Keywords.}
Stochastic logistic model, Markov chain, Explicit scheme, Strong convergence,   Stochastic permanence, Stationary distribution

\vspace{3mm} \noindent {\bf 2000 MR Subject Classification.} ~60H10.

\end{abstract}

\section{ Introduction}\label{one1}
  In this manuscript, we consider the numerical  approximation of stochastic logistic model with environmental fluctuations
  described by the following
  switching diffusion system (SDS)
\be\label{logistic1}
\mathrm{d}x(t)=x(t)\Big[\big(b(r(t))-a(r(t))x(t)\big)\mathrm{d}t + \sigma(r(t))\mathrm{d}B(t)\Big]
\ee
with an initial value $x(0)=x_0\in \mathbb{R}_+:=(0,+\infty)$, $r(0)=\ell\in \mathbb{S}:=\{1, 2,\ldots,m\}$, $x(t)$ is the population size at time $t$.
 $r(t )$ is a right-continuous Markov chain with finite state space $\mathbb{S}$ and the generator $\Gamma=(\gamma_{ij})_{m\times m}$ satisfying $\gamma_{ij}\geq 0$ for $i\neq j$, $\sum_{j\in \mathbb{S}}\gamma_{ij}=0$ for each $i\in \mathbb{S}$.  $B(t)$ is a scalar standard Brownian motion, which is independent of $r(t)$.
 $b(i)$ and $a(i)$ represents the intrinsic growth rate and the intraspecific competition coefficient in regime $i$, respectively.
 $\sigma^2(i)$ is a constant representing the intensity of
the white noise in regime $i$.

This model plays an important role in biomathematics applications,
the dynamical behaviors of SDS \eqref{logistic1} and its related stochastic models  have been investigated recently in \cite{hu20181AMM,Li2016JMAA,Li2009JCAM,Li2009DCDS}. It is known that a unique strong solution exists for SDS \eqref{logistic1}, and that non-negativity of the initial value is preserved, see e.g. \cite{Li2011,Li2016JMAA}.
 Once we attempt to simulate SDS \eqref{logistic1} using classical discretization methods, see e.g.
\cite{eKloeden1992book}, we face three  difficulties:
\begin{itemize}

\item[$\bullet$] In general, these methods do not preserve positivity and therefore are not well
defined when directly applied to SDS \eqref{logistic1}.

\item[$\bullet$] The  drift  term is not globally Lipschitz continuous and therefore standard assumptions
required for  strong convergence, see e.g. \cite{eKloeden1992book}, do not hold.

\item[$\bullet$]
Despite the convergence analysis, how to approximate long-time behaviors of SDS \eqref{logistic1} is significant and challenging.
\end{itemize}
 Our primary objective is to construct easily implementable preserving positivity numerical solutions and prove that they converge to the true solution of the original SDS \eqref{logistic1},   moreover, realize the approximation of long-time dynamical properties including the stochastic permanence, extinctive and stability in distribution.

In recent years, a few Euler-Maruyama (EM) discretization schemes have been developed for diffusion systems and SDSs including the implicit EM method \cite{Higham02,eKloeden1992book}, the tamed EM method \cite{Hutzenthaler12a,Sabanis13a,Sabanis16a,Nguyen18}, the tamed Milstein method \cite{Wang13a}, the stopped EM method \cite{Liu13}  and the truncated EM method \cite{Li2019IMA,mao15a}, to mention a few. In these EM methods,
the approximation can potentially escape the domain of the exact solutions of systems. Consequently,  in order to close the gap, a lot of effort has focused on deriving schemes staying in restricted domains for diffusion systems with non-Lipschitz continuous coefficients \cite{Alfonsi05,Alfonsi13, Neuenkirch14,Dereich12, Chassagneux16}.
 Several  modified  EM  methods have been developed such as the implicit schemes \cite{Dereich12, Neuenkirch14} and the explicit EM  schemes \cite{Chassagneux16,Liu13},
 in the context of mathematical finance, a thorough overview of these can be found in \cite{Kloeden13}.
 A now classical trick is to apply a suitable Lamperti transform in order to obtain
diffusion systems  with constant diffusion coefficient, thereby translating all the non-smoothness to
the drift. In the context of non-globally Lipschitz coefficients, this idea, introduced
by Alfonsi \cite{Alfonsi05}, was further exploited in \cite{Alfonsi13, Neuenkirch14} to obtain strong convergence rates
for implicit ``Lamperti-Euler" schemes, in particular for the CIR and the Ait-Sahalia
models, and for scalar diffusion systems with one-sided Lipschitz continuous drift and constant
diffusion coefficients\cite{Neuenkirch14}. Recently,
 in the context of ecology, Mao, Wei and Wiriyakraikul \cite{Mao2021} have established a positive preserving truncated EM method for stochastic Lotka-Volterra competition model but without any convergence rate of the algorithm.  Chen, Gan and Wang \cite{CGW2021} have proposed the Lamperti smoothing truncation scheme that can preserve the domain of the original SDEs and proved a mean-square convergence rate of order one. These modified EM methods have shown their abilities to approximate the solutions of nonlinear diffusion systems. However, to the best of our knowledge,  these modified EM methods still cannot handle the convergence of  nonlinear SDS \eqref{logistic1}.

Motivated by Lamperti transform \cite{Alfonsi05, Neuenkirch14} and truncation approaches \cite{Li2019IMA,mao15a},
our key idea is to transform the original SDS \eqref{logistic1} using the Lamperti transformation
into an new SDS, i.e.,
applying
 It\^{o}'s formula to $y(t)=\log x(t)$ gives
\begin{align}\la{log_eq4}
\mathrm{d} y(t)
= \big(\beta(r(t))-a(r(t)) \mathrm{e}^{y(t)}\big)\mathrm{d}t+\sigma(r(t)) \mathrm{d}B(t),
\end{align}
where $\beta(i):=b(i)-\frac{\sigma^2(i)}{2}$ for any $i\in \mathbb{S}$.
 The transformed SDS \eqref{log_eq4} is
then approximated by a truncation EM scheme
and transforming back yields the   preserving positivity  numerical schemes for the original SDS \eqref{logistic1}, which has a computational cost of the same order as the
classical  EM scheme.  This allows us to prove rate  of convergence for the original SDS \eqref{logistic1}, and the numerical solutions keep the underlying excellent properties of the exact solution  of SDS \eqref{logistic1}.
  Here, we extend that work in several ways:
\begin{itemize}

\item[$\bullet$]  Constructing an easily implementable scheme to maintain the positive path of the exact solution for nonlinear SDS  \eqref{logistic1}. The scheme shares some of the features of the truncation schemes family.

\item[$\bullet$] The explicit EM approximate solution will converge to the exact solution with order 1/2  for nonlinear SDS   \eqref{logistic1}.

\item[$\bullet$] Considering the maximum error in the discretization points, we prove that the explicit EM scheme for the stochastic logistic models without regime switching  strongly converges with rate one.

\item[$\bullet$] Without extra restrictions the numerical solutions of the appropriate explicit scheme stay in step of dynamical properties with the exact solutions.
\end{itemize}

The rest of the paper is organized as follows. Section \ref{two2} gives some preliminary results on certain
properties of the exact solutions.  Section \ref{section3} constructs an explicit scheme, and optimal convergence rate is obtained.   Section \ref{strong2} focuses on the analyses the stochastic permanence and extinction of the SDS \eqref{logistic1}. The other explicit scheme is constructed preserving the stochastic permanence and extinction.
 Section \ref{sections5} analyses the stability of the SDS \eqref{logistic1} in distribution yielding an invariant measure   $\mu(\cdot\times\cdot)$, and explicit scheme
preserving the stability in distribution and a numerical invariant measure, which tends to $\mu(\cdot\times\cdot)$ as the step
size tends to 0. Section \ref{section numerical} presents a couple of examples to illustrate our results. Section \ref{section_7} reconstructs an explicit scheme,
 and yields the strong convergence with rate one. Some  examples  are given to illustrate the availability of this scheme.

\section{Preliminaries}\label{two2}

Suppose that both $r(\cdot)$ and $B(\cdot)$ are
 defined on the complete
probability space $( \Omega, \cal{F}$, $\{ {\cal{F}}_{t}\} _{t \geq 0}, \PP )$ with a filtration ${\{ {\cal{F}}_{t}\}} _{t \geq 0}$ satisfying the usual conditions (i.e., it is right continuous and $\mathcal{F}_0$ contains all $\mathbb{P}$-null sets). For each $i \in \mathbb{S}$, both $a(i)$ and $b(i)$ are nonnegative constants, $\sigma(i)$ is a constant.   Let   $\mathbb{E}$ denotes the expectation corresponding to $\PP$ and $|\cdot|$ denote the Euclidean norm in $\RR$. The generator of $\{r(t)\}_{t\geq 0}$ is denoted by $\Gamma$, so that for a sufficiently small $\delta>0$,
$$
\mathbb{P}\{r(t+\delta)=j|r(t)=i\}=\left\{
\begin{array}{ll}
\gamma_{ij}\delta+o(\delta),&~~\hbox{if}~~~i\neq j,\\
1+\gamma_{ii}\delta+o(\delta),&~~\hbox{if}~~~i=j,
\end{array}
\right.
$$
where $o(\delta)$ satisfies $\lim_{\delta\rightarrow 0} {o(\delta)}/{\delta}=0$. Here $\gamma_{ij}\geq 0$ is the transition rate from $i$ to $j$ if $i\neq j$ while $\gamma_{ii}=-\sum_{i\neq j}\gamma_{ij}$. It is well known that almost every sample path of $r(t)$ is right-continuous step functions with a finite number of simple jumps in any finite interval of $\mathbb{R}_+$ (cf. \cite{mao11a}).
 If $\mathbb{D}$ is a set, its indicator function is denoted by $I_{\mathbb{D}}$, namely $I_{\mathbb{D}}(x)=1$ if $x\in \mathbb{D}$ and $0$ otherwise.
For convenience, we let $C_u$ and $C$ denote two generic positive real constants respectively, whose value may change in different appearances, where $C_u$ is  dependent on $u$. And  let $\mathcal{N}(0,1)$ denotes   the standard normal distribution. For any $c=(c(1),\ldots, c(m))$(or $=(c_1,\ldots, c_m)$), define $\hat{c}=\min_{i\in \mathbb{S}}c(i)$,  $\check{c}=\max_{i\in \mathbb{S}}c(i)$ and $|\breve{c}|=\max_{i\in \mathbb{S}}|c(i)|$. We state a useful lemma which can be found in \cite{Li2011}.

\begin{lemma}[\!\!\cite{Li2011}]\la{log_le1}
There exists a unique continuous positive solution $x(t)$ to SDS \eqref{logistic1} for any initial value $x(0)=x_0 > 0$ and $r(0)=\ell\in \mathbb{S}$, which is global and represented by
\begin{align*}
x(t)=\frac{x_0\exp\bigg\{\dis\int_{0}^{t}\beta(r(s)) \mathrm{d}s
+\int_{0}^{t}\sigma(r(s)) \mathrm{d}B(s)\bigg\}}{1+x_0\dis\int_{0}^{t}a(r(s))
\exp\bigg\{\int_{0}^{s}\beta(r(u)) \mathrm{d}u+\int_{0}^{s}\sigma(r(u)) \mathrm{d}B(u)\bigg\}\mathrm{d}s}.
\end{align*}
\end{lemma}
By virtue of Lemma \ref{log_le1}, for any $p>0$, the solution $x(t)$ of SDS \eqref{logistic1}  with any initial value $x_0\in \mathbb{R}_+$, $ \ell\in \mathbb{S}$, satisfies
\begin{align}\la{log_eq3}
\sup_{0\leq t\leq T}\mathbb{E}\big[x^p(t)\big]
\leq  C_T,~~\forall~T>0.
\end{align}
 Now, we give the boundedness of its inverse moment. The inverse moment plays an important role in the analysis of convergence rate for the numerical scheme.
\begin{lemma}\la{inverse_x}
For any $p>0$,  then we have
\begin{align*}
\sup_{0\leq t\leq T}\mathbb{E}\big[x^{-p}(t)\big]
\leq  C_T,~~\forall~T>0.
\end{align*}
\end{lemma}
{\bf  Proof.}~~By virtue of Lemma \ref{log_le1}, the solution $x(t)$ with positive initial value will remain in $\mathbb{R}_{+}$ for all $t\geq 0$ with probability 1. Define $
U(t)=x^{-1}(t)
$ on $t\geq 0$,
we derive from \eqref{logistic1} that
\begin{align*}
\mathrm{d} U( t )
=& \Big[ - \Big( \beta (r( t ) ) - \frac{\sigma ^ { 2 }(r( t ) )}{2} \Big) U ( t ) + a(r( t ) ) \Big] \mathrm{d} t - \sigma (r( t ) ) U ( t ) \mathrm{d} B ( t ).
\end{align*}
Define a Lyapunov function $V(u)=(1+u)^{p}$ for any $p>0$. Using the method of Lyapunov function analysis, we could obtain the required assertion. The left proof is rather standard and hence is omitted.
\eproof
%
%
\begin{remark}
By  virtue of Lemma \ref{inverse_x} and \eqref{log_eq3}, for any $q\in \mathbb{R}$,  then the transformed SDS \eqref{log_eq4} has the following exponential integrability property
\begin{align*}
\sup_{0\leq t\leq T}\mathbb{E}\big[\mathrm{e}^{q y(t)}\big]
\leq  C_T,~~\forall~T>0.
\end{align*}
 Obviously, we also have the following property
\begin{align}\la{integrability}
\sup_{0\leq t\leq T}\mathbb{E}\big[|y(t)|^p\big]\leq \sup_{0\leq t\leq T}\mathbb{E}\big[\mathrm{e}^{p |y(t)|}\big]
\leq C_T,
\end{align}
for any $p>0$ and $T>0$.
\end{remark}

\section{Convergence rate}\la{section3}
In this section, we aim to construct an easily implementable explicit scheme and show the rate of convergence. The rate is optimal similar to the standard results of the  EM scheme for SDSs with globally Lipschitz coefficients, see \cite[p.115]{Mao06}. Given a stepsize $\Delta>0$ and let $t_k=k\Delta$, $r_k=r(t_k)$ for $k\geq0$, and one-step transition probability matrix
$P(\Delta)=\big(P_{ij}(\Delta)\big)_{m\times m}=\exp(\Delta \Gamma)$. The discrete Markov chain $\{r_k, k=0, 1,\ldots\}$ can be simulated by the techniques in \cite[p.112]{Mao06}.\! Throughout the article,
 $C$ and $C_T$ are independent of $\Delta$ and $k$.

To define appropriate numerical solutions, we firstly propose an explicit scheme to approximate the exact solution of SDS \eqref{log_eq4}. For any given stepsize $\Delta\in(0, 1)$, define a  truncated EM   scheme by
\begin{align}\la{log_eq5*}
\left\{
\begin{array}{ll}
Z_0=\log x_0,&\\
\bar{Z}_{k+1}=Z_{k}+\big(\beta(r_k)-a(r_k) \mathrm{e}^{Z_k}\big)\Delta +\sigma(r_k) \Delta B_k,& \\
Z_{k+1}=\bar{Z}_{k+1}\wedge \log(K\Delta^{-\theta } ),
\end{array}
\right.
\end{align}
for any integer $k\geq0$, where $\Delta B_k=B(t_{k+1})-B(t_k)$,  $K\geq x_0\vee 1$ is a constant independent of the iteration order $k$ and the stepsize $\Delta$, we use the convention $\theta =+\infty$ if  $\check{a}=0$ and $\theta \in (0, 1/2]$ otherwise.
Transforming back, i.e.
\begin{align}\la{TEM_1}
X_k=\mathrm{e}^{Z_k},~~~~~~k=0,1,\ldots,
\end{align}
gives a strictly positive approximation of the original  SDS \eqref{logistic1}.
 Obviously,  we have
\begin{align}\la{logeq:4.1}
Z_{k}=\bar{Z}_{k}\wedge \log(K\Delta^{-\theta } )\leq \bar{Z}_{k},~~~~ 0< X_k\leq K\Delta^{-\theta }.
\end{align}
To proceed, we define $\bar{Z}_{\Delta}(t)$ and $X_{\Delta}(t)$ by
$
\bar{Z}_{\Delta}(t):=\bar{Z}_k, X_{\Delta}(t):=X_k,~ \forall~t\in[t_k,t_{k+1}).
$
\begin{remark}
The $\check{a}=0$ implies that  $a(i)\equiv 0$, then we have $\theta =+\infty$ and  $K\Delta^{-\theta }\equiv +\infty$ for any $\Delta\in (0, 1]$. Thus, $\bar{Z}_k=Z_k$ and \eqref{logeq:4.1} hold always for any $\Delta\in (0, 1]$.
\end{remark}


In order to study the
rate of convergence of numerical solutions $\{X_k\}_{k\geq 0}$, we first give  the following lemmas.
\begin{lemma}\la{log_le3.1}
For any $p>0$, the truncated EM scheme defined by \eqref{log_eq5*} has the property that
\begin{align*}
\sup_{\Delta\in (0, 1)}\sup_{0\leq k\leq \lfloor T/\Delta\rfloor}\mathbb{E}\big[\mathrm{e}^{pZ_{k}}\big]
\leq \sup_{\Delta\in (0, 1)}\sup_{0\leq k\leq \lfloor T/\Delta\rfloor}\mathbb{E}\big[\mathrm{e}^{p\bar{Z}_{k}}\big]
\leq  C_T,~~\forall~T>0,
\end{align*}
 where $\lfloor T/\Delta\rfloor$  represents the integer part of $T/\Delta$.
\end{lemma}
{\bf  Proof.}~~
Since $a(i)\geq0$ for any $i\in \mathbb{S}$, we know that
\begin{align*}
\bar{Z}_{k+1}\leq Z_{k}+ \beta(r_k) \Delta +\sigma(r_k) \Delta B_k\leq \bar{Z}_{k}+ [\check{\beta}]^+ \Delta +\sigma(r_k) \Delta B_k.
\end{align*}
Then
 $$
\bar{Z}_{k}\leq  \log x_0+[\check{\beta}]^+k\Delta
  +\sum_{j=0}^{k-1} \sigma(r_{j})\Delta B_{j}
$$
for any integer $k\geq 1$. Thus, for any $p>0$  we have
\begin{align}\la{log_eq9}
\mathbb{E}\big[\mathrm{e}^{p\bar{Z}_{k}}\big]
\leq&\mathbb{E}\Big[\exp\Big(p Z_0 +p[\check{\beta} ]^{+}T+p \sum_{j=0}^{k-1} \sigma(r_{j})\Delta B_{j}\Big)\Big]
= x_0^p\mathrm{e}^{p T[\check{\beta}]^{+}}\mathbb{E}
\Big[\exp\Big(p \sum_{j=0}^{k-1} \sigma(r_{j})\Delta B_{j}\Big)\Big]\nn\\
=&x_0^p\mathrm{e}^{p T[\check{\beta}]^{+}}
\mathbb{E}\bigg[\exp\Big(\frac{p^2}{2} \sum_{j=0}^{k-1}  \sigma^2(r_j)\Delta \Big)\exp\Big(-\frac{p^2}{2}\sum_{j=0}^{k-1} \sigma^2(r_{j})\Delta+p \sum_{j=0}^{k-1} \sigma(r_{j})\Delta B_{j}\Big)\bigg]\nn\\
\leq&x_0^p\exp\Big(p[\check{\beta} ]^{+}T
+\frac{p^2|\breve{\sigma}|^2}{2}T\Big)
\mathbb{E}\bigg[\exp\Big(-\frac{p^2}{2}\sum_{j=0}^{k-1} \sigma^2(r_{j})\Delta+p \sum_{j=0}^{k-1} \sigma(r_{j})\Delta B_{j}\Big)\bigg].
\end{align}
On the other hand,
\begin{align*}
\exp\Big(-\frac{p^2}{2}\sum_{j=0}^{k-1} \sigma^2(r_{j})\Delta+p \sum_{j=0}^{k-1} \sigma(r_{j})\Delta B_{j}\Big)
=\prod_{j=0}^{k-1}\Theta_j,
\end{align*}
where
$$
\Theta_j=\exp\Big(-\frac{p^2}{2} \sigma^2(r_j)\Delta+p \sigma(r_{j})\Delta B_{j}\Big).
$$
Then we have
\begin{align*}
\mathbb{E}\bigg[\prod_{j=0}^{k-1}\Theta_j\bigg]
= \mathbb{E}\bigg[\mathbb{E}\bigg(\prod_{j=0}^{k-1}\Theta_j\bigg|\mathcal{F}_{t_{k-1}}\bigg)\bigg]
= \mathbb{E}\bigg[\prod_{j=0}^{k-2}\Theta_j\mathbb{E}\big(\Theta_{k-1}\big|\mathcal{F}_{t_{k-1}}\big)\bigg].
\end{align*}
Obviously,
\begin{align*}
\mathbb{E}\big(\Theta_{k-1}\big|\mathcal{F}_{t_{k-1}}\big)
=\mathbb{E}\bigg[\sum_{i\in \mathbb{S}}I_{\{r_k=i\}}\Theta^i_{k-1}\Big|\mathcal{F}_{t_{k-1}}\bigg]
=\sum_{i\in \mathbb{S}}I_{\{r_{k-1}=i\}}\mathbb{E}\big(\Theta^i_{k-1}\big|\mathcal{F}_{t_{k-1}}\big),
\end{align*}
where
$$
\Theta^i_{k-1}=\exp\Big(-\frac{p^2}{2} \sigma^2(i)\Delta+p \sigma(i)\Delta B_{k-1}\Big),~~~~~~~i\in \mathbb{S}.
$$
Note that $\Delta B_{k-1}=B(t_{k})-B(t_{k-1})$
is independent of $\mathcal{F}_{t_{k-1}}$, by \cite[Lemma 3.2, p. 104]{Mao06}, we can derive that
\begin{align*}
\mathbb{E}\big(\Theta^i_{k-1}\big|\mathcal{F}_{t_{k-1}}\big)
=\exp\Big(-\frac{p^2}{2} \sigma^2(i)\Delta\Big)\mathbb{E}\Big[\exp\Big( p \sigma(i)\Delta B_{k-1}\Big)\Big]=1.
\end{align*}
Hence,
$
\mathbb{E}\big[\prod_{j=0}^{k-1}\Theta_j\big]
 = \mathbb{E}\big[\prod_{j=0}^{k-2}\Theta_j\big].
$
Repeating this procedure, we obtain that
\begin{align*}
\mathbb{E}\bigg[\exp\Big(-\frac{p^2}{2}\sum_{j=0}^{k-1} \sigma^2(r_{j})\Delta+p \sum_{j=0}^{k-1} \sigma(r_{j})\Delta B_{j}\Big)\bigg]
=1.
\end{align*}
The above equality together with \eqref{log_eq9} implies
\begin{align*}
\mathbb{E}\big[\mathrm{e}^{p\bar{Z}_{k}}\big]
\leq x_0^p\exp\Big(p[\check{\beta} ]^{+}T
+\frac{p^2|\breve{\sigma}|^2}{2}T\Big).
\end{align*}
The proof is complete.\eproof

\begin{lemma}\la{leB}
For any $L\geq0$ and integer $m\geq 0$, we have
\begin{align*}
\E\Big[\mathrm{e}^{\frac{|\Delta B_k|^2}{4\Delta}}|\Delta B_k|^{m}\Big]
 = \sqrt{\frac{2}{\pi}} \Gamma\Big(\frac{m+1}{2}\Big)\big(4\Delta\big)^{\frac{m}{2}},~~~~~
 \E\big[\mathrm{e}^{L |\Delta B_k|}    |\Delta B_k|^m\big|\mathcal{F}_{t_k}\big]
 \leq C \Delta^{\frac{m}{2}},
\end{align*}
where $\Gamma(\cdot)$ is the Gamma function.
\end{lemma}
{\bf  Proof.}~~
For any integer $m\geq 0$ and due to $\Delta B_k\sim \mathcal{N}(0, \Delta)$, we deduce that
 \begin{align*}
\E\Big[\mathrm{e}^{\frac{|\Delta B_k|^2}{4\Delta}}|\Delta B_k|^{m}\Big]
= &\frac{2}{\sqrt{2\pi \Delta}}\int_{0}^{+\infty}x^{m}\mathrm{e}^{-\frac{x^2}{4\Delta}} \mathrm{d}x
=  \frac{4\Delta}{\sqrt{2\pi \Delta}}\int_{0}^{+\infty}(4\Delta)^{\frac{m-1}{2}}x^{\frac{m-1}{2}}\mathrm{e}^{-x} \mathrm{d}x\nn\\
= &  \frac{ \big(4\Delta\big)^{\frac{m+1}{2}}}{\sqrt{2\pi \t}}\int_{0}^{+\infty}x^{\frac{m+1}{2}-1}\mathrm{e}^{-x} \mathrm{d}x= \sqrt{\frac{2}{\pi}} \Gamma\Big(\frac{m+1}{2}\Big)\big(4\Delta\big)^{\frac{m}{2}} ,
\end{align*}
which implies that
\begin{align*}
 & \E\Big[ \exp\big(L |\Delta B_k|\big) |\Delta B_k|^m \big|\mathcal{F}_{t_k}\Big]
\leq  \E\Big[\exp\big( L^2\Delta+\frac{|\Delta B_k|^2}{4\Delta}\big) |\Delta B_k|^{m} \Big|\mathcal{F}_{t_k}\Big]\nn\\
=&\mathrm{e}^{ L^2\Delta} \E\Big[\mathrm{e}^{\frac{|\Delta B_k|^2}{4\Delta}}|\Delta B_k|^{m}\big|\mathcal{F}_{t_k}\Big]
= \mathrm{e}^{  L^2\Delta}  \E\Big[\mathrm{e}^{\frac{|\Delta B_k|^2}{4\Delta}}|\Delta B_k|^{3m}\Big] \leq C \Delta^{\frac{m}{2}}.
\end{align*}
The proof is complete.
\eproof

\begin{lemma}\la{inv_Z}
For any $p> 0$, the truncated EM scheme defined by \eqref{log_eq5*} has the property that
\begin{align*}
\sup_{\Delta\in (0, 1)}\sup_{0\leq k\leq \lfloor T/\Delta\rfloor}\mathbb{E}\big[\mathrm{e}^{-pZ_{k}}\big]
\leq  C_T,~~~\forall~T>0.
\end{align*}
\end{lemma}
{\bf  Proof.}~~Using the Taylor formula, we obtain that
\begin{align}\la{yhf0602}
\mathrm{e}^{-\bar{Z}_{k+1}}\leq& \mathrm{e}^{-Z_{k}}-\mathrm{e}^{-Z_{k}}(\bar{Z}_{k+1}-Z_{k})+\frac{1}{2}\mathrm{e}^{-Z_{k}}(\bar{Z}_{k+1}-Z_{k})^2
+\frac{1}{6}\mathrm{e}^{-Z_{k}}\mathrm{e}^{|\bar{Z}_{k+1}-Z_{k}|}|\bar{Z}_{k+1}-Z_{k}|^3\nn\\
\leq & \mathrm{e}^{-Z_{k}}\bigg[1+\big(a(r_k) \mathrm{e}^{Z_k}-\beta(r_k)\big)\Delta-\sigma(r_k) \Delta B_k+\frac{1}{2} (\bar{Z}_{k+1}-Z_{k})^2\nn\\
&~~~~~~+\frac{1}{6} \exp{\Big(\big(\check{a}  \mathrm{e}^{Z_k}+|\breve{\beta}|\big)\Delta+|\breve{\sigma}| |\Delta B_k|\Big)}\Big(\big(\check{a}  \mathrm{e}^{Z_k}+|\breve{\beta}|\big)\Delta+|\breve{\sigma}| |\Delta B_k|\Big)^3\bigg]\nn\\
\leq & \mathrm{e}^{-Z_{k}}\bigg[1+\big(a(r_k) \mathrm{e}^{Z_k}-\beta(r_k)\big)\Delta-\sigma(r_k) \Delta B_k
+\frac{\sigma^2(r_k) (\Delta B_k)^2}{2}\nn\\
&~~~~~~+ \big(\check{a}  K \Delta^{-\theta }+|\breve{\beta}|  \big)^2 \Delta^2
-\sigma(r_k)\Delta B_k\big(a(r_k) \mathrm{e}^{Z_k}-\beta(r_k)\big)\Delta \nn\\
&~~~~~~+\frac{1}{6} \exp{\Big(\big(\check{a}  K\Delta^{-\theta }+|\breve{\beta}|\big)\Delta+|\breve{\sigma}| |\Delta B_k|\Big)}\Big(\big(\check{a}  K\Delta^{-\theta }+|\breve{\beta}|\big)\Delta+|\breve{\sigma}| |\Delta B_k|\Big)^3\bigg]\nn\\
\leq & \mathrm{e}^{-Z_{k}}\bigg[1+\big(a(r_k) \mathrm{e}^{Z_k}-\beta(r_k)\big)\Delta
+ C \Delta^{2(1-\theta )}
+\frac{\sigma^2(r_k) (\Delta B_k)^2}{2}-\sigma(r_k) \Delta B_k\nn\\
&~~~~~~
-\sigma(r_k)\Delta B_k\big(a(r_k) \mathrm{e}^{Z_k}-\beta(r_k)\big)\Delta\nn\\
&~~~~~~+\frac{1}{6}\big(\check{a}  K\vee|\breve{\beta}|\vee|\breve{\sigma}|\big)^3\mathrm{e}^{(\check{a}  K\vee|\breve{\beta}|) \Delta^{1-\theta }} \exp\big( |\breve{\sigma}||\Delta B_k|\big)\big( \Delta^{1-\theta }+ |\Delta B_k|\big)^3\bigg]\nn\\
\leq & \mathrm{e}^{-Z_{k}}\Big[1+\big(a(r_k) \mathrm{e}^{Z_k}-\beta(r_k)\big)\Delta+
C \Delta^{2(1-\theta )}
+\frac{\sigma^2(r_k) (\Delta B_k)^2}{2}-\sigma(r_k) \Delta B_k\nn\\
&~~~~~~
-\sigma(r_k)\Delta B_k\big(a(r_k) \mathrm{e}^{Z_k}-\beta(r_k)\big)\Delta
+\mathcal{U}_k\Big],
\end{align}
where
$$
\mathcal{U}_k=C\exp\big( |\breve{\sigma}| |\Delta B_k|\big)\big( \Delta^{3(1-\theta )}+  |\Delta B_k|^3\big),
$$
which implies that
\begin{align}
(1+\mathrm{e}^{-\bar{Z}_{k+1}})^p\leq (1+\mathrm{e}^{-Z_{k}})^p(1+\zeta_k)^p,
\end{align}
where
\begin{align*}
\zeta_k=& \frac{\mathrm{e}^{-Z_{k}}}{1+\mathrm{e}^{-Z_{k}}}\Big[ \big(a(r_k) \mathrm{e}^{Z_k}-\beta(r_k)\big)\Delta+C \Delta^{2(1-\theta )}
+\frac{\sigma^2(r_k) (\Delta B_k)^2}{2}-\sigma(r_k) \Delta B_k\nn\\
&~~~~~~~~~~~~~~~~~~
-\sigma(r_k)\Delta B_k\big(a(r_k) \mathrm{e}^{Z_k}-\beta(r_k)\big)\Delta
+\mathcal{U}_k\Big],
\end{align*}
and we can see that $\varsigma_k>-1$. For the given constant $p>2$, choose an integer $m$ such that $2m<p\leq2(m+1)$. It follows from \cite[Lemma 3.3]{li2018jde} and \eqref{logeq:5.45} that
\begin{align}\la{log}
&\E \Big[\big(1+\mathrm{e}^{-\bar{Z}_{k+1}}\big)^{p}  \big|\mathcal{F}_{t_k}\Big]\nn\\
\leq&\big(1+\mathrm{e}^{-Z_{k}}\big)^{p} \bigg\{1+ p \E\big[\zeta_k\big|\mathcal{F}_{t_k}\big] + \frac{p(p-1)}{2} \E\big[\zeta_k^2\big|\mathcal{F}_{t_k}\big] + \mathbb{E}\big[P_m(\zeta_k)\zeta_k^3 \big|\mathcal{F}_{t_k}\big]\bigg\},
\end{align}
where $P_m(x)$ represents a $m$th-order polynomial of $x$  with coefficients depending only on $p$, and $m$ is an
integer. Noticing that the increment $\t B_k$ is independent of $\mathcal{F}_{t_k}$,  we
derive that
\begin{align}\la{properties}
\mathbb{E} \big(|\Delta B_{k} |^{2j}\big|\mathcal{F}_{t_k}\big)=(2j-1)!!\Delta^{j},~~~
\mathbb{E} \big((\Delta B_{k} )^{2j-1}\big|\mathcal{F}_{t_k}\big)=0,~~j=1,2,\ldots
\end{align}
and using \eqref{logeq:4.1} and Lemma \ref{leB}, we compute
\begin{align}\la{Y_4}
\mathbb{E}\big[\zeta_k \big|\mathcal{F}_{t_k}\big]
=&\frac{\mathrm{e}^{-Z_{k}}}{1+\mathrm{e}^{-Z_{k}}}\Big[ \big(a(r_k) \mathrm{e}^{Z_k}-\beta(r_k)\big)\Delta+C \Delta^{2(1-\theta )}
+\frac{\sigma^2(r_k)\Delta}{2}
+\mathbb{E}\big[\mathcal{U}_k\big|\mathcal{F}_{t_k}\big]\Big]\nn\\
\leq&\frac{\mathrm{e}^{-Z_{k}}}{1+\mathrm{e}^{-Z_{k}}}\Big( a(r_k) \mathrm{e}^{Z_k}\Delta+C \Delta \Big)\leq C\t.
\end{align}
and
\begin{align} \la{Y_5}
\mathbb{E}\big[\zeta_k^2 \big|\mathcal{F}_{t_k}\big]
=& \frac{\mathrm{e}^{-2Z_{k}}}{(1+\mathrm{e}^{-Z_{k}})^2}\mathbb{E}\Big[ \Big( a(r_k) \mathrm{e}^{Z_k}\t-\beta(r_k) \Delta+C \Delta^{2(1-\theta )}
+\frac{\sigma^2(r_k) (\Delta B_k)^2}{2}\nn\\
&~~~~~~~~~~~~~~~~~~
- \sigma(r_k) \Delta B_k-\sigma(r_k)\Delta B_k\big(a(r_k) \mathrm{e}^{Z_k}-\beta(r_k)\big)\Delta
+\mathcal{U}_k\Big)^2\Big|\mathcal{F}_{t_k}\Big]\nn\\
\leq&\frac{C\mathrm{e}^{-2Z_{k}}}{(1+\mathrm{e}^{-Z_{k}})^2}\Big[ a^2(r_k) \mathrm{e}^{2Z_k}\t^2 + \Delta^{\frac{3}{2}}
+ \sigma^2(r_k)\t
+\mathbb{E}\big[\mathcal{U}_k^2\big|\mathcal{F}_{t_k}\big]\Big]\leq C\t.
\end{align}
To estimate $\mathbb{E}\big[P_m(\zeta_k)\zeta_k^3 \big|\mathcal{F}_{t_k}\big]$, we begin with $\mathbb{E}\big[\zeta_k^3 \big|\mathcal{F}_{t_k}\big]$. Using \eqref{logeq:4.1},  \eqref{properties} and  Lemma \ref{leB} we obtain
 \begin{align*}
\mathbb{E}\big[\zeta_k^3 \big|\mathcal{F}_{t_k}\big]
=& \frac{\mathrm{e}^{-3Z_{k}}}{(1+\mathrm{e}^{-Z_{k}})^3}\mathbb{E}\Big[ \Big( a(r_k) \mathrm{e}^{Z_k}\t-\beta(r_k) \Delta+C \Delta^{2(1-\theta )}
+\frac{\sigma^2(r_k) (\Delta B_k)^2}{2}\nn\\
&~~~~~~~~~~~~~~~~~~
- \sigma(r_k) \Delta B_k-\sigma(r_k)\Delta B_k\big(a(r_k) \mathrm{e}^{Z_k}-\beta(r_k)\big)\Delta
+\mathcal{U}_k\Big)^3\Big|\mathcal{F}_{t_k}\Big]\nn\\
\leq& \frac{C\mathrm{e}^{-3Z_{k}}}{\big(1+\mathrm{e}^{-Z_{k}}\big)^3}
\Big(  a^3(r_k)\mathrm{e}^{3Z_k}\t^3+ \t^{\frac{3}{2}}  +\mathbb{E}\big[\mathcal{U}_k^3\big|\mathcal{F}_{t_k}\big]\Big)
\leq  C\t^{\frac{3}{2}}.
\end{align*}
 On the other hand, we can use the same method to derive that
\begin{align*}
\mathbb{E}\big[\zeta_k^3 \big|\mathcal{F}_{t_k}\big]
\geq   -C\t^{\frac{3}{2}}.
\end{align*}
 Thus, both of the above inequalities imply $\mathbb{E}\big[c_0\zeta_k^3\big|\mathcal{F}_{t_k}\big]\leq o(\t)$ for any constant $c_0$, where $c_j$ represents the coefficient of $\zeta_k^j$ term in polynomial $P_m(\zeta_k)$. We can also show that
\begin{align*}
\mathbb{E}\big[|c_j\zeta_k^{3+j}|\big|\mathcal{F}_{t_k}\big]\leq o(\t)
\end{align*}
for any $j\geq 1$. These imply
\begin{align}\la{Y_8}
\mathbb{E}\big[P_m(\zeta_k)\zeta_k^3 \big|\mathcal{F}_{t_k}\big]\leq o(\t).
\end{align}
Combining \eqref{log}, \eqref{Y_4}, \eqref{Y_5} and \eqref{Y_8},   we obtain that
\begin{align*}
 \E \Big[\big(1+\mathrm{e}^{-\bar{Z}_{k+1}}\big)^{p}  \big|\mathcal{F}_{t_k}\Big]
\leq \big(1+\mathrm{e}^{-Z_{k}}\big)^{p}\big(1+C\t\big)
\end{align*}
for any integer $0\leq k\leq \lfloor T/\Delta\rfloor$.  Obviously,
 \begin{align}
 \E \Big[\big(1+\mathrm{e}^{-\bar{Z}_{k+1}}\big)^{p}  \Big]
\leq \big(1+C\t\big)\E \Big[\big(1+\mathrm{e}^{-Z_{k}}\big)^{p}\Big].
\end{align}
Define
$\Omega_k=\big\{\bar{Z}_{k}> \log(K\t^{-\theta})\big\}$.  Using the Chebyshev  inequality, we can  see that
\begin{align*}
 \E \Big[\big(1+\mathrm{e}^{-Z_{k}}\big)^{p}  \Big]=& \E \Big[\big(1+\mathrm{e}^{-Z_{k}}\big)^{p}I_{\Omega^c_k}  \Big]+\E \Big[\big(1+\mathrm{e}^{-Z_{k}}\big)^{p}I_{\Omega_k}  \Big]\nn\\
=& \E \Big[\big(1+\mathrm{e}^{-\bar{Z}_{k}}\big)^{p}I_{\Omega^c_k}  \Big]+\E \Big[\big(1+K^{-1}\t^{\theta}\big)^{p}I_{\Omega_k}  \Big]\nn\\
\leq & \E \Big[\big(1+\mathrm{e}^{-\bar{Z}_{k}}\big)^{p}  \Big]+2^p\mathbb{P}\big\{\bar{Z}_{k}> \log(K\t^{-\theta})\big\}\nn\\
\leq &\big(1+C\t\big)\E \Big[\big(1+\mathrm{e}^{-Z_{k-1}}\big)^{p}\Big]
+2^p\frac{\mathbb{E}\mathrm{e}^{\bar{Z}_{k}/\theta}}{K^{1/\theta}\t^{-1}}.
\end{align*}
 It follows from the result  of  Lemma \ref{log_le3.1} that
\begin{align*}
 \E \Big[\big(1+\mathrm{e}^{-Z_{k}}\big)^{p}  \Big]
\leq &\big(1+C\t\big)\E \Big[\big(1+\mathrm{e}^{-Z_{k-1}}\big)^{p}\Big]
+C_T\t\nn\\
\leq &\big(1+C\t\big)^k \big(1+x^{-1}_{0}\big)^{p}
+C_T\t\sum_{j=0}^{k-1}\big(1+C\t\big)^j\nn\\
\leq &\mathrm{e}^{Ck\t} \big(1+x^{-1}_{0}\big)^{p}
+C_Tk\t \mathrm{e}^{Ck\t} \leq e^{C_T}\Big[\big(1+x^{-1}_{0}\big)^{p}+C_T\Big].
\end{align*}
The proof is complete.
\eproof

\begin{remark}
By   virtue of Lemmas \ref{log_le3.1} and \ref{inv_Z}, for any $q\in \mathbb{R}$,  then the truncated EM scheme defined by \eqref{log_eq5*} has  the following exponential integrability property
\begin{align*}
\sup_{\Delta\in (0, 1)}\sup_{0\leq k\leq \lfloor T/\Delta\rfloor}\mathbb{E}\big[\mathrm{e}^{q Z_k}\big]
\leq  C_T,~~\forall~T>0.
\end{align*}
 Obviously, we also have the following property
\begin{align}\la{num_inv}
\sup_{\Delta\in (0, 1)}\sup_{0\leq k\leq \lfloor T/\Delta\rfloor}\mathbb{E}\big[|Z_k|^p\big]\leq \sup_{\Delta\in (0, 1)}\sup_{0\leq k\leq \lfloor T/\Delta\rfloor}\mathbb{E}\big[\mathrm{e}^{p |Z_k|}\big]
\leq C_T 
\end{align}
for any $p>0$ and $T>0$.
\end{remark}

In order to show that the numerical scheme defined by  \eqref{TEM_1}  perform the dynamical behaviors of exact solutions perfectly,
 we further require the chosen $\theta \in (0, 1/2)$ if $\check{a}>0$.
 By
\eqref{logeq:4.1}, for any $p>0$,
\begin{align*}
 \mathrm{e}^{Z_{k}}\leq \mathrm{e}^{\bar{Z}_{k}}, ~~~ Z^p_{k}=\big[\bar{Z}_{k}\wedge \log(K\Delta^{-\theta } )\big]^{p}\leq \bar{Z}^{p}_{k},~~~ X_k=\mathrm{e}^{Z_{k}}\leq K\Delta^{-\theta }.
\end{align*}
Moreover, to study the
rate of convergence of numerical solutions $\{X_k\}_{k\geq 0}$, we also need to study  the EM method to \eqref{log_eq4}, which are defined as follows:
  For any given stepsize $\Delta\in (0, 1]$,
\begin{align}\la{log_eq5}
\left\{
\begin{array}{ll}
Y_0=\log x_0,&\\
Y_{k+1}=Y_{k}+\big(\beta(r_k)-a(r_k) \mathrm{e}^{Y_k}\big)\Delta +\sigma(r_k) \Delta B_k,&
\end{array}
\right.
\end{align}
for any integer $k\geq 0$, where  $\t B_k=B(t_{k+1})-B(t_{k})$. Transforming back, i.e.
\begin{align}\la{TEM2}
X^{E}_k=\mathrm{e}^{Y_k},~~~~~~k=0,1,\ldots,
\end{align}
gives a strictly positive approximation of the original  SDS \eqref{logistic1}.
 In addition, we also  need the following lemma, the proof of which can be found in Appendix A.
\begin{lemma}\la{log_th2}
The EM  method  defined by \eqref{log_eq5} has the property that
\begin{align}\la{logeq2.14}
\sup_{0\leq k\leq \lfloor T/\Delta\rfloor}\mathbb{E}\Big[  |Y_{k}-y(t_{k})|^{2}\Big]
\leq& C_T\Delta
\end{align}
for any $T>0$ and $\Delta\in (0, 1]$.
\end{lemma}

Define
$
\tau_{L}=\inf\{k\Delta\geq 0: x(t_k)\geq \mathrm{e}^{L}\},~ \rho_{\Delta}=\inf\big\{k\Delta\geq 0: \bar{Z}_{\Delta}(t_k)\geq \log(K\Delta^{-\theta } )\big\}.
$
By  \eqref{log_eq3}  we have
\begin{align*}
\mathrm{e}^{pL}\mathbb{P}\{\tau_{L}\leq T\}\leq \! \mathbb{E}\Big[x^p(T\wedge \tau_{L})I_{\{\tau_{L}\leq T\}}\Big]\!+\mathbb{E}\Big[x^p(T\wedge \tau_{L})I_{\{\tau_{L}> T\}}\Big]=\mathbb{E}x^p(T\wedge \tau_{L})\leq C_T.
\end{align*}
It is easy to see that
\begin{align} \la{logeq3.5}
\mathbb{P}\{\tau_{L}\leq T\}\leq  C_T \mathrm{e}^{-pL},~~~~~\forall~p>0,
\end{align}
  where $C_T$  is  a positive constant independent of $L$.   Moreover, by  virtue of Lemma \ref{log_le3.1}, we have
\begin{align*}
 \E\Big[ \exp{\big(p \log(K\Delta^{-\theta } )  \big)}I_{\{\rho_{\Delta}  \leq  T \}}\Big]
 \leq& \E\Big[ \exp{\big(p\bar{Z}_{\Delta}( T \wedge \rho_{\Delta})  \big)}\Big] =\mathbb{E}\Big[\exp{\big(p\bar{Z}_{  \lfloor\frac{T\wedge  \rho_{\Delta}}{\t}\rfloor}\big)}\Big]\nn\\
 \leq&\sup_{  0\leq k\leq \lfloor T/\t \rfloor}\mathbb{E}\big[\exp{\big(p\bar{Z}_{k}\big)}\big]
\leq C_T,
\end{align*}
implies that
\begin{align} \la{logeq3.6}
 \mathbb{P} {\big\{\rho_{\t} \leq   T \big\}}\leq    C_T K^{-p}\Delta^{\theta p},~~~~\forall~p>0.
\end{align}

\begin{theorem}\la{logTh3.6}
The truncated EM scheme defined by \eqref{log_eq5*} has the property that
\begin{align*}
\sup_{0\leq k\leq\lfloor {T}/{\Delta}\rfloor}\mathbb{E}\Big [  |Z_{k}-y(t_{k})|^{2}\Big]
\leq  C_T\Delta
\end{align*}
for any $\Delta\in (0, 1)$ and $T>0$.
\end{theorem}
{\bf  Proof.}~~Define $\bar{\theta}_{ \t}=\tau_{\log(K\Delta^{-\theta } )} \wedge \rho_{\t}$, $  \Omega_1:= \{\omega: \bar{\theta}_{\t} > T\}$, $\bar{u}_k=Z_k-y(t_{k}),$ for  any $k\Delta\in[0, T]$,
where  $\tau_L$ and $\rho_{\t}$  are defined by \eqref{logeq3.5}, and \eqref{logeq3.6},  respectively.
For any $\bar{p}>2$, using the Young inequality we obtain that
\begin{align}\la{logeq3.7}
\E|\bar{ u}_{k}|^2 =& \E\lf(|\bar{ u}_{k}|^2 I_{\Omega_1}\rt)
+ \E\lf(|\bar{ u}_{k}|^2 I_{\Omega_1^c}\rt) \nonumber \\
 \leq &\E\lf(|\bar{ u}_{k}|^2 I_{\Omega_1}\rt)
+ \frac{2\t}{\bar{p}}\E\lf(|\bar{ u}_{k}|^{\bar{p}}  \rt)+ \frac{\bar{p}-2}{\bar{p}\t^{2/(\bar{p}-2)}} \PP(\Omega_1^c).
\end{align}
It follows from the results of   \eqref{integrability}  and \eqref{num_inv} that
\begin{align}\la{logeq3.8}
 \frac{2\t}{\bar{p}}\E\lf(|\bar{ u}_{k}|^{\bar{p}}  \rt) \leq  C_T\t.
\end{align}
 It follows from  \eqref{logeq3.5}, and \eqref{logeq3.6} that
\begin{align}\la{logeq3.9}
 \frac{\bar{p}-2}{\bar{p}\t^{2/(\bar{p}-2)}} \PP(\Omega_1^c)  \leq &  \frac{\bar{p}-q}{\bar{p}\t^{2/(\bar{p}-2)}}  \Big( \PP {\{\tau_{\log(K\Delta^{-\theta } )} \leq T\}} +\PP {\{\rho_{  \t} \leq T\}}  \Big)\nonumber\\
 \leq & \frac{2(\bar{p}-2)}{\bar{p}\t^{2/(\bar{p}-2)}}   \frac{C_T}{ K^{\delta}\Delta^{-\delta\theta }}
 =\frac{2C_T(\bar{p}-2)}{K^{\delta}\bar{p}}
 \Delta^{ \delta\theta -\frac{2}{\bar{p}-2}}   \leq  C_T\t,
\end{align}
where $ \delta\geq  \bar{p}/ \theta (\bar{p}-2)$.
Inserting \eqref{logeq3.8}, \eqref{logeq3.9}  and \eqref{logeq2.14} into \eqref{logeq3.7} yields
\begin{align*}
\E|\bar{ u}_{k}|^2
 \leq\E\lf(|\bar{ u}_{k}|^2 I_{\Omega_1}\rt)
+  C_T\Delta
\leq&
\sup_{0\leq k\leq\lfloor \frac{T\wedge \rho_{\Delta}}{\Delta}\rfloor}\mathbb{E}\Big[  |Z_{k}-y(t_{k})|^{2}\Big]+  C_T\Delta\nn\\
=&
\sup_{0\leq k\leq\lfloor \frac{T\wedge \rho_{\Delta}}{\Delta}\rfloor}\mathbb{E}\Big[  |Y_{k}-y(t_{k})|^{2}\Big]+  C_T\Delta
\leq  C_T\Delta .
\end{align*}
The proof is complete.  \eproof
As a consequence we also obtain the same convergence order for the approximation of the original SDS \eqref{logistic1} by $X_k:=\mathrm{e}^{Z_k}$.
\begin{theorem}\la{log_th3.2}
For any $0<p<1$, there exists a constant $C_T>0$ such that
\begin{align*}
 \sup_{0\leq k\leq \lfloor T/\Delta\rfloor}\mathbb{E}\Big[ |x(t_k)-X_{k}|^{2p}\Big]
\leq& C_T\Delta^{p}
\end{align*}
for any $\Delta\in (0, 1)$ and $T>0$.
\end{theorem}
{\bf  Proof.}~~   Using H\"{o}lder's inequality,   we have
\begin{align*}
 \mathbb{E}\bigg[  |x(t_k)-\mathrm{e}^{Z_k}|^{2p}\bigg]
=&\mathbb{E}\bigg[ |\mathrm{e}^{y(t_k)}-\mathrm{e}^{Z_k}|^{2p}\bigg]
\leq \mathbb{E}\bigg[ \big(\mathrm{e}^{2p y(t_k)}+ \mathrm{e}^{2pZ_k}\big) |y(t_k)-Z_k|^{2p}\bigg]\nn\\
\leq&\bigg[\mathbb{E} \big(\mathrm{e}^{2p y(t_k)}+ \mathrm{e}^{2pZ_k}\big)^{1-p} \bigg]^{1/(1-p)} \bigg[\mathbb{E}  |y(t_k)-Z_k|^{2}\bigg]^{p}.
\end{align*}
Thus, by applying Theorem \ref{logTh3.6}, Lemma \ref{log_le3.1} and \eqref{log_eq3}, we infer that
\begin{align*}
\mathbb{E}\bigg[  |x(t_k)-\mathrm{e}^{Z_k}|^{2p}\bigg] \leq C_T\Delta^{p}
\end{align*}
for any $\Delta\in (0, 1)$. The proof is complete. \eproof

\begin{remark}
As a consequence we also obtain the same convergence order for the approximation of the original SDS \eqref{logistic1} by $X^E_k=\mathrm{e}^{Y_k}$, where $Y_k$ is defined by \eqref{log_eq5}.  Again, for numerical solutions $\{X^E_k\}_{k\geq 0}$, by applying Lemmas \ref{log_th2} and \ref{log_le4}, for any  $0<p<1$ we have
\begin{align*}
 \sup_{k=0,\ldots, \lfloor T/\Delta\rfloor}\mathbb{E}\Big[ |x(t_k)-X^{E}_k|^{2p}\Big]
\leq& C_T\Delta^{p}
\end{align*}
for any $\Delta\in (0, 1]$ and $T>0$.
\end{remark}

\section{Stochastic permanence}\label{strong2}
In this section, we focus on the stochastic permanence and extinction.   Firstly, we establish the criterion on the  stochastic permanence and extinction of the exact solution  of SDS \eqref{logistic1}. Then we show that the numerical scheme defined by  \eqref{TEM_1} keep this property very well.
From this  section  we always assume   $r(t)$ is {\it irreducible},
 namely, the following linear equation
 \begin{eqnarray}\label{eq:a1.2}
\pi \Gamma=0,\ \ \ \ \ \ \ \  \ \sum_{i=1}^{m}\pi_i=1,
\end{eqnarray}
has a unique solution $\pi=(\pi_1, \dots, \pi_m)\in {\mathbb R}^{1\times m}$  satisfying $\pi_i>0$ for each $i\in \mathbb{S}$. This solution is termed a stationary distribution. Then the rank of $\Gamma$ is $m-1$. It follows that null space of $\Gamma$ is one dimensional spanned by $\I_m:=(1,   \cdots, 1)^T\in\mathbb{R}^m$.
Consider the linear equation
\begin{eqnarray}\label{eq:a2.1}
\Gamma \nu=\xi,
\end{eqnarray}
where $\nu$ and $\xi\in \mathbb{R}^m$.

\begin{lemma}[\!\!{\cite{yz09}}]\label{Le:Yin} The following assertions hold.
\begin{itemize}
\item[$\mathrm{(1)}$]Equation (\ref{eq:a2.1}) has a solution if and only if $\pi \xi =0$.

\item[$\mathrm{(2)}$] Suppose that $\nu_1$ and $\nu_2$ are two solutions of (\ref{eq:a2.1}). Then $\nu_1-\nu_2=\kappa \I_m$ for some $\kappa\in \mathbb{R}$.

\item[$\mathrm{(3)}$] Any solution of (\ref{eq:a2.1}) can be written as $\nu=\kappa \I_m+h_0$, where $\kappa\in \mathbb{R}$ is an arbitrary
constant, and $h_0\in \mathbb{R}^m$ is the unique solution of (\ref{eq:a2.1}) satisfying $\pi h_0=0$.
\end{itemize}
\end{lemma}

For the   definitions of stochastic permanence and its relatives (see e.g., Li and Yin \cite[Definitions 2.1-2.3]{Li2016JMAA}).
 We begin with the following lemmas and make use of it to obtain the stochastically
ultimate upper boundedness of SDS \eqref{logistic1}.
\begin{lemma}[\!\!\cite{Wangr17}]\la{L:2}
If $
\pi a:= \sum_{i\in \mathbb{S}}\pi_i a_i>0
$ hold,  for any $\eta\in (0,1)$ sufficiently small,
 the solution  $x(t)$ of SDS \eqref{logistic1} with any initial value
  $x_0\in \mathbb{R}_+$, $\ell\in \mathbb{S}$ has the property that
\begin{align*}
\limsup_{t\rightarrow \infty}\mathbb{E}\big[\log^\eta \big(x(t)\vee 1\big) \big]\leq C.
\end{align*}
  \end{lemma}
In the special case where $\hat{a}>0$,  we cite the following lemma from literature.
\begin{lemma}[\!\!\cite{Li2011}]\la{log_le_le_le2}
If $\hat{a}>0$, for any $p>0$,
 the solution  $x(t)$ of SDS \eqref{logistic1} with any initial value
  $x_0\in \mathbb{R}_+$, $\ell\in \mathbb{S}$ has the property that
\begin{align*}
\limsup_{t\rightarrow \infty}\mathbb{E}\big[x^p(t)\big]
\leq  C.
\end{align*}
\end{lemma}

Next we continue to consider the case $\pi a =0$. As we know, either $\pi a>0$ or $\pi a =0$  because of each $a(i) \equiv0$. The $\pi a=0$ implies that  $\check{a}=0$, then SDS \eqref{logistic1} degenerates into
\be\label{logeq:4.19}
\mathrm{d}x(t)=  b(r(t)) x(t)\mathrm{d}t + \sigma(r(t))x(t)\mathrm{d}B(t).
\ee
We give the limit of the moment of linear SDS \eqref{logeq:4.19} for small $p$, which is stronger
than the stochastically ultimate upper boundedness.
\begin{lemma}[\!\!\cite{Li2016JMAA}]\la{log:Le5.4}
If $
\pi \beta:= \sum_{i\in \mathbb{S}}\pi_i \beta_i<0
$ hold, then for any $\rho\in(0,1)$ sufficiently small,
the solution  $x(t)$ of SDS \eqref{logeq:4.19} has the property that
\begin{align*}
\lim_{t\rightarrow \infty}\mathbb{E}\big[ x^{\rho}(t) \big]=0.
\end{align*}
  \end{lemma}
Moreover, the following result is a direct consequence of \eqref{log_eq4}.
\begin{lemma}[\!\!\cite{Li2016JMAA}]\la{log:Le5.5}
The solution  $x(t)$ of SDS \eqref{logistic1} with any initial value $x_0\in\mathbb{R}_+$, $\ell\in \mathbb{S}$ satisfies
\begin{align*}
\limsup_{t\rightarrow \infty}\frac{y(t)}{t}=\limsup_{t\rightarrow \infty}\frac{\log x(t)}{t}\leq\pi \beta~~~~a.s.
\end{align*}
In particular, when $\pi a=0$,
\begin{align*}
\limsup_{t\rightarrow \infty}\frac{y(t)}{t}=\limsup_{t\rightarrow \infty}\frac{\log x(t)}{t}=\pi \beta~~~~a.s.
\end{align*}
  \end{lemma}

Now we look for the ultimate lower boundary of the moment of solutions in order
to obtain the stochastic permanence.
\begin{lemma}[\!\!\cite{Li2016JMAA}]\la{log:Le5.8}
If $
\pi \beta:= \sum_{i\in \mathbb{S}}\pi_i \beta_i>0
$ hold, for any $\vartheta\in(0,1)$ sufficiently small,
the solution  $x(t)$ of SDS \eqref{logistic1} has the property that
\begin{align}\la{logyhf1}
\limsup_{t\rightarrow \infty}\mathbb{E}\big[ x^{-\vartheta}(t) \big]\leq C.
\end{align}
  \end{lemma}

\begin{theorem}[\!\!\cite{Li2016JMAA}]\la{log_permanence}
If $\pi a>0$ and $\pi \beta>0$ hold, SDS \eqref{logistic1} is stochastically permanent.
\end{theorem}

\begin{theorem}[\!\!\cite{Li2016JMAA}]\label{log_Th*5.3}
Suppose that  $\pi \beta\neq0$. Then,
 \begin{itemize}
\item  the solution of \eqref{logistic1} are
stochastically permanent if and only if $ \pi a>0,~\pi \beta>0$;
\item the solution of  \eqref{logistic1} are
almost surely extinctive if and only if $\pi \beta<0$;
    \item   almost all paths of \eqref{logistic1} increase at an exponential rate if and only if $\pi a=0, ~\pi \beta>0$.
     \end{itemize}
\end{theorem}

\subsection{Stochastic permanence of numerical solution}
In order to approximate the stochastic permanence  of SDS \eqref{logistic1} we need to construct the appropriate scheme such that the numerical solutions must be both stochastically ultimately upper bounded and lower bounded. In this subsection, we first give the definitions of the stochastic permanence and the stochastically ultimate boundedness of numerical solutions to the SDS \eqref{logistic1}.
\begin{defn}
A time-discretization $\{X_k\}$, with stepsize $\Delta\in(0,1]$, of the solution to the SDS \eqref{logistic1} is said to be  stochastically ultimately
upper bounded, if for any $\nu\in(0,1)$, there exist a positive constant $\chi=\chi(\nu)$ such
that for any initial value  $x_0\in \mathbb{R}_+, \ell\in \mathbb{S}$ satisfies
$$
\limsup_{k\rightarrow +\infty}\mathbb{P}\big\{X_k>\chi\big\}<\nu.
$$
\end{defn}

\begin{defn}
A time-discretization $\{X_k\}$, with stepsize $\Delta\in(0,1]$, of the solution to the SDS \eqref{logistic1} is said to be stochastically ultimately
lower bounded, if for any $\nu\in(0,1)$, there exist a positive constant $\chi=\chi(\nu)$ such
that for any initial value  $x_0\in \mathbb{R}_+, \ell\in \mathbb{S}$   satisfies
$$
\limsup_{k\rightarrow +\infty}\mathbb{P}\big\{X_k<\chi\big\}<\nu.
$$
\end{defn}

\begin{defn}
A time-discretization $\{X_k\}$, with stepsize $\Delta\in(0,1]$, of the solution to the SDS \eqref{logistic1}   is said to be stochastically permanent if its time-discretization solutions
are both stochastically ultimately upper bounded and stochastically ultimately
lower bounded.
\end{defn}

For convenience, denote by $\mathcal{G}_{t_k}$ the $\sigma$-algebra generated by $\{\mathcal{F}_{t_k}, r_{k+1}\}$. Obviously, $\mathcal{F}_{t_k}\subseteq\mathcal{G}_{t_k}$.
 We begin with a criterion on asymptotic upper boundedness of the moment, and make use of it to obtain the stochastically ultimate upper boundedness of the numerical solutions.
\begin{lemma}\la{L:3}
Under the condition  of Lemma \ref{L:2},
there is a constant $\Delta^{*}_1\in (0, 1)$ such that the scheme
 \eqref{TEM_1} has the property that
\begin{align*}
\sup_{k \geq0}\mathbb { E } \left[ \log^{\eta}\big( X_{k}\vee 1\big) \right]=\sup_{k \geq0}\mathbb { E } \left[ \big(Z_{k}\vee 0\big)^{\eta} \right]\leq\sup_{k \geq0}\mathbb { E } \left[ \big( \bar{Z}_{k}\vee 0 \big)^{\eta}\right]
 \leq  C
\end{align*}
for any $\Delta\in (0, \Delta^{*}_1]$, where $\eta$ is defined in Lemma \ref{L:2}.
\end{lemma}
{\bf  Proof.}~~Note that
 \begin{eqnarray*}
\pi[-a+(\pi a)\I_m]=0,\ \ \ \ \ \sum_{i=1}^{m}\pi_i=1.
\end{eqnarray*}
It follows from Lemma \ref{Le:Yin} (1) that the equation
 \begin{eqnarray*}
\dis \Gamma c=-a+(\pi a)\I_m
\end{eqnarray*}
has a solution $c=(c_1,   \cdots, c_m)^T\in \mathbb{R}^m$. Thus we have
 \begin{eqnarray}\label{eq:2.7}
a(i)+\sum_{j=1}^{m}\gamma_{i j}c_j=\pi a>0,\ \ \ i\in \mathbb{S}.
\end{eqnarray}
Using the well-known  Taylor formula we get
\begin{align}\la{log:eq:5.8}
\mathrm{e}^{\bar{Z}_{k+1}}= \mathrm{e}^{Z_{k}}+\mathrm{e}^{Z_{k}}(\bar{Z}_{k+1}-Z_{k})+\frac{1}{2}\mathrm{e}^{Z_{k}}(\bar{Z}_{k+1}-Z_{k})^2
+\frac{1}{6}\mathrm{e}^{\bar{\xi}_k}(\bar{Z}_{k+1}-Z_{k})^3,
\end{align}
where $\bar{\xi}_k\in (\bar{Z}_{k+1}\wedge Z_{k}, \bar{Z}_{k+1}\vee Z_{k})$.   For $\omega\in\{\bar{Z}_{k+1}< Z_{k}\}$, we have
 \begin{align}\la{xuyao1}
\mathrm{e}^{\bar{\xi}_k}(\bar{Z}_{k+1}-Z_{k})^3\leq 0.
 \end{align}
 On the other hand, for  $\omega\in\{\bar{Z}_{k+1}\geq Z_{k}\}$, we have \begin{align}\la{T_w2}
&  \frac{1}{6}\mathrm{e}^{\bar{\xi}_k}(\bar{Z}_{k+1}-Z_{k})^3
\leq \frac{1}{6}\mathrm{e}^{\bar{Z}_{k+1}}(\bar{Z}_{k+1}-Z_{k})^3
= \frac{1}{6}\mathrm{e}^{Z_{k}}\mathrm{e}^{\bar{Z}_{k+1}-Z_{k}}(\bar{Z}_{k+1}
-Z_{k})^3\nn\\
=&\frac{1}{6}\mathrm{e}^{Z_{k}}\exp{\Big(\big( \beta(r_k)-a(r_k) \mathrm{e}^{Z_k}\big)\Delta +\sigma(r_k) \Delta B_k\Big)}\Big(\big(\beta(r_k)-a(r_k) \mathrm{e}^{Z_k}\big)\Delta +\sigma(r_k) \Delta B_k\Big)^3\nn\\
\leq&\frac{2}{3}\mathrm{e}^{Z_{k}} \mathrm{e}^{\check{\beta}  \Delta} \big(|\check{\beta}| \vee |\breve{\sigma}| \big)^3\exp{\Big( |\breve{\sigma}| |\Delta B_k|\Big)}\Big(  \Delta^3 +  |\Delta B_k|^3\Big)=:\mathrm{e}^{Z_{k}}\bar{\mathcal{U}}_k>0,
\end{align}
where
$$
\bar{\mathcal{U}}_k:=\frac{2}{3} \mathrm{e}^{\check{\beta}  \Delta} \big(|\check{\beta}| \vee |\breve{\sigma}| \big)^3\exp{\Big( |\breve{\sigma}| |\Delta B_k|\Big)}\Big(  \Delta^3 +  |\Delta B_k|^3\Big).
$$
Therefore, we derive from \eqref{log:eq:5.8}-\eqref{T_w2} that for any integer $k\geq 0$,
\begin{align}\la{log:eq:5.9}
\!\!\mathrm{e}^{\bar{Z}_{k+1}}-\mathrm{e}^{Z_{k}}\leq& \mathrm{e}^{Z_{k}}\Big[ \bar{Z}_{k+1}-Z_{k} +\frac{1}{2} (\bar{Z}_{k+1}-Z_{k})^2
+\bar{\mathcal{U}}_k\Big]\nn\\
=&\mathrm{e}^{Z_{k}}\Big[\big(\beta(r_k)-a(r_k) \mathrm{e}^{Z_k}\big)\Delta +\sigma(r_k) \Delta B_k+\frac{\sigma^2(r_k)}{2} \big(\Delta B_k\big)^2\nn\\
&~~ +\sigma(r_k) \Delta B_k\big(\beta(r_k)-a(r_k) \mathrm{e}^{Z_k}\big) \Delta
+\frac{1}{2}\big(\beta(r_k)-a(r_k) \mathrm{e}^{Z_k}\big)^2\Delta^2+\bar{\mathcal{U}}_k\Big].
\end{align}
Choose a constant $0<\eta_0\leq 1$ such that for each $0<\eta\leq \eta_0$,   $$\xi^{\eta,c}_i:=1-c_i \eta>0,  \ i=1,\cdots,m.$$
By the Markov property (see, \cite[Lemma 3.2]{li2018jde} for more details), we derive that
\begin{align}\la{logeq:5.9}
\E\big[\xi^{\eta,c}_{r_{k+1}}\big|\mathcal{F}_{t_k}\big]=\xi^{\eta,c}_{r_{k}}+\sum_{j\in \mathbb{S}}\xi^{\eta,c}_{j}\left(\gamma_{r_k j}\t+o(\t)\right).
\end{align}
Moreover, one observes
\begin{align*}
 \log\big(1+\mathrm{e}^{\bar{Z}_{k+1}}\big)
\leq   \log\big(1+\mathrm{e}^{Z_{k}}\big)
+\frac{\mathrm{e}^{\bar{Z}_{k+1}}-\mathrm{e}^{Z_{k}}}{1+\mathrm{e}^{Z_{k}}}.
\end{align*}
Then using the above inequality and \eqref{log:eq:5.9},   we have
 \begin{align}\la{logeq:5.7}
\big[1+\log\big(1+\mathrm{e}^{\bar{Z}_{k+1}}\big)\big]^{\eta}
\leq &\big[1+\log\big(1+\mathrm{e}^{Z_{k}}\big)\big]^{\eta}\big(1+\varsigma_k\big)^{\eta},
\end{align}
where
 \begin{align*}
\varsigma_k=&\big[1+\log\big(1+\mathrm{e}^{Z_{k}}\big)\big]^{-1} \frac{\mathrm{e}^{Z_{k}}}{1+\mathrm{e}^{Z_{k}}}\Big[\big(\beta(r_k)-a(r_k) \mathrm{e}^{Z_k}\big)\Delta +\sigma(r_k) \Delta B_k+\frac{\sigma^2(r_k)}{2} \big(\Delta B_k\big)^2\nn\\
&~~~~~~~~~~~~~~~~~~~~~~~~+\sigma(r_k) \Delta B_k\big(\beta(r_k)-a(r_k) \mathrm{e}^{Z_k}\big) \Delta
+\frac{1}{2}\big(\beta(r_k)-a(r_k) \mathrm{e}^{Z_k}\big)^2\Delta^2+\bar{\mathcal{U}}_k\Big],
\end{align*}
and we can see that $\varsigma_k>-1$.   For any $0<\eta\leq 1$,  by  virtue of \cite[Lemma 3.3]{li2018jde},
we derive from \eqref{logeq:5.7} that
\begin{align}\la{logeq:*5.9}
&\E \Big[\big[1+\log\big(1+\mathrm{e}^{\bar{Z}_{k+1}}\big)\big]^{\eta}\xi^{\eta,c}_{r_{k+1}} \big|\mathcal{F}_{t_k}\Big]\nn\\
\leq&\big[1+\log\big(1+\mathrm{e}^{Z_{k}}\big)\big]^{\eta} \bigg\{\E\big[\xi^{\eta,c}_{r_{k+1}}\big|\mathcal{F}_{t_k}\big]+\eta \E\big[\varsigma_k\xi^{\eta,c}_{r_{k+1}}\big|\mathcal{F}_{t_k}\big] + \frac{\eta(\eta-1)}{2} \E\big[\varsigma_k^2\xi^{\eta,c}_{r_{k+1}}\big|\mathcal{F}_{t_k}\big]\nn\\
 &~~~~~~~~~~~~~~~~~~~~~~ ~~~~~+ \frac{\eta(\eta-1)(\eta-2)}{6} \E\big[\varsigma_k^3\xi^{\eta,c}_{r_{k+1}}
 \big|\mathcal{F}_{t_k}\big]\bigg\}\nn\\
 \leq&\big[1+\log\big(1+\mathrm{e}^{Z_{k}}\big)\big]^{\eta} \bigg\{\E\big[\xi^{\eta,c}_{r_{k+1}}\big|\mathcal{F}_{t_k}\big]+\eta \E\big[\varsigma_k\xi^{\eta,c}_{r_{k+1}}\big|\mathcal{F}_{t_k}\big] +  \E\big[\varsigma_k^3\xi^{\eta,c}_{r_{k+1}}
 \big|\mathcal{F}_{t_k}\big]\bigg\}.
\end{align}
 Using the properties
\begin{align}\la{logeq:5.10}
\mathbb{E} \big(|\Delta B_{k} |^{2j}\big|\mathcal{G}_{t_k}\big)=(2j-1)!!\Delta^{j},~~~
\mathbb{E} \big((\Delta B_{k} )^{2j-1}\big|\mathcal{G}_{t_k}\big)=0,~~j=1,2,\ldots
\end{align}
 we deduce that
\begin{align*}
\E\big[\varsigma_k \big|\mathcal{G}_{t_k}\big]=&\big[1+\log\big(1+\mathrm{e}^{Z_{k}}\big)\big]^{-1} \frac{\mathrm{e}^{Z_{k}}}{1+\mathrm{e}^{Z_{k}}}\Big[\big(\beta(r_k)-a(r_k) \mathrm{e}^{Z_k}\big)\Delta + \frac{\sigma^2(r_k)}{2} \Delta  \nn\\
&~~~~~~~~~~~~~~~~~~~~~~~~~~~~~~~~~~~~
+\frac{1}{2}\big(\beta(r_k)-a(r_k) \mathrm{e}^{Z_k}\big)^2\Delta^2+\E\big(\bar{\mathcal{U}}_k\big|\mathcal{G}_{t_k}\big)\Big]\nn\\
=&\big[1+\log\big(1+\mathrm{e}^{Z_{k}}\big)\big]^{-1} \frac{\mathrm{e}^{Z_{k}}}{1+\mathrm{e}^{Z_{k}}}\Big[\big(b(r_k)-a(r_k) \mathrm{e}^{Z_k}\big)\Delta    \nn\\
&~~~~~~~~~~~~~~~~~~~~~~~~~~~~~~~~~~~~
+ \big( |\breve{\beta}|^2+\check{a}^2  K^2\Delta^{-2\theta }\big) \Delta^2+\E\big(\bar{\mathcal{U}}_k\big|\mathcal{G}_{t_k}\big)\Big].
\end{align*}
By virtue of Lemma \ref{leB}, we obtain
 $
\E\big(\bar{\mathcal{U}}_k\big|\mathcal{G}_{t_k}\big)   \leq C\Delta^{\frac{3}{2}}.
$
  Then
\begin{align*}
\E\big[\varsigma_k \big|\mathcal{G}_{t_k}\big]\leq&\big[1+\log\big(1+\mathrm{e}^{Z_{k}}\big)\big]^{-1} \frac{\mathrm{e}^{Z_{k}}}{1+\mathrm{e}^{Z_{k}}}\Big[ \big(\check{b} -a(r_k) \mathrm{e}^{Z_k}\big)\Delta
+ C\Delta^{2(1-\theta )}+C\Delta^{\frac{3}{2}}\Big].
\end{align*}
Making use of the above inequality and \eqref{logeq:5.9} yields
\begin{align}\la{logeq:5.14}
& \E\big[\varsigma_k\xi^{\eta,c}_{r_{k+1}}\big|\mathcal{F}_{t_k}\big]
 =   \E\big[\xi^{\eta,c}_{r_{k+1}}\E\big(\varsigma_k\big|\mathcal{G}_{t_k}\big)\big|\mathcal{F}_{t_k}\big]\nn\\
\leq& \big[1+\log\big(1+\mathrm{e}^{Z_{k}}\big)\big]^{-1} \frac{\mathrm{e}^{Z_{k}}}{1+\mathrm{e}^{Z_{k}}}\Big[ \big(\check{b} -a(r_k) \mathrm{e}^{Z_k}\big)\Delta
+ C\Delta^{2(1-\theta )}+C\Delta^{\frac{3}{2}}\Big]
\E\big[\xi^{\eta,c}_{r_{k+1}}\big|\mathcal{F}_{t_k}\big]\nn\\
=& \big[1+\log\big(1+\mathrm{e}^{Z_{k}}\big)\big]^{-1} \frac{\mathrm{e}^{Z_{k}}}{1+\mathrm{e}^{Z_{k}}}\Big[ \big(\check{b} -a(r_k) \mathrm{e}^{Z_k}\big)\Delta\nn\\
&~~~~~~~~~~~~~~~~~~~~~~~~~~~+ C\Delta^{2(1-\theta )}+C\Delta^{\frac{3}{2}}\Big]
 \Big[\xi^{\eta,c}_{r_{k}}+\sum_{j\in \mathbb{S}}\xi^{\eta,c}_{j}\left(\gamma_{r_k j}\t+o(\t)\right)\Big]\nn\\
\leq& \big[1+\log\big(1+\mathrm{e}^{Z_{k}}\big)\big]^{-1} \frac{\mathrm{e}^{Z_{k}}}{1+\mathrm{e}^{Z_{k}}}\Big[ \big(\check{b} -a(r_k) \mathrm{e}^{Z_k}\big)\xi^{\eta,c}_{r_{k}}\Delta + C\Delta^{2(1-\theta )}+C\Delta^{\frac{3}{2}}\Big]
.
\end{align}
On the other hand,  it follows from \eqref{logeq:5.7} that
\begin{align*}
\E\big[\varsigma^3_k \big|\mathcal{G}_{t_k}\big]
=&\big[1+\log\big(1+\mathrm{e}^{Z_{k}}\big)\big]^{-3} \frac{\mathrm{e}^{3Z_{k}}}{ (1+\mathrm{e}^{Z_{k}} )^3}\E\Bigg\{\sigma^3_1(r_k) \big(\Delta B_k\big)^3\big(1+\big(\beta(r_k)-a(r_k) \mathrm{e}^{Z_k}\big) \Delta\big)^3\nn\\
&+3\sigma^2(r_k) \big(\Delta B_k\big)^2\big(1+\big(\beta(r_k)-a(r_k) \mathrm{e}^{Z_k}\big) \Delta\big)^2\Big[\big(\beta(r_k)-a(r_k) \mathrm{e}^{Z_k}\big)\Delta\nn\\
 &~~~~~~~~~~~~~~~~~~~~~~~~~~~+ \frac{\sigma^2(r_k)}{2} \big(\Delta B_k\big)^2 +\frac{1}{2}\big(\beta(r_k)-a(r_k) \mathrm{e}^{Z_k}\big)^2\Delta^2+\bar{\mathcal{U}}_k\Big] \nn\\
 &+3\sigma(r_k)  \Delta B_k \big(1+\big(\beta(r_k)-a(r_k) \mathrm{e}^{Z_k}\big) \Delta\big) \Big[\big(\beta(r_k)-a(r_k) \mathrm{e}^{Z_k}\big)\Delta\nn\\
 &~~~~~~~~~~~~~~~~~~~~~~~~~~~+ \frac{\sigma^2(r_k)}{2} \big(\Delta B_k\big)^2 +\frac{1}{2}\big(\beta(r_k)-a(r_k) \mathrm{e}^{Z_k}\big)^2\Delta^2+\bar{\mathcal{U}}_k\Big]^2\nn\\
&+\Big[\big(\beta(r_k)-a(r_k) \mathrm{e}^{Z_k}\big)\Delta+ \frac{\sigma^2(r_k)}{2} \big(\Delta B_k\big)^2\nn\\
 &~~~~~~~~~~~~~~~~~~~~~~~~~~~ +\frac{1}{2}\big(\beta(r_k)-a(r_k) \mathrm{e}^{Z_k}\big)^2\Delta^2+\bar{\mathcal{U}}_k\Big]^3\Big|\mathcal{G}_{t_k}\Bigg\}.
\end{align*}
Using the properties \eqref{logeq:5.10} and
\begin{align*}
 \E\Big[\big(\Delta B_k\big)\mathrm{e}^{  |\breve{\sigma}||\Delta B_k|} \big( \Delta^{3}+  |\Delta B_k|^3\big)  \Big] =0,~~~~~
 \E\Big[\big(\Delta B_k\big)\mathrm{e}^{  2|\breve{\sigma}||\Delta B_k|} \big( \Delta^{3}+  |\Delta B_k|^3\big)^2 \Big] =0,
 \end{align*}
we deduce that
\begin{align*}
\E\big[\varsigma^3_k \big|\mathcal{G}_{t_k}\big]
=&\big[1+\log\big(1+\mathrm{e}^{Z_{k}}\big)\big]^{-3} \frac{\mathrm{e}^{3Z_{k}}}{ (1+\mathrm{e}^{Z_{k}} )^3}\E\Bigg\{3\sigma^2(r_k) \big(\Delta B_k\big)^2\big(1+\big(\beta(r_k)-a(r_k) \mathrm{e}^{Z_k}\big) \Delta\big)^2 \nn\\
&\times\Big[\big(\beta(r_k)-a(r_k) \mathrm{e}^{Z_k}\big)\Delta + \frac{\sigma^2(r_k)}{2} \big(\Delta B_k\big)^2 +\frac{1}{2}\big(\beta(r_k)-a(r_k) \mathrm{e}^{Z_k}\big)^2\Delta^2+\bar{\mathcal{U}}_k\Big] \nn\\
&+\Big[\big(\beta(r_k)\!-a(r_k) \mathrm{e}^{Z_k}\big)\Delta\!+ \frac{\sigma^2(r_k)}{2} \big(\Delta B_k\big)^2 \!+\frac{1}{2}\big(\beta(r_k)\!-a(r_k) \mathrm{e}^{Z_k}\big)^2\Delta^2\!+\bar{\mathcal{U}}_k\Big]^3\Big|\mathcal{G}_{t_k}\Bigg\}\nn\\
 \leq&\big[1+\log\big(1+\mathrm{e}^{Z_{k}}\big)\big]^{-3}  \E\Bigg\{9|\breve{\sigma}|^2 \big(\Delta B_k\big)^2\big(1+|\breve{\beta}|^2\Delta^2+\check{a}^2K^2\Delta^{2(1-\theta )}\big)  \nn\\
&\times\Big[ |\breve{\beta}| \Delta + \frac{\sigma^2(r_k)}{2} \big(\Delta B_k\big)^2 +  |\breve{\beta}|^2\Delta^2+\check{a}^2  K^2\Delta^{2(1-\theta )} +\bar{\mathcal{U}}_k\Big] \nn\\
&+\Big[|\breve{\beta}|  \Delta+ \frac{\sigma^2(r_k)}{2} \big(\Delta B_k\big)^2 +  |\breve{\beta}|^2\Delta^2+\check{a}^2  K^2\Delta^{2(1-\theta )} +\bar{\mathcal{U}}_k\Big]^3\Big|\mathcal{G}_{t_k}\Bigg\}\nn\\
 \leq&\big[1+\log\big(1+\mathrm{e}^{Z_{k}}\big)\big]^{-3}  \Big(C\Delta^2+9|\breve{\sigma}|^2 \E\big[\big(\Delta B_k\big)^2\bar{\mathcal{U}}_k\big] +10\E\big[\bar{\mathcal{U}}^3_k \big]\Big).
\end{align*}
By virtue of Lemma \ref{leB}, we obtain
\begin{align*}
\E\Big[\big(\Delta B_k\big)^2 \bar{\mathcal{U}}_k \Big|\mathcal{G}_{t_k}
 \Big]\leq&C\Big\{\Delta^{3}\E\Big[\mathrm{e}^{\frac{ |\Delta B_k|^2}{4\Delta}} |\Delta B_k |^2\Big]+ \E\Big[\mathrm{e}^{\frac{ |\Delta B_k|^2}{4\Delta}} |\Delta B_k|^5\Big]\bigg\}\leq C\Delta^{\frac{5}{2}},
\end{align*}
and
\begin{align*}
 \E\Big[ \bar{\mathcal{U}}^3_k\Big|\mathcal{G}_{t_k}
 \Big]
\leq &C\E\Big[ \mathrm{e}^{\frac{ |\Delta B_k|^2}{4\Delta}}\Big(  \Delta^{9}+  |\Delta B_k|^9\Big)\Big]\leq C\Delta^{\frac{9}{2}},
\end{align*}
which implies
\begin{align}\la{logeq:5.16}
& \E\big[\varsigma_k^3\xi^{\eta,c}_{r_{k+1}}\big|\mathcal{F}_{t_k}\big]
 =   \E\big[\xi^{\eta,c}_{r_{k+1}}\E\big(\varsigma_k^3\big|\mathcal{G}_{t_k}\big)\big|\mathcal{F}_{t_k}\big]\nn\\
\leq& C\big[1+\log\big(1+\mathrm{e}^{Z_{k}}\big)\big]^{-3}\t^2
\E\big[\xi^{\eta,c}_{r_{k+1}}\big|\mathcal{F}_{t_k}\big]
 \leq  C\big[1+\log\big(1+\mathrm{e}^{Z_{k}}\big)\big]^{-3} \t^2.
\end{align}
Combining \eqref{logeq:*5.9}, \eqref{logeq:5.14} and \eqref{logeq:5.16},  we derive from \eqref{logeq:5.9}  that
\begin{align*}
&\E \Big[\big[1+\log\big(1+\mathrm{e}^{\bar{Z}_{k+1}}\big)\big]^{\eta}\xi^{\eta,c}_{r_{k+1}} \big|\mathcal{F}_{t_k}\Big]\nn\\
\leq&\big[1+\log\big(1+\mathrm{e}^{Z_{k}}\big)\big]^{\eta} \bigg\{\xi^{\eta,c}_{r_{k}}+\sum_{j\in \mathbb{S}}\xi^{\eta,c}_{j}\left(\gamma_{r_k j}\t+o(\t)\right)+C   \big[1+\log\big(1+\mathrm{e}^{Z_{k}}\big)\big]^{-3}  \Delta^2\nn\\
 &~~~~~~~~~~~~~~~~~  +\eta \big[1+\log\big(1+\mathrm{e}^{Z_{k}}\big)\big]^{-1} \frac{\mathrm{e}^{Z_{k}}}{1+\mathrm{e}^{Z_{k}}}\Big[ \big(\check{b} -a(r_k) \mathrm{e}^{Z_k}\big)\xi^{\eta,c}_{r_{k}}\Delta + C\Delta^{2(1-\theta )}+C\Delta^{\frac{3}{2}}\Big] \bigg\}\nn\\
 \leq& \xi^{\eta,c}_{r_{k}}\big[1+\log\big(1+\mathrm{e}^{Z_{k}}\big)\big]^{\eta} \bigg\{1+\frac{1}{\xi^{\eta,c}_{r_{k}}}\sum_{j\in \mathbb{S}}\xi^{\eta,c}_{j} \gamma_{r_k j}\t +o(\t) \nn\\
 &~~~~~~~~~~~~~~~~~~~~~~ ~~~~~+\eta \big[1+\log\big(1+\mathrm{e}^{Z_{k}}\big)\big]^{-1} \frac{\mathrm{e}^{Z_{k}}}{1+\mathrm{e}^{Z_{k}}}  \big(\check{b} -a(r_k) \mathrm{e}^{Z_k}\big) \Delta  \bigg\}\nn\\
=&\xi^{\eta,c}_{r_{k}}\big[1+\log\big(1+\mathrm{e}^{Z_{k}}\big)\big]^{\eta} \bigg\{1+\frac{1}{\xi^{\eta,c}_{r_{k}}}\sum_{j\in \mathbb{S}}\xi^{\eta,c}_{j} \gamma_{r_k j}\t+o(\t) \nn\\
 &~~~~~~~~~~~~~~~~~~-\eta \big[1+\log\big(1+\mathrm{e}^{Z_{k}}\big)\big]^{-1} \frac{a(r_k)\mathrm{e}^{2Z_{k}}}{1+\mathrm{e}^{Z_{k}}}    \Delta+\eta \big[1+\log\big(1+\mathrm{e}^{Z_{k}}\big)\big]^{-1} \frac{\check{b} \mathrm{e}^{Z_{k}}}{1+\mathrm{e}^{Z_{k}}} \Delta
    \bigg\}\nn\\
\leq&\xi^{\eta,c}_{r_{k}}\big[1+\log\big(1+\mathrm{e}^{Z_{k}}\big)\big]^{\eta} \bigg\{1+\frac{1}{\xi^{\eta,c}_{r_{k}}}\sum_{j\in \mathbb{S}}\xi^{\eta,c}_{j} \gamma_{r_k j}\t+o(\t)
\nn\\
&~~~~~~~~~~~~~~~~~~~~ ~-\eta a(r_k) \big[1+\log\big(1+\mathrm{e}^{Z_{k}}\big)\big]^{-1} \big(1+\mathrm{e}^{Z_{k}}\big)      \Delta +2\eta a(r_k) \big[1+\log\big(1+\mathrm{e}^{Z_{k}}\big)\big]^{-1} \Delta \nn\\
&~~~~~~~~~~~~~~~~~~~~ ~-\eta a(r_k) \big[1+\log\big(1+\mathrm{e}^{Z_{k}}\big)\big]^{-1} \frac{1}{1+\mathrm{e}^{Z_{k}}}    \Delta+\eta \check{b}\big[1+\log\big(1+\mathrm{e}^{Z_{k}}\big)\big]^{-1}   \Delta
   \bigg\}.
\end{align*}
Then it follows from the inequality $u  >  \log u$ for any $u > 0$ and  non-positivity of $a(r_k)$  that
\begin{align}\la{eq_4.17}
&\E \Big[\big[1+\log\big(1+\mathrm{e}^{\bar{Z}_{k+1}}\big)\big]^{\eta}\xi^{\eta,c}_{r_{k+1}} \big|\mathcal{F}_{t_k}\Big]\nn\\
 \leq&\xi^{\eta,c}_{r_{k}}\big[1+\log\big(1+\mathrm{e}^{Z_{k}}\big)\big]^{\eta} \bigg\{1+\frac{1}{\xi^{\eta,c}_{r_{k}}}\sum_{j\in \mathbb{S}}\xi^{\eta,c}_{j} \gamma_{r_k j}\t
-\eta a(r_k) \frac{1+\log\big(1+\mathrm{e}^{Z_{k}}\big)-1}{ 1+\log\big(1+\mathrm{e}^{Z_{k}}\big)} \Delta +o(\t)\nn\\
&~~~~~~~~~~~~~~~~~~~~ ~+2\eta a(r_k) \big[1+\log\big(1+\mathrm{e}^{Z_{k}}\big)\big]^{-1} \Delta  +\eta \check{b}\big[1+\log\big(1+\mathrm{e}^{Z_{k}}\big)\big]^{-1}   \Delta  \bigg\}\nn\\
=&\xi^{\eta,c}_{r_{k}}\big[1+\log\big(1+\mathrm{e}^{Z_{k}}\big)\big]^{\eta} \bigg\{1+\frac{1}{\xi^{\eta,c}_{r_{k}}}\sum_{j\in \mathbb{S}}\xi^{\eta,c}_{j} \gamma_{r_k j}\t
-\eta a(r_k)  \Delta+\eta \check{b}\big[1+\log\big(1+\mathrm{e}^{Z_{k}}\big)\big]^{-1}   \Delta \nn\\
&~~~~~~+o(\t)+3\eta a(r_k) \big[1+\log\big(1+\mathrm{e}^{Z_{k}}\big)\big]^{-1} \Delta
      \bigg\}.
\end{align}
By the properties of the generator, we have
\begin{align}\la{eq_4.19}
\frac{1}{\xi^{\eta,c}_{i}}\sum_{j\in \mathbb{S}}\xi^{\eta,c}_{j} \gamma_{i j}=&\frac{1}{1-c_i \eta}\sum_{j=1}^m(1-c_j \eta) \gamma_{i j}\nn\\
=&\frac{\eta}{1-c_i \eta}\sum_{j\neq i}  \gamma_{i j}(c_i-c_j)\nn\\
=&-\frac{\eta}{1-c_i \eta}\sum_{j=1}^m \gamma_{i j} c_j\nn\\
=&-\eta\bigg(\sum_{j=1}^m \gamma_{i j} c_j+\frac{c_i \eta}{1-c_i \eta}\sum_{j=1}^m \gamma_{i j} c_j\bigg).
\end{align}
It follows from \eqref{eq:2.7}, \eqref{eq_4.17} and \eqref{eq_4.19} that
\begin{align}\la{eq_4.20}
&\E \Big[\big[1+\log\big(1+\mathrm{e}^{\bar{Z}_{k+1}}\big)\big]^{\eta}\xi^{\eta,c}_{r_{k+1}} \big|\mathcal{F}_{t_k}\Big]\nn\\
\leq&\xi^{\eta,c}_{r_{k}}\big[1+\log\big(1+\mathrm{e}^{Z_{k}}\big)\big]^{\eta} \bigg\{1-\eta\bigg(\pi a+\frac{c_{r_k} \eta}{1-c_{r_k} \eta}\sum_{j=1}^m \gamma_{{r_k} j} c_j\bigg)\t+o(\t)
 \nn\\
&~~~~~~~~~~~~~~~~~~~~~+\eta \check{b}\big[1+\log\big(1+\mathrm{e}^{Z_{k}}\big)\big]^{-1}   \Delta+3\eta a(r_k) \big[1+\log\big(1+\mathrm{e}^{Z_{k}}\big)\big]^{-1} \Delta
     \bigg\}.
\end{align}
Choose a constant $0<\eta_1\leq \eta_0$ such that for any $0<\eta\leq \eta_1$,
$$
\pi a+\frac{c_i \eta}{1-c_i \eta}\sum\limits_{j=1}^{m}\gamma_{ij}c_j>0,\ \ i\in \mathbb{S}.
$$
Now, choose a positive constant $\kappa=(\kappa(\eta))<1$ sufficiently small such that it satisfies
\begin{eqnarray}\label{eq:2.16}
\hbar_i:=\pi a+\frac{c_i \eta}{1-c_i \eta}\sum\limits_{j=1}^{m}\gamma_{ij}c_j-\frac{\kappa}{\eta}>0.
\end{eqnarray}
Then, by \eqref{eq_4.20} and \eqref{eq:2.16}, we have
\begin{align*}
&\E \Big[\big[1+\log\big(1+\mathrm{e}^{\bar{Z}_{k+1}}\big)\big]^{\eta}\xi^{\eta,c}_{r_{k+1}} \big|\mathcal{F}_{t_k}\Big]\nn\\
\leq  &\xi^{\eta,c}_{r_{k}}\big[1+\log\big(1+\mathrm{e}^{Z_{k}}\big)\big]^{\eta} \bigg\{1-\kappa \t-\eta\hbar_k\t
 +\eta \check{b}\big[1+\log\big(1+\mathrm{e}^{Z_{k}}\big)\big]^{-1}   \Delta \nn\\
&~~~~~~+o(\t)+3\eta a(r_k) \big[1+\log\big(1+\mathrm{e}^{Z_{k}}\big)\big]^{-1} \Delta
     \bigg\}.
\end{align*}
 Choose    $ {\Delta}_1 ^*\in(0, 1)$ sufficiently small
such that
${\Delta}_1^*< 2/\kappa$ and  $o\left({\Delta}_1^*\right)\leq \kappa \eta {\Delta}_1^*/2.
$ For any $\Delta\in(0, {\Delta}_1^*]$, yields
\begin{align*}
&\E \Big[\big[1+\log\big(1+\mathrm{e}^{\bar{Z}_{k+1}}\big)\big]^{\eta}\xi^{\eta,c}_{r_{k+1}} \big|\mathcal{F}_{t_k}\Big]\nn\\
 \leq&\xi^{\eta,c}_{r_{k}}\big[1+\log\big(1+\mathrm{e}^{Z_{k}}\big)\big]^{\eta} \Big(1-\frac{\kappa}{2}\Delta \Big) \nn\\
  &-\eta \xi^{\eta,c}_{r_{k}}\bigg\{\hbar_k\big[1+\log\big(1+\mathrm{e}^{Z_{k}}\big)\big]^{\eta} -\big(3\check{a}+ \check{b}\big)  \big[1+\log\big(1+\mathrm{e}^{Z_{k}}\big)\big]^{\eta-1} \bigg\}\Delta\nn\\
      \leq&\Big(1-\frac{\kappa}{2}\Delta \Big)\big[1+\log\big(1+\mathrm{e}^{\bar{Z}_{k}}\big)\big]^{\eta}\xi^{\eta,c}_{r_{k}} +C \Delta
\end{align*}
for any integer $k\geq 0$. Repeating this procedure arrives at
\begin{align*} 
 \E \Big[\Big(1+\log\big(1+\mathrm{e}^{\bar{Z}_{k}}\big)\Big)^{\eta}\xi^{\eta,c}_{r_{k}}\Big]
 \leq& \mathrm{e}^{- \frac{\kappa}{2} k\t}\Big(1+\log\big(1+x_{0} \big)\Big)^{\eta}\xi^{\eta,c}_{\ell}+\frac{2C}{\kappa}\left[1-\left(1-  \frac{\kappa}{2}  \t\right)^{k}\right].
\end{align*}
Therefore,
$
\sup\limits_{k\geq0}\E \Big[ \log^{\eta}\big(1+\mathrm{e}^{\bar{Z}_{k}}\big) \Big] \leq C.
$
 The desired assertion follows.\eproof

The proofs of both below lemmas can be found in Appendix B.
\begin{lemma}\la{log_le*6.10}
If $\hat{a}>0$, for any $p>0$ and $\Delta\in (0, 1)$, the truncated EM scheme defined by \eqref{TEM_1} has the property that
\begin{align}\la{1}
\sup_{  k \geq 0} \mathbb{E}\big[ X^{p}_{k}\big]=\sup_{  k \geq 0} \mathbb{E}\big[\mathrm{e}^{pZ_{k}}\big]
\leq \sup_{  k \geq 0} \mathbb{E}\big[\mathrm{e}^{p\bar{Z}_{k}}\big]
\leq  C.
\end{align}
\end{lemma}

\begin{theorem}
Under the condition  of Lemma \ref{L:3}, the numerical solutions $X_{k}$ are stochastically ultimately upper bounded.
\end{theorem}
{\bf  Proof.}~~The proof is an application of Chebyshev's inequality, so we omit it.\eproof
Next we continue to consider the case $\pi a =0$, we can get the following results.
\begin{lemma}\la{log:Le5.6}
Under the condition  of Lemma \ref{log:Le5.4},
there is a constant $\Delta_2^{*}\in (0, 1)$ such that the scheme
 \eqref{TEM_1} has the property that
\begin{align*}
 \lim_{k\rightarrow \infty}\E \big[ X^{\rho}_{k} \big]= \lim_{k\rightarrow \infty}\E \big[ \mathrm{e}^{\rho Z_{k}} \big] =0
\end{align*}
for any $\Delta\in (0, \Delta_2^{*}]$, where $\rho$ is defined in Lemma \ref{log:Le5.4}.
  \end{lemma}

 Moreover, we can also get the following result.
\begin{lemma}\la{*log:Le5.8}
For any $\Delta \in (0,1)$ and initial value $(x_0, \ell)\in \mathbb{R}_+ \times \mathbb{S}$, the scheme
 \eqref{TEM_1} has the property that
\begin{align*}
\limsup_{k\rightarrow \infty}\frac{\log(X_{k})}{k\Delta}=\limsup_{k\rightarrow \infty}\frac{Z_{k}}{k\Delta}\leq\limsup_{k\rightarrow \infty}\frac{\bar{Z}_{k}}{k\Delta}
\leq\pi \beta~~~~~ a.s.
\end{align*}
In particularly, when $\pi a=0$,
\begin{align*}
\limsup_{k\rightarrow \infty}\frac{\log(X_{k})}{k\Delta}=\limsup_{k\rightarrow \infty}\frac{Z_{k}}{k\Delta}=\limsup_{k\rightarrow \infty}\frac{\bar{Z}_{k}}{k\Delta}
=\pi \beta~~~~~ a.s.
\end{align*}
  \end{lemma}
{\bf  Proof.}~~By the scheme
 \eqref{log_eq5*}, we have
\begin{align*}
\bar{Z}_{k+1}=&Z_{k}+\big(\beta(r_k)-a(r_k) \mathrm{e}^{Z_k}\big)\Delta +\sigma(r_k) \Delta B_k
\leq \bar{Z}_{k}+ \beta(r_k) \Delta +\sigma(r_k) \Delta B_k\nn\\
\leq&Z_{k-1}+ \beta(r_{k-1}) \Delta +\sigma(r_{k-1}) \Delta B_{k-1}+ \beta(r_k) \Delta +\sigma(r_k) \Delta B_k\nn\\
\leq&Z_{0}+ \sum_{j=0}^{k}\beta(r_j) \Delta +\sum_{j=0}^{k}\sigma(r_j) \Delta B_j.
\end{align*}
By the strong law of large numbers for martingales (see \cite[Theorem 1.6]{Mao06}), we  have
$$
\lim_{k\rightarrow \infty}\frac{\sum_{j=0}^{k-1}\sigma(r_j)\Delta B_j}{k\Delta}=0 ~~~~a.s.
$$
Then, by the ergodic property of the Markov chain (see, e.g., \cite{Anderson1991}), we compute
\begin{align*}
\lim_{k \rightarrow \infty } \frac { 1 } {k }  \sum_{j=0}^{k-1} \beta(r_j) \Delta & = \sum _ { i \in \mathbb { S } } \pi_ { i}\beta(i)\Delta=\pi \beta \Delta~~~ a. s.
\end{align*}
which implies
\begin{align*}
\limsup_{k\rightarrow \infty}\frac{Z_{k}}{k\Delta}
\leq\pi \beta~~~~~~~~~~~a.s.
\end{align*}
Particularly, when $\pi a=0$,  we have $a(\cdot)\equiv0$ and
\begin{align*}
Z_{k+1}
= Z_{k}+ \beta(r_k) \Delta +\sigma(r_k) \Delta B_k
=  Z_{0}+\sum_{i=0}^{k} \beta(r_k) \Delta + \sum_{i=0}^{k}\sigma(r_k)\Delta B_i,
\end{align*}
then the required assertion follows. \eproof

On the other hand, to show the numerical solutions $X_k$ defined by \eqref{TEM_1} is  stochastically ultimately lower bounded, we need the following lemma.
\begin{lemma}\la{log:Le5.9}
Under the condition  of Lemma \ref{log:Le5.8} and $\hat{a}>0$,
there is a constant $\Delta_3^{*}\in (0, 1)$ such that the scheme
 \eqref{TEM_1} has the property that
\begin{align*}
\sup\limits_{k\geq0}\mathbb E \big[ X^{-\vartheta}_{k}\big]=\sup\limits_{k\geq0}\mathbb E \big[\mathrm{e}^{-\vartheta Z_{k}}\big]\leq C
\end{align*}
for any $\Delta\in (0, \Delta_3^{*}]$, where $\vartheta$ is defined in Lemma \ref{log:Le5.8}.
  \end{lemma}
{\bf  Proof.}~~By \eqref{yhf0602}, we have
\begin{align}\la{logeq:5.45}
(1+\mathrm{e}^{-\bar{Z}_{k+1}})\leq (1+\mathrm{e}^{-Z_{k}})(1+\zeta_k),
\end{align}
where $\zeta_k$ is defined by \eqref{yhf0602}.
It follows from Lemma \ref{Le:Yin} (1) that the system of equations
\begin{eqnarray*}
\Gamma d=- \beta+(\pi \beta)\I_{m}
\end{eqnarray*}
has a solution $d=(d_1, \cdots, d_m)^T\in \mathbb{R}^m$. Then we have
\begin{eqnarray}\label{yhf1}
   \beta(i)+\sum_{j=1}^{m}\gamma_{ij}d_j= \pi \beta>0.
\end{eqnarray}
Choose a constant $0<\vartheta_1<1$ such that for each $0<\vartheta\leq \vartheta_1$,
\begin{eqnarray*}
\xi^{\vartheta,d}_i:=1-d_i \vartheta>0,\ \ \ i\in \mathbb{S}.
\end{eqnarray*}
For any $0<\vartheta<1$, by  virtue of \cite[Lemma 3.3]{li2018jde}, it follows from \eqref{logeq:5.45} that
\begin{align}\la{logeq:5.46}
&\E \Big[\big(1+\mathrm{e}^{-\bar{Z}_{k+1}}\big)^{\vartheta}\xi^{\vartheta,d}_{r_{k+1}} \big|\mathcal{F}_{t_k}\Big]\nn\\
\leq&\big(1+\mathrm{e}^{-Z_{k}}\big)^{\vartheta} \bigg\{\E\big[\xi^{\vartheta,d}_{r_{k+1}}\big|\mathcal{F}_{t_k}\big]+ \vartheta \E\big[\zeta_k\xi^{\vartheta,d}_{r_{k+1}}\big|\mathcal{F}_{t_k}\big] + \frac{\vartheta(\vartheta-1)}{2} \E\big[\zeta_k^2\xi^{\vartheta,d}_{r_{k+1}}\big|\mathcal{F}_{t_k}\big]\nn\\
 &~~~~~~~~~~~~~~~~~~~~~~ + \frac{\vartheta(\vartheta-1)(\vartheta-2)}{6} \E\big[\zeta_k^3\xi^{\vartheta,d}_{r_{k+1}}
 \big|\mathcal{F}_{t_k}\big]\bigg\}.
\end{align}
Then, making use of the techniques in the proof of Lemma \ref{log:Le5.6} as well as Lemma \ref{leB}  and \eqref{logeq:5.10} yields
\begin{align*}
\E\big[\mathcal{U}_k\big|\mathcal{G}_{t_k}\big]
\leq C\big(\Delta^{3(1-\theta )} + \Delta^{\frac{3}{2}}\big),
\end{align*}
and
\begin{align*}
\E\big[\zeta_k \big|\mathcal{G}_{t_k}\big]\leq&
\frac{\mathrm{e}^{-Z_{k}}}{1+\mathrm{e}^{-Z_{k}}}\Big( \big(a(r_k) \mathrm{e}^{Z_k}-\beta(r_k)\big)\Delta+ C\Delta^{2(1-\theta )}
+\frac{\sigma^2(r_k) \Delta  }{2}
+\E\big[\mathcal{U}_k\big|\mathcal{G}_{t_k}\big]\Big).
\end{align*}
  By the Markov property, we derive that
\begin{align}\la{log:eq:5.49}
\!\!\E\big[\zeta_k  \xi^{\vartheta,d}_{r_{k+1}} \big|\mathcal{F}_{t_k}\big]\leq&
\frac{\mathrm{e}^{-Z_{k}}}{1+\mathrm{e}^{-Z_{k}}}\Big[ \Big(a(r_k) \mathrm{e}^{Z_k}-\beta(r_k)+\frac{\sigma^2(r_k)  }{2}\Big)\Delta+ C \Delta^{2(1-\theta )} \Big]\E\big[ \xi^{\vartheta,d}_{r_{k+1}} \big|\mathcal{F}_{t_k}\big]\nn\\
\leq&
\frac{\mathrm{e}^{-Z_{k}}}{1+\mathrm{e}^{-Z_{k}}}\Big[ \xi^{\vartheta,d}_{r_{k}}\Big(a(r_k) \mathrm{e}^{Z_k}-\beta(r_k)+\frac{\sigma^2(r_k)  }{2}\Big)\Delta\Big]+ C \Delta^{2(1-\theta )}\nn\\
=&
  \xi^{\vartheta,d}_{r_{k}}\Big(-\beta(r_k)+\frac{\sigma^2(r_k)  }{2}\Big)\Delta-\frac{1}{1+\mathrm{e}^{-Z_{k}}} \xi^{\vartheta,d}_{r_{k}}\Big(-\beta(r_k)
  +\frac{\sigma^2(r_k)  }{2}\Big)\Delta \nn\\
&+\frac{ a(r_k)}{1+\mathrm{e}^{-Z_{k}}}  \xi^{\vartheta,d}_{r_{k}} \Delta + C \Delta^{2(1-\theta )}\nn\\
\leq&-
  \xi^{\vartheta,d}_{r_{k}}\Big(\beta(r_k)\!-\frac{\sigma^2(r_k)  }{2}\Big)\Delta
 +\frac{ \beta(r_k)+a(r_k)}{1+\mathrm{e}^{-Z_{k}}}  \xi^{\vartheta,d}_{r_{k}} \Delta + C \Delta^{2(1-\theta )}.
\end{align}
Using the techniques in the proof of Lemma \ref{L:3}, we show that
\begin{align*}
\E\big[\zeta_k^2 \big|\mathcal{G}_{t_k}\big]=& \frac{\mathrm{e}^{-2Z_{k}}}{(1+\mathrm{e}^{-Z_{k}})^2}\E\bigg[\bigg( \big(a(r_k) \mathrm{e}^{Z_k}-\beta(r_k)\big)\Delta+ C \Delta^{2(1-\theta )}
+\frac{\sigma^2(r_k) (\Delta B_k)^2}{2}\nn\\
&~~~~~~~~~~~~~~~~~~~~ -\sigma(r_k) \Delta B_k
-\sigma(r_k)\Delta B_k\big(a(r_k) \mathrm{e}^{Z_k}-\beta(r_k)\big)\Delta
+\mathcal{U}_k\bigg)^2\Big|\mathcal{G}_{t_k}\bigg]\nn\\
\geq& \frac{\mathrm{e}^{-2Z_{k}}}{(1+\mathrm{e}^{-Z_{k}})^2}\E\bigg[ \sigma^2(r_k) (\Delta B_k)^2-2\sigma(r_k) \Delta B_k \Big( \big(a(r_k) \mathrm{e}^{Z_k}-\beta(r_k)\big)\Delta
+ C \Delta^{2(1-\theta )}\nn\\
&~~~~~~~~~~~~~~~~~~
+\frac{\sigma^2(r_k) (\Delta B_k)^2}{2} -\sigma(r_k)\Delta B_k\big(a(r_k) \mathrm{e}^{Z_k}-\beta(r_k)\big)\Delta
+\mathcal{U}_k\Big) \Big|\mathcal{G}_{t_k}\bigg]\nn\\
=& \frac{\mathrm{e}^{-2Z_{k}}}{(1+\mathrm{e}^{-Z_{k}})^2}\bigg[ \sigma^2(r_k) \Delta +2\sigma^2(r_k) \big(a(r_k) \mathrm{e}^{Z_k}-\beta(r_k)\big)\Delta^2
-2\sigma(r_k)\E\Big( \Delta B_k \mathcal{U}_k \big|\mathcal{G}_{t_k}\Big)\bigg]\nn\\
\geq& \frac{\mathrm{e}^{-2Z_{k}}}{(1+\mathrm{e}^{-Z_{k}})^2}\Big[\sigma^2(r_k) \Delta -2|\breve{\sigma}|^2  \big(\check{a}  K\Delta^{-\theta }+|\breve{\beta}|\big)\Delta^2
 \Big]\nn\\
 \geq& \frac{(1+\mathrm{e}^{-Z_{k}}-1)^2}{(1+\mathrm{e}^{-Z_{k}})^2} \sigma^2(r_k) \Delta -C\Delta^{2-\theta },
\end{align*}
which implies that
\begin{align}\la{log:eq:5.50}
\E\big[\zeta_k^2 \xi^{\vartheta,d}_{r_{k+1}} \big|\mathcal{F}_{t_k}\big]
 \geq& \frac{(1+\mathrm{e}^{-Z_{k}}-1)^2}{(1+\mathrm{e}^{-Z_{k}})^2} \sigma^2(r_k)\xi^{\vartheta,d}_{r_{k}} \Delta -C\Delta^{2-\theta }\nn\\
=& \Big(1-\frac{2}{ 1+\mathrm{e}^{-Z_{k}} }+\frac{1}{(1+\mathrm{e}^{-Z_{k}})^2}\Big) \sigma^2(r_k)\xi^{\vartheta,d}_{r_{k}} \Delta -C\Delta^{2-\theta }.
\end{align}
In addition,
\begin{align*}
\E\Big[\big(\Delta B_k\big) \mathcal{U}_k \Big|\mathcal{G}_{t_k}
 \Big]=0,~~~
\E\Big[\big(\Delta B_k\big) \mathcal{U}_k^2\Big|\mathcal{G}_{t_k}
 \Big] =0.
 \end{align*}
This together with \eqref{logeq:5.10} as well as Lemma \ref{leB}  yields
\begin{align*}
\E\big[\zeta_k^3 \big|\mathcal{G}_{t_k}\big]=& \frac{\mathrm{e}^{-3Z_{k}}}{(1+\mathrm{e}^{-Z_{k}})^3}\E\bigg[\bigg( \big(a(r_k) \mathrm{e}^{Z_k}-\beta(r_k)\big)\Delta+C \Delta^{2(1-\theta )}
+\frac{\sigma^2(r_k) (\Delta B_k)^2}{2}\nn\\
&~~~~~~~~~~~~~~~~~~-\sigma(r_k) \Delta B_k
-\sigma(r_k)\Delta B_k\big(a(r_k) \mathrm{e}^{Z_k}-\beta(r_k)\big)\Delta
+\mathcal{U}_k\bigg)^3\Big|\mathcal{G}_{t_k}\bigg]\nn\\
=& \frac{\mathrm{e}^{-3Z_{k}}}{(1+\mathrm{e}^{-Z_{k}})^3}\E\bigg[
-\sigma^3(r_k) \big(\Delta B_k\big)^3\Big(1+
 \big(a(r_k) \mathrm{e}^{Z_k}-\beta(r_k)\big)\Delta\Big)^3
\nn\\
&+3\sigma^2(r_k) \big(\Delta B_k\big)^2\Big(1+
 \big(a(r_k) \mathrm{e}^{Z_k}-\beta(r_k)\big)\Delta\Big)^2\nn\\
&~~~~~~~~~~~~~~~~~\times \Big( \big(a(r_k) \mathrm{e}^{Z_k}-\beta(r_k)\big)\Delta+ C \Delta^{2(1-\theta )}
+\frac{\sigma^2(r_k) (\Delta B_k)^2}{2}+\mathcal{U}_k
 \Big)\nn\\
 &-3\sigma(r_k) \big(\Delta B_k\big) \Big(1+
 \big(a(r_k) \mathrm{e}^{Z_k}-\beta(r_k)\big)\Delta\Big) \nn\\
&~~~~~~~~~~~~~~~~~\times \Big( \big(a(r_k) \mathrm{e}^{Z_k}-\beta(r_k)\big)\Delta+C \Delta^{2(1-\theta )}
+\frac{\sigma^2(r_k) (\Delta B_k)^2}{2}+\mathcal{U}_k
 \Big)^2\nn\\
&+\Big( \big(a(r_k) \mathrm{e}^{Z_k}-\beta(r_k)\big)\Delta+ C \Delta^{2(1-\theta )}
+\frac{\sigma^2(r_k) (\Delta B_k)^2}{2}+\mathcal{U}_k
 \Big)^3\Big|\mathcal{G}_{t_k}\bigg]\nn\\
 =& \frac{\mathrm{e}^{-3Z_{k}}}{(1+\mathrm{e}^{-Z_{k}})^3}\E\bigg[
 3\sigma^2(r_k) \big(\Delta B_k\big)^2\Big(1+
 \big(a(r_k) \mathrm{e}^{Z_k}-\beta(r_k)\big)\Delta\Big)^2\nn\\
&~~~~~~~~~~~~~~~~~\times \Big( \big(a(r_k) \mathrm{e}^{Z_k}-\beta(r_k)\big)\Delta+ C \Delta^{2(1-\theta )}
+\frac{\sigma^2(r_k) (\Delta B_k)^2}{2}+\mathcal{U}_k
 \Big)\nn\\
&+\Big( \big(a(r_k) \mathrm{e}^{Z_k}-\beta(r_k)\big)\Delta+ C \Delta^{2(1-\theta )}
+\frac{\sigma^2(r_k) (\Delta B_k)^2}{2}+\mathcal{U}_k
 \Big)^3\Big|\mathcal{G}_{t_k}\bigg]\nn\\
\leq& \frac{\mathrm{e}^{-3Z_{k}}}{(1+\mathrm{e}^{-Z_{k}})^3} \bigg[
 3\sigma^2(r_k)\Big(1+
 \big(a(r_k) \mathrm{e}^{Z_k}-\beta(r_k)\big)\Delta\Big)^2\nn\\
&~~~~~~~~~~~ \times \Big( \big(a(r_k) \mathrm{e}^{Z_k}-\beta(r_k)\big)\Delta^2+ C \Delta^{3-2\theta }
+\frac{3\sigma^2(r_k) \Delta^2}{2}+ \E\big[\big(\Delta B_k\big)^2\mathcal{U}_k\big|\mathcal{G}_{t_k}\big]
 \Big)\nn\\
&+16\Big( \big|a(r_k) \mathrm{e}^{Z_k}-\beta(r_k)\big|^3\Delta^3+ C \Delta^{6(1-\theta )}
+\frac{5!!}{8}\sigma^6(r_k) \Delta^3+\E\big[\mathcal{U}_k^3\big|\mathcal{G}_{t_k}\big]
 \Big)\bigg]\nn\\
\leq& \frac{\mathrm{e}^{-3Z_{k}}}{(1+\mathrm{e}^{-Z_{k}})^3} \bigg[
 C\big( \Delta^{2-\theta }+   \Delta^{3-2\theta }+   \Delta^{2}
 +  \Delta^{4-3\theta }+\Delta^{\frac{5}{2}}
 \big)\nn\\
&~~~~~~~~~~~~~~~~~~~~~+C\big(\Delta^{3(1-\theta )}+  \Delta^{6(1-\theta )}+  \Delta^{3}
 + \Delta^{9(1-\theta )}+ \Delta^{\frac{9}{2}}
 \big)\bigg],
\end{align*}
which implies that
\begin{align}\la{log:eq:5.54}
 \E\big[\zeta_k^3 \xi^{\vartheta,d}_{r_{k+1}} \big|\mathcal{F}_{t_k}\big]
\leq   C
  \big(\Delta^{2-\theta  }+\Delta^{3(1-\theta )} \big).
\end{align}
Substituting  \eqref{log:eq:5.49}, \eqref{log:eq:5.50} and \eqref{log:eq:5.54}  into \eqref{logeq:5.46},  we derive  that
\begin{align}\la{yhf3}
&\E \Big[\big(1+\mathrm{e}^{-\bar{Z}_{k+1}}\big)^{\vartheta}\xi^{\vartheta,d}_{r_{k+1}} \big|\mathcal{F}_{t_k}\Big]\nn\\
\leq&\big(1+\mathrm{e}^{-Z_{k}}\big)^{\vartheta} \bigg\{\xi^{\vartheta,d}_{r_{k}}+\sum_{j\in \mathbb{S}}\xi^{\vartheta,d}_{j}\left(\gamma_{r_k j}\t+o(\t)\right) + C\Delta^{2(1-\theta )}
  + C \Delta^{2-\theta  }   \nn\\
&~~~~~~~~~~~~~~~~~~-
  \vartheta\xi^{\vartheta,d}_{r_{k}}\Big(\beta(r_k) -\frac{\sigma^2(r_k)  }{2}\Big)\Delta
 +\frac{ \vartheta(\beta(r_k)+a(r_k))}{1+\mathrm{e}^{-Z_{k}}}  \xi^{\vartheta,d}_{r_{k}} \Delta \nn\\
&~~~~~~~~~~~~~~~~~~+\frac{\vartheta(\vartheta-1)}{2}\Big(1
-\frac{2}{ 1+\mathrm{e}^{-Z_{k}} }+\frac{1}{(1+\mathrm{e}^{-Z_{k}})^2}\Big) \sigma^2(r_k)\xi^{\vartheta,d}_{r_{k}} \Delta  \bigg\}\nn\\
\leq&\big(1+\mathrm{e}^{-Z_{k}}\big)^{\vartheta} \bigg\{\xi^{\vartheta,d}_{r_{k}}+\sum_{j\in \mathbb{S}}\xi^{\vartheta,d}_{j} \gamma_{r_k j}\t+o(\t)-
  \vartheta\xi^{\vartheta,d}_{r_{k}}\Big(\beta(r_k) -\frac{\sigma^2(r_k)  }{2}-\frac{\vartheta(\vartheta-1)}{2}\sigma^2(r_k)\Big)\Delta    \nn\\
&~~~~~~~~~~~~~
 +\frac{\vartheta}{1+\mathrm{e}^{-Z_{k}}}\Big( \beta(r_k)+a(r_k)-(\vartheta-1)\sigma^2(r_k) \Big)  \xi^{\vartheta,d}_{r_{k}} \Delta  + \frac{\vartheta(\vartheta-1)\sigma^2(r_k)}{2(1+\mathrm{e}^{-Z_{k}})^2}  \xi^{\vartheta,d}_{r_{k}} \Delta \bigg\}\nn\\
\leq&\big(1+\mathrm{e}^{-Z_{k}}\big)^{\vartheta} \xi^{\vartheta,d}_{r_{k}}\bigg[1+\vartheta\Big(\frac{1}{\vartheta\xi^{\vartheta,d}_{r_{k}}}\sum_{j\in \mathbb{S}}\xi^{\vartheta,d}_{j} \gamma_{r_k j}- \beta(r_k)
+\frac{\vartheta \sigma^2(r_k) }{2}\Big)
   \t +o(\t)\bigg]\nn\\
&+\vartheta\Big(a(r_k)+\beta(r_k)
-(\vartheta-1)\sigma^2(r_k)
 \Big) \big(1+\mathrm{e}^{-Z_{k}}\big)^{\vartheta-1}\xi^{\vartheta,d}_{r_{k}}\Delta
 + \frac{\vartheta(\vartheta-1)\sigma^2(r_k) }{2(1+\mathrm{e}^{-Z_{k}})^{2-\vartheta}}  \xi^{\vartheta,d}_{r_{k}} \Delta.
\end{align}
By the properties of the generator, we have
\begin{align} \la{yhf4}
\frac{1}{\vartheta\xi^{\vartheta,d}_{i}}\sum_{j\in \mathbb{S}}\xi^{\vartheta,d}_{j} \gamma_{i j}= -\bigg(\sum_{j=1}^m \gamma_{i j} d_j+\frac{d_i \vartheta}{1-d_i \vartheta}\sum_{j=1}^m \gamma_{i j} d_j\bigg).
\end{align}
It follows from \eqref{yhf1}, \eqref{yhf3} and \eqref{yhf4} that
\begin{align*}
&\E \Big[\big(1+\mathrm{e}^{-\bar{Z}_{k+1}}\big)^{\vartheta}\xi^{\vartheta,d}_{r_{k+1}} \big|\mathcal{F}_{t_k}\Big]\nn\\
\leq&\big(1+\mathrm{e}^{-Z_{k}}\big)^{\vartheta} \xi^{\vartheta,d}_{r_{k}}\bigg[1-\vartheta\Big(
\sum_{j=1}^m \gamma_{r_{k} j} d_j+\frac{d_{r_{k}} \vartheta}{1-d_{r_{k}} \vartheta}\sum_{j=1}^m \gamma_{{r_{k}} j} d_j+ \beta(r_k)
-\frac{\vartheta \sigma^2(r_k) }{2}\Big)
   \t +o(\t)\bigg]\nn\\
&+\vartheta\Big(a(r_k)+\beta(r_k)
-(\vartheta-1)\sigma^2(r_k)
 \Big) \big(1+\mathrm{e}^{-Z_{k}}\big)^{\vartheta-1}\xi^{\vartheta,d}_{r_{k}}\Delta
 + \frac{\vartheta(\vartheta-1)\sigma^2(r_k) }{2(1+\mathrm{e}^{-Z_{k}})^{2-\vartheta}}  \xi^{\vartheta,d}_{r_{k}} \Delta\nn\\
=&\big(1+\mathrm{e}^{-Z_{k}}\big)^{\vartheta} \xi^{\vartheta,d}_{r_{k}}\bigg[1-\vartheta\Big(
 \pi\beta+\frac{d_{r_{k}} \vartheta}{1-d_{r_{k}} \vartheta}\sum_{j=1}^m \gamma_{{r_{k}} j} d_j
-\frac{\vartheta \sigma^2(r_k) }{2}\Big)
   \t +o(\t)\bigg]\nn\\
&+\vartheta\Big(a(r_k)+\beta(r_k)
-(\vartheta-1)\sigma^2(r_k)
 \Big) \big(1+\mathrm{e}^{-Z_{k}}\big)^{\vartheta-1}\xi^{\vartheta,d}_{r_{k}}\Delta
 + \frac{\vartheta(\vartheta-1)\sigma^2(r_k) }{2(1+\mathrm{e}^{-Z_{k}})^{2-\vartheta}}  \xi^{\vartheta,d}_{r_{k}} \Delta.
\end{align*}
Choose a small constant $0<\vartheta_1\leq \vartheta_0$ such that for any $0<\vartheta\leq \vartheta_1$,
\begin{eqnarray*}
 \pi \beta+\vartheta\bigg(\frac{d_i}{1-d_i \vartheta}\sum\limits_{j=1}^{m}\gamma_{ij}d_j-\frac{\sigma^2(i)}{2} \bigg)>0,~~~~~i\in \mathbb{S}.
\end{eqnarray*}
Now, choose a positive constant $\bar{\lambda}=\bar{\lambda}(\vartheta)<1$ sufficiently small such that it satisfies
\begin{eqnarray*}
 \pi \beta+\vartheta\bigg(\frac{d_i}{1-d_i \vartheta}\sum\limits_{j=1}^{m}\gamma_{ij}d_j-\frac{\sigma^2(i)}{2} \bigg)-\frac{\bar{\lambda}}{\vartheta}>0,~~~~~i\in \mathbb{S}.
\end{eqnarray*}
Then  we have
\begin{align*}
&\E \Big[\big(1+\mathrm{e}^{-\bar{Z}_{k+1}}\big)^{\vartheta}\xi^{\vartheta,d}_{r_{k+1}} \big|\mathcal{F}_{t_k}\Big]\nn\\
\leq& \big(1+\mathrm{e}^{-Z_{k}}\big)^{\vartheta} \xi^{\vartheta,d}_{r_{k}}\Big(1- \frac{\bar{\lambda}}{2}
   \t +o(\t)\Big)-\frac{\bar{\lambda}}{2}\big(1+\mathrm{e}^{-Z_{k}}\big)^{\vartheta} \xi^{\vartheta,d}_{r_{k}}\t\nn\\
&+\vartheta\Big(a(r_k)+\beta(r_k)
-(\vartheta-1)\sigma^2(r_k)
 \Big) \big(1+\mathrm{e}^{-Z_{k}}\big)^{\vartheta-1}\xi^{\vartheta,d}_{r_{k}}\Delta
 + \frac{\vartheta(\vartheta-1)\sigma^2(r_k) }{2(1+\mathrm{e}^{-Z_{k}})^{2-\vartheta}}  \xi^{\vartheta,d}_{r_{k}} \Delta\nn\\
 \leq& \Big(1- \frac{\bar{\lambda}}{2}
   \t +o(\t)\Big)\big(1+\mathrm{e}^{-Z_{k}}\big)^{\vartheta} \xi^{\vartheta,d}_{r_{k}}+C\t.
\end{align*}
  Choose    $ {\Delta}_3 ^*\in(0, 1) $ sufficiently small
such that
$\Delta_3^*< 4/\bar{\lambda}$ and   $o\left(\Delta_3^*\right)\leq \bar{\lambda} \Delta_3^*/4.
$ For any $\Delta\in(0, \Delta_3^*]$  yields
\begin{align*}
 \E \Big[\big(1+\mathrm{e}^{-\bar{Z}_{k+1}}\big)^{\vartheta}\xi^{\vartheta,d}_{r_{k+1}} \big|\mathcal{F}_{t_k}\Big]
\leq \Big(1-\frac{\bar{\lambda}}{4}  \t \Big)\big(1+\mathrm{e}^{-Z_{k}}\big)^{\vartheta}\xi^{\vartheta,d}_{r_{k}}+C  \Delta
\end{align*}
for any integer $k\geq 0$.  Obviously,
\begin{align*}
 \E \Big[\big(1+\mathrm{e}^{-\bar{Z}_{k+1}}\big)^{\vartheta}\xi^{\vartheta,d}_{r_{k+1}} \Big]
\leq \Big(1-\frac{\bar{\lambda}}{4}  \t \Big)\E\Big[\big(1+\mathrm{e}^{-Z_{k}}\big)^{\vartheta}\xi^{\vartheta,d}_{r_{k}}\Big]+C  \Delta.
\end{align*}
Define $\Omega_k=\{\bar{Z}_k>\log(K\Delta^{-\theta})\}$, we have
\begin{align}\la{Omega}
\big(1+\mathrm{e}^{-Z_{k}}\big)^{\vartheta}
\xi^{\vartheta,d}_{r_{k}}=&\big(1+\mathrm{e}^{-Z_{k}}\big)^{\vartheta}\xi^{\vartheta,d}_{r_{k}}I_{\Omega_k}
+\big(1+\mathrm{e}^{-\bar{Z}_{k}}\big)^{\vartheta}\xi^{\vartheta,d}_{r_{k}}I_{\Omega_k^c}\nn\\
\leq &\big(1+K^{-1}\t^{\theta}\big)^{\vartheta}\xi^{\vartheta,d}_{r_{k}}I_{\Omega_k}
+\big(1+\mathrm{e}^{-\bar{Z}_{k}}\big)^{\vartheta}\xi^{\vartheta,d}_{r_{k}}\nn\\
\leq &C I_{\Omega_k}
+\big(1+\mathrm{e}^{-\bar{Z}_{k}}\big)^{\vartheta}\xi^{\vartheta,d}_{r_{k}}.
\end{align}
By Chebyshev's inequality and Lemma \ref{log_le*6.10},
\begin{align*}
 \E \Big[\big(1+\mathrm{e}^{-Z_{k}}\big)^{\vartheta}\xi^{\vartheta,d}_{r_{k}} \Big]
\leq &
C \E\big[I_{\Omega_k}\big]
+\E\Big[\big(1+\mathrm{e}^{-\bar{Z}_{k}}\big)^{\vartheta}\xi^{\vartheta,d}_{r_{k}}\Big]\nn\\
\leq&\Big(1-\frac{\bar{\lambda}}{4}  \t \Big)\E\Big[\big(1+\mathrm{e}^{-Z_{k-1}}\big)^{\vartheta}\xi^{\vartheta,d}_{r_{k-1}}\Big]+C  \Delta+C \mathbb{P}\big\{\bar{Z}_k>\log(K\Delta^{-\theta})\big\}\nn\\
\leq&\Big(1-\frac{\bar{\lambda}}{4}  \t \Big)\E\Big[\big(1+\mathrm{e}^{-Z_{k-1}}\big)^{\vartheta}\xi^{\vartheta,d}_{r_{k-1}}\Big]+C  \Big(\Delta+  \frac{\E\big[\mathrm{e}^{\bar{Z}_k/\theta}\big]}{K^{1/\theta}\Delta^{-1}}\Big)\nn\\
\leq&\Big(1-\frac{\bar{\lambda}}{4}  \t \Big)\E\Big[\big(1+\mathrm{e}^{-Z_{k-1}}\big)^{\vartheta}\xi^{\vartheta,d}_{r_{k-1}}\Big]+C  \Delta.
\end{align*}
Repeating this procedure arrives at
\begin{align*}
\E \Big[\big(1+\mathrm{e}^{-Z_{k}}\big)^{\vartheta}\xi^{\vartheta,d}_{r_{k}} \Big]
 \leq& \mathrm{e}^{- \frac{\bar{\lambda}}{4} k \t} (1+ x^{-1}_{0} )^{\vartheta}\xi^{\vartheta,d}_{\ell} +\frac{4C}{\bar{\lambda}}\left[1-\left(1-  \frac{\bar{\lambda}}{4}  \t\right)^{k}\right]\leq C
\end{align*}
for any integer $k\geq 0$.   Therefore,
$
\sup\limits_{k\geq0}\mathbb E \big[\mathrm{e}^{-\vartheta Z_{k}}\big]\leq C.
$
 The desired assertion follows.

\begin{theorem}
Under the condition  of Lemma \ref{log:Le5.9}, the numerical solutions $X_{k}$ are stochastically ultimately  lower bounded.
\end{theorem}
{\bf  Proof.}~~The proof is an application of Chebyshev's inequality, so we omit it.\eproof

\begin{theorem}\la{log:nu_permanence}
Under the condition  of Theorem \ref{log_permanence} and $\hat{a}>0$,  for any $\Delta\in (0, \Delta_1^{*}\wedge \Delta_3^{*}]$,  the numerical solutions $X_{k}$ are  stochastically permanent.
\end{theorem}

Moreover, we obtain the following improved necessary
and sufficient conditions for the dynamical behaviors of   the numerical solutions $X_{k}$ defined by  \eqref{TEM_1}.
 \begin{theorem}\label{log_Th*5.4}
Suppose that  $\pi \beta\neq0$. For any $\Delta\in (0, \Delta_1^{*}\wedge \Delta_3^{*}]$, if $\hat{a}>0$,
 \begin{itemize}
\item  the numerical solutions $X_{k}$ are
stochastically permanent if and only if $\pi \beta>0$;
\item the numerical solutions $X_{k}$ are
 almost surely extinctive if and only if $\pi \beta<0$.
     \end{itemize}
     In particular, if $\pi a=0$,
      \begin{itemize}
      \item the numerical solutions $X_{k}$ are
 almost surely extinctive if and only if $\pi \beta<0$;
         \item   almost all paths of $X_{k}$ increase at an exponential rate if and only if $\pi \beta>0$.
     \end{itemize}
\end{theorem}

\section{Stability in distribution}\la{sections5}
In this section,  we first give sufficient conditions that guarantee SDS  \eqref{logistic1} is asymptotically stable in distribution.  Then we show that the  explicit schemes \eqref{TEM_1} can approximate the invariant measure of SDS \eqref{logistic1} effectively.  From
this section as a standing assumption, we always assume
  $\hat{a}>0$.
  For the convenience of invariant measure study we introduce some notations. We write $(x^{x_0, \ell}_t,r^{\ell}_t)$ in lieu of $(x(t), r(t))$ to highlight the initial data $(x(0), r(0)) =(x_0,\ell)$. Following \cite[p.212]{Mao06}, we denote by $\mathcal{P}(\mathbb{R}_+\times  \mathbb{S})$ the space of all probability measures on $\mathbb{R}_+\times  \mathbb{S}$ and for $\bar{\mu}, \bar{\nu}\in\mathcal{P}(\mathbb{R}_+\times  \mathbb{S})$ define
$$
d _ { \mathbb{L} } \left(\bar{\mu} , \bar{\nu} \right) : = \sup _ { H \in \mathbb{ L } } \left| \sum _ { i = 1 } ^ { m } \int _ { \mathbb { R } _ { + }   } H\left( x , i \right) \bar{\mu} \left( dx , i \right) - \sum _ { i= 1 } ^ { m } \int _ { \mathbb { R } _ { + }   } H\left( x   , i\right) \bar{\nu} \left( d x   , i\right) \right|
$$
where
$$
\mathbb { L } : = \Big\{H: \mathbb { R } _ { + }  \times  \mathbb{S} \rightarrow \mathbb { R } : \left| H\left( x , i \right) - H\left( y  , j \right) \right| \leq \left| x  - y  \right| + |i - j| ,~~  | H( \cdot ,\cdot ) | \leq 1 \Big\}.
$$
By  virtue of Lemma \ref{log_le1}, SDS \eqref{logistic1} has a unique continuous positive solution $(x^{x_0, \ell}_t,r^{\ell}_t)$, which is a time-homogeneous Markov process. Let $\mathbf{P}_{t}(x_0, \ell;\mathrm{d}x\times \{i\})$ denote the transition probability of the process $(x^{x_0, \ell}_t,r^{\ell}_t)$.
\begin{defn}[\!\!\cite{Mao06}] The process $(x(t),r(t))$ is said to be asymptotically stable in distribution if there exists a probability
measure $\mu(\cdot\times\cdot)$ on $\mathbb { R } _ { + }  \times  \mathbb{S} $ such that the transition probability $\mathbf{P}_{t}(x_0, \ell;\cdot\times\cdot)$ of $(x(t),r(t))$ converges weakly to $\mu(\cdot\times\cdot)$ as
 $t\rightarrow \infty$ for every $(x_0,\ell)\in \mathbb { R } _ { + }  \times  \mathbb{S}$. SDS \eqref{logistic1} is said to be asymptotically stable in distribution if $(x(t),r(t))$ is asymptotically stable in distribution.
\end{defn}

It is easy to observe that Theorem \ref{log_permanence} guarantees that for any $(x_0,\ell)\in \mathbb{R}_+\times\mathbb{S}$, the family of transition probabilities $\{\mathbf{P}_{t}(x_0, \ell;\cdot\times\cdot): t\geq 0\}$ is tight. That is, for any $\varepsilon>0$ there is a compact subset $\mathbb{K}_+=\mathbb{K}_+(\varepsilon, x_0,\ell)$ of $\mathbb{R}_+$ such that
\begin{align}\la{log_eq_tight}
\mathbf{P}_{t}\big(x_0, \ell;\mathbb{K}_+\times\mathbb{S}\big)\geq 1-\varepsilon~~~\forall t\geq 0.
\end{align}
Next we give the existence and uniqueness of the invariant measure for the solution  $(x^{x_0, \ell}_t,r^{\ell}_t)$  of SDS \eqref{logistic1}.

\begin{lemma}[\!\!\cite{hu20181AMM}]\la{log_le*7.4}
The solutions of SDS \eqref{logistic1} satisfy
 $$\lim_{t\rightarrow \infty}\mathbb{E}\big| x^{x_0,\ell}_t- x^{\bar{x}_0,\ell}_t\big|= 0$$
 uniformly in $(x_0,\bar{x}_0,\ell)\in \mathbb{K}_+\times\mathbb{K}_+\times\mathbb{S}$, where $x^{x_0,\ell}_t$, $x^{\bar{x}_0,\ell}_t$ are respectively two solutions of SDS \eqref{logistic1} with initial values $x_0$, $\bar{x}_0$ for any compact subset $\mathbb{K}_+$ of $\mathbb{R}_+$, and $\ell\in \mathbb{S}$.
  \end{lemma}

\begin{lemma}[\!\!\cite{hu20181AMM}]\la{log_le*7.5}
For any compact subset $\mathbb{K}_+$ of $\mathbb{R}_{+}$ and  any $T\geq 0$,
$$
\sup _ { \left( x _ { 0 } , \ell \right) \in \mathbb{K}_+ \times \mathbb{S} } \mathbb{E} \left[ \sup _ { 0 \leq t \leq T } x^2( t )  \right]\leq C_T,
$$
where $x( t )$   is the solution of SDS \eqref{logistic1} with the initial value  $(x _ { 0 },\ell)\in \mathbb{K}_+\times \mathbb{S}$.
\end{lemma}

\begin{lemma}\la{log_le_le7.6}
If $\pi \beta>0$ hold, for any compact subset
$\mathbb{K}_+$ of $ \mathbb { R } _ { + } $,
$$
\sup _ { \left( x _ { 0 } , \ell \right) \in \mathbb{K}_+ \times \mathbb{S} } \mathbb{E} \left[ \sup _ { 0 \leq t \leq T } x^{-\vartheta}( t )  \right] \leq C_T, ~~~~~\forall~T \geq 0,
$$
where $\vartheta$ is given by   Lemma \ref{log:Le5.8}, $x( t )$   is the solution of SDS \eqref{logistic1} with the initial value   $(x _ { 0 },\ell)\in \mathbb{K}_+\times \mathbb{S}$.
\end{lemma}
{\bf  Proof.}~~
Borrowing the proof method of \cite[Lemma 3.5]{Li2009DCDS} we can get the desired result but omit the details to avoid duplication.\eproof

By Lemma \ref{log_le_le7.6} and Jensens's inequality,
$
\big(\mathbb{E}x^{\vartheta}(t)\big)^{-1}\leq \mathbb{E} \big[ x^{-\vartheta} (t) \big]\leq C_T.
$
So
\begin{align}\la{log_eq*eq7.8}
\mathbb{E}\Big[\sup_{0\leq t\leq T}x^{\vartheta}(t)\Big]\geq \mathbb{E}x^{\vartheta}(t)\geq \frac{1}{C_T}>0.
\end{align}
Using techniques in the proofs of \cite[Lemmas 5.6 and 5.7]{Mao06}, we obtain the following lemmas.
\begin{lemma}\la{log_le*7.7}
Under the condition  of Theorem \ref{log_permanence},
$$
\lim_{t\rightarrow \infty} d_{\mathbb{L}}\big(\mathbf{P}_{t}(x_0, \ell;\cdot\times\cdot), \mathbf{P}_{t}(\bar{x}_0,\bar{\ell};\cdot\times\cdot)\big) =0
$$
uniformly in $x_0,\bar{x}_0\in \mathbb{K}_+$ and $\ell,\bar{\ell}\in \mathbb{S}$.
  \end{lemma}

\begin{lemma}\la{log_le*7.8}
Under the condition  of   Lemma \ref{log_le*7.7}, for any $(x_0,\ell)\in \mathbb{R}_+\times \mathbb{S}$, $\big\{\mathbf{P}_{t}(x_0, \ell;\cdot\times\cdot)\big\}_{t\geq 0}$ is Cauchy in the space $\mathcal{P}(\mathbb{R}_+\times  \mathbb{S})$ with metric $d_{\mathbb{L}}$.
  \end{lemma}

\begin{theorem}\la{log_th*6.1}
Under the condition  of   Lemma \ref{log_le*7.7},   SDS  \eqref{logistic1} is asymptotically stable in distribution.
\end{theorem}
{\bf  Proof.}~~ For any $x\in \mathbb{R}_+$, define $v(x):=\log(x)$. For any
  compact subset
$\mathbb{K}_+$ of $ \mathbb { R } _ { + } $, and any $(x_0,\ell)\in \mathbb{K}_+\times \mathbb{S}$, define $\bar{\mathbf{P}}_{t}(v(x_0), \ell; \mathbb{D}\times\{i\}):=\mathbf{P}_{t}(x_0, \ell;v^{-1}(\mathbb{D})\times\{i\})~\forall~\mathbb{D}\in \mathscr{B}(\mathbb{R}), i\in \mathbb{S}$. Then for any $i\in \mathbb{S}$,
$
 \bar{\mathbf{P}}_{t}(v(x_0), \ell; \cdot\times i)
$
 is the transform of
$
\mathbf{P}_{t}(x_0, \ell; \cdot\times i)
$
 corresponding to $v(\cdot)$.
 By the well-known Chebyshev inequality, it is easy to observe that Lemmas \ref{log_le_le_le2} and \ref{log:Le5.8} guarantees that for any $(v(x_0),\ell)\in \mathbb{R} \times\mathbb{S}$, the family of transition probabilities $\{\bar{\mathbf{P}}_{t}(v(x_0), \ell;\cdot\times\cdot): t\geq 0\}$ is tight on $\mathbb{R} \times\mathbb{S}$. Since $(\mathbb{R}, |\cdot|)$ is complete
and separable, the tightness of $\{\bar{\mathbf{P}}_{t}(v(x_0), \ell;\cdot\times\cdot): t\geq 0\}$ on $\mathbb{R} \times\mathbb{S}$ is equivalent to  relatively compactness (see \cite[Theorems 6.1, 6.2]{Billingsley1968}). Then any sequence $\{\bar{\mathbf{P}}_{t_n}(v(x_0), \ell;\cdot\times\cdot): n\geq 0\}$  ($t_n\rightarrow \infty$ as $n\rightarrow \infty$)
  has a weak convergent subsequence denoted by $\{\bar{\mathbf{P}}_{t_n}(v(x_0), \ell;\cdot\times\cdot): n\geq 0\}$ with
some notation abuse.  Assume its weak limit is an invariant measure $\bar{\mu}(\cdot\times\cdot)$ on $\mathbb{R} \times\mathbb{S}$.  Define $\mu(\mathbb{K}_+\times\cdot)=\bar{\mu}(v(\mathbb{K}_+)\times\cdot)~\forall~\mathbb{K}\in \mathscr{B}(\mathbb{R}_+)$. Then $\mu(\cdot\times\cdot)$ is an invariant measure on  $\mathbb{R}_+\times \mathbb{S}$, and the corresponding further subsequence of  $\{\mathbf{P}_{t_n}(x_0, \ell;\cdot\times\cdot): n\geq 0\}$ converges weakly to $\mu(\cdot\times\cdot)$ on $\mathbb{R}_+ \times\mathbb{S}$.
The following  proof is the same way as the proof of \cite[Theorem 5.43]{Mao06}, so we omit it here. \eproof

\subsection{Stability in distribution of numerical solution}
We write $(X^{x_0,\ell}_k, r^{\ell}_k)$ in lieu of $(X_k, r_k)$ to highlight the initial data $(X_0, r_0) =(x_0,  \ell)$. By  \eqref{TEM_1}, we know that   $X^{x_0,\ell}_k=\mathrm{e}^{Z^{y_0,\ell}_k}$ and $x_0=\mathrm{e}^{y_0}$,
where $Z_k^{y_0,\ell}$  denote the numerical solutions defined by \eqref{log_eq5*} with initial data $(y_0,\ell)$. Similar to that of  \cite[Theorem 6.14]{Mao06}, we can prove the following result.
\begin{lemma}
$\{(Z_k, r_k)\}_{k\geq 0}$ is a time homogeneous Markov chain.
  \end{lemma}

It is easy to observe that Theorem \ref{log:nu_permanence} guarantees that for any $(x_0,\ell)\in \mathbb{R}_+\times\mathbb{S}$, the family of transition probabilities $\{\mathbf{P}^{\Delta}_{k}(x_0, \ell;\cdot\times\cdot): k\geq 0\}$ is tight. That is, for any $\varepsilon>0$ there is a compact subset $\mathbb{K}_+=\mathbb{K}_+(\varepsilon, x_0,\ell)$ of $\mathbb{R}_+$ such that
\begin{align}\la{log_eq_nu_tight}
\mathbf{P}^{\Delta}_{k}\big(x_0, \ell;\mathbb{K}_+\times\mathbb{S}\big)\geq 1-\varepsilon~~~\forall k\geq 0.
\end{align}
To show the numerical solutions $(X^{x_0, \ell}_k, r^{\ell}_k)$ defined by \eqref{TEM_1} is asymptotically stable in distribution and admit a unique invariant measure $\mu^{\Delta}(\cdot\times\cdot)\in \mathcal{P}(\mathbb{R}_+\times \mathbb{S})$, we need the following  three lemmas, the proofs of which can be found in Appendix B.
\begin{lemma}\la{log_le*6.11}
The numerical solutions defined by \eqref{TEM_1} has the property that
 $$\lim_{k\rightarrow \infty}\mathbb{E}\big| X^{x_0,\ell}_k- X^{\bar{x}_0,\ell}_k\big|=\lim_{k\rightarrow \infty}\mathbb{E}\big| \mathrm{e}^{Z^{y_0,\ell}_k}-  \mathrm{e}^{Z^{\bar{y}_0,\ell}_k}\big|= 0$$
 uniformly in $(x_0,\bar{x}_0,\ell)\in \mathbb{K}_+\times\mathbb{K}_+\times\mathbb{S}$, for any compact subset $\mathbb{K}_+$ of $\mathbb{R}_+$.
  \end{lemma}

\begin{lemma}\la{log_le*6.12}
For any compact subset $ \mathbb{K}_+ $ of $\mathbb { R } _ { + } $,
 $$\sup_{(x_0,\ell)\in \mathbb{K}_+\times\mathbb{S}}\mathbb{E}\Big[\sup_{0\leq k\Delta \leq T} X^{2}_k \Big]\leq C_T,  ~~~~ \forall~ T\geq 0,$$
  where $X_k$ is the numerical solutions defined by \eqref{TEM_1} with the initial value $(x_0,\ell)\in \mathbb{K}_+\times\mathbb{S}$.
  \end{lemma}

\begin{lemma}\la{log_le*6.13}
Under the condition Lemma \ref{log_le_le7.6},
 for any compact subset $ \mathbb{K}_+\!$ of $\mathbb { R } _ { + } $,
$$
\sup _ { \left( x _ { 0 } , \ell \right) \in \mathbb{K}_+ \times \mathbb{S} } \mathbb{E} \left[ \sup _ { 0 \leq k\Delta \leq T } X_k^{-\vartheta}   \right] \leq C_T, ~~~~~\forall~T\geq 0,
$$
 where $X_k$ is the numerical solutions defined by \eqref{TEM_1} with  the  initial value $(x_0,\ell)\in \mathbb{K}_+\times\mathbb{S}$.
\end{lemma}
By Lemma \ref{log_le*6.13} and Jensens's inequality we can get
\begin{align}\la{log_eq*eq6.32}
\mathbb{E}\Big[\sup_{0\leq k\Delta\leq T}X^{\vartheta}_k\Big]\geq \mathbb{E}X^{\vartheta}_k\geq \frac{1}{C_T}>0.
\end{align}
Using techniques in the proofs of \cite[Lemmas 6.11, 6.12 and 6.16]{Mao06}, we obtain the following three lemmas, the proofs of which are straightforward, so are omitted.
\begin{lemma}\la{log_le*6.14}
Under the condition  of Theorem \ref{log_th*6.1}, for any $\Delta\in (0, \Delta^*\wedge \Delta_3^*)$,
$$
\lim_{k\rightarrow \infty} d_{\mathbb{L}}\big(\mathbf{P}^{\Delta}_{k}(x_0, \ell;\cdot\times\cdot), \mathbf{P}^{\Delta}_{k}(\bar{x}_0,\bar{\ell};\cdot\times\cdot)\big) =0
$$
uniformly in $x_0,\bar{x}_0\in \mathbb{K}_+$ and $\ell,\bar{\ell}\in \mathbb{S}$, for any compact subset $\mathbb{K}_+$ of $\mathbb{R}_+$.
  \end{lemma}

\begin{lemma}\la{log_le*6.15}
Under the condition  of Lemma \ref{log_le*6.14},
$\big\{\mathbf{P}^{\Delta}_{k}(x_0, \ell; \cdot\times\cdot)\big\}_{k\geq 0}$ is Cauchy in the space $\mathcal{P}(\mathbb{R}_+\times  \mathbb{S})$ with metric $d_{\mathbb{L}}$.
  \end{lemma}

\begin{lemma}\la{log_le*6.16}
Fix any $(x_0,\ell)\in \mathbb{R}_+\times\mathbb{S}$. Then for any given $T>0$ and $\varepsilon>0$, there is a $\Delta^{**}\in (0,1)$, which is sufficiently small, such that
$$
d_{\mathbb{L}}\big(\mathbf{P}^{\Delta}_{k}(x_0,  \ell ;\cdot\times\cdot), \mathbf{P}_{k\Delta}(x_0, \ell ;\cdot\times\cdot)\big)<\varepsilon
$$
provided $\Delta\in (0, \Delta^{**})$ and $k\Delta\leq T$.
\end{lemma}

\begin{theorem}\la{log_th*6.2}
Under the condition  of Theorem \ref{log_th*6.1}, for any $\Delta\in (0, \Delta^*\wedge \Delta_3^*)$,
 the  numerical solutions $(X^{x_0, \ell}_k, r^{\ell}_k)$  is asymptotically stable in distribution and admit a unique invariant measure $\mu^{\Delta}(\cdot\times\cdot)\in \mathcal{P}(\mathbb{R}_+\times \mathbb{S})$.
\end{theorem}
{\bf  Proof.}~~ For any $x\in \mathbb{R}_+$, define $v(x):=\log(x)$. For any
  compact subset
$\mathbb{K}_+$ of $ \mathbb { R } _ { + } $, and any $(x_0,\ell)\in \mathbb{K}_+\times \mathbb{S}$, define $\bar{\mathbf{P}}^{\Delta}_{k}(v(x_0), \ell; \mathbb{D}\times\{i\}):=\mathbf{P}^{\Delta}_{k}(x_0, \ell;v^{-1}(\mathbb{D})\times\{i\})~\forall~\mathbb{D}\in \mathscr{B}(\mathbb{R}), i\in \mathbb{S}$. Then for any $i\in \mathbb{S}$,
$
 \bar{\mathbf{P}}^{\Delta}_{k}(v(x_0), \ell; \cdot\times i)
$
 is the transform of
$
\mathbf{P}^{\Delta}_{k}(x_0, \ell; \cdot\times i)
$
 corresponding to $v(\cdot)$.
 By the well-known Chebyshev inequality, it is easy to observe that Lemmas \ref{L:3} and \ref{log:Le5.9} guarantees that for any $(v(x_0),\ell)\in \mathbb{R} \times\mathbb{S}$, the family of transition probabilities $\{\bar{\mathbf{P}}^{\Delta}_{k}(v(x_0), \ell;\cdot\times\cdot): k\geq 0\}$ is tight on $\mathbb{R} \times\mathbb{S}$. Since $(\mathbb{R}, |\cdot|)$ is complete
and separable, the tightness of $\{\bar{\mathbf{P}}^{\Delta}_{k}(v(x_0), \ell;\cdot\times\cdot): k\geq 0\}$ on $\mathbb{R} \times\mathbb{S}$ is equivalent to  relatively compactness. Then any sequence $\{\bar{\mathbf{P}}^{\Delta}_{k_n}(v(x_0), \ell;\cdot\times\cdot): n\geq 0\}$  ($k_n\rightarrow \infty$ as $n\rightarrow \infty$)
  has a weak convergent subsequence denoted by $\{\bar{\mathbf{P}}^{\Delta}_{k_n} (v(x_0), \ell;\cdot\times\cdot): n\geq 0\}$ with
some notation abuse.  Assume its weak limit is an invariant measure $\bar{\mu}^{\Delta}(\cdot\times\cdot)$ on $\mathbb{R} \times\mathbb{S}$.  Define $\mu^{\Delta}(\mathbb{K}_+\times\cdot)
=\bar{\mu}^{\Delta}(v(\mathbb{K}_+)\times\cdot)~\forall~\mathbb{K}_+\in \mathscr{B}(\mathbb{R}_+)$.
The following  proof is the same way as the proof of \cite[Theorem 6.19]{Mao06}, so we omit it here. \eproof

We can now show that the numerical stationary distribution will weakly converge to the stationary distribution of the exact solutions.
\begin{theorem}\la{log_th*6.3}
Under the condition  of Theorem \ref{log_th*6.1},
$$
\lim_{\Delta\rightarrow 0}d_{\mathbb{L}}\big(\mu(\cdot\times\cdot), \mu^{\Delta}(\cdot\times\cdot)\big)=0.
$$
\end{theorem}
The proof of this theorem is standard (see, e.g. \cite[Theorem 6.23]{Mao06}), and hence is omitted to avoid repetition.

\section{Numerical examples}\la{section numerical}
In order to illustrate the efficiency of numerical schemes we consider a number of examples and present some simulations.
 First, we will show that the classical EM method will not be able to reproduce the  dynamical properties of the SDS \eqref{logistic1}.
 To show this, recall that the classical EM method applied to \eqref{logistic1} produces
\begin{align}\la{EM_sds1}
\left\{
\begin{array}{ll}
X^{\Delta}_0=  x_0,&\\
X^{\Delta}_{k+1}=X^{\Delta}_{k}+X^{\Delta}_{k}\Big[\big(b(r_k)-a(r_k) X^{\Delta}_k\big)\Delta +\sigma(r_k) \Delta B_k\Big].&
\end{array}
\right.
\end{align}
We choose a number $\bar{\Delta}\in (0, (0.4\mathrm{e}-1)/\check{b}]$, the following lemma shows that  for any given stepsize $\Delta\in (0, \bar{\Delta}]$ and any initial value $(x_0,\ell)\in \mathbb{R}_+\times\mathbb{S}$,  the  numerical solutions $\{|X^{\Delta}_{k}|\}_{k\geq 1}$
will tend to infinity super-exponentially with a positive probability.

\begin{lemma}\la{appendis_L1}
Let $\{X^{\Delta}_{k}\}_{k\geq 1}$ be defined by \eqref{EM_sds1}. Suppose for any $i\in \mathbb{S}$,
$a(i)-0.5|\sigma(i)|\geq1.4$ and $\sigma(i)\neq 0$. Then the conditional probability
\begin{align*}
\mathbb{P}\bigg(|X^{\Delta}_{k+1}|\geq  \frac{\exp(2^{k})}{\Delta},~\forall~k\geq 1\bigg||X^{\Delta}_{1}|\geq  \frac{\mathrm{e}}{\Delta}\bigg)
\geq\exp\Big(-
\frac{4\mathrm{e}^{-2^{-3}}}{\mathrm{e}^{4}-1}\Big).
\end{align*}
\end{lemma}
The proof of this lemma can be found in Appendix C. It then follows from Lemma \ref{appendis_L1} and \eqref{app_eq0} that
\begin{align*}
&\mathbb{P}\Big(|X^{\Delta}_{k}|\geq  \frac{\exp(2^{k-1})}{\Delta},~\forall~k\geq 1\Big)\nn\\ =&\mathbb{P}\Big(|X^{\Delta}_{1}|\geq  \frac{\mathrm{e}}{\Delta}\Big)\mathbb{P}\bigg(|X^{\Delta}_{k+1}|\geq\frac{\exp(2^{k})}{\Delta},~\forall~k\geq 1\bigg||X^{\Delta}_{1}|\geq  \frac{\mathrm{e}}{\Delta}\bigg)>0.
\end{align*}
In other words, $\{|X^{\Delta}_{k}|\}_{k\geq 1}$ will tend to infinity faster than ${\exp(2^{k-1})}/{\Delta}$ with a positive probability.
 However, our theory established in the previous sections shows that the scheme \eqref{TEM_1} can reproduce the dynamical properties of the SDS \eqref{logistic1} very well.

To illustrate our theory, as well as to compare to the simulations of the classical EM method, we shall illustrate these conclusions through the following examples.
\begin{expl}\la{exp1}
{\rm
In this example   we consider SDS \eqref{logistic1} with the Markov chain  $r(t)$ is on the state space $\mathbb{S}=\{1,2\}$ with the generator
$$
\Gamma=\left(
  \begin{array}{ccc}
    -8 & 8\\
    2 & -2\\
  \end{array}
\right),
$$
and the coefficients in each state
are given in Table \ref{T1}.  By solving the linear equation \eqref{eq:a1.2} we obtain the unique stationary (probability) distribution
$\pi=(\pi_1,\pi_2)=(0.2,0.8)$.
\renewcommand\arraystretch{0.8}
\begin{table}[!htbp]
\centering
\begin{tabular}{|c|c|c|c|c|}
\hline
\diagbox[width=7em,trim=r]{\small{States}}{\small{Coefficients}}
&~~~~~~\small{$b(i)$}~~~~~~&~~~~~~\small{$a(i)$}~~~~~~&~~~~~~\small{$\sigma(i)$}
~~~~~~&\small{$\beta(i)=b(i)-0.5\sigma^2(i)$}\\
\hline
$i=1$&2&1.8&0.8&1.68\\
\hline
$i=2$&1&2.5&2&-1\\
\hline
\end{tabular}
  \caption{Values of the coefficients in Example \ref{exp1}}\la{T1}
\vspace{-2em}\end{table}

Compute
$$
\pi a>0,~~~~~~~~\pi \beta=\pi_1\beta(1)+\pi_2\beta(2)=-0.4640.
$$
Therefore, by Theorem  \ref{log_Th*5.3},
\begin{align*}
\limsup_{t\rightarrow \infty}\frac{\log x(t)}{t}\leq\pi \beta<0~~~~a.s.
\end{align*}
In other words, this says that the extinction of the population happens.  However, by virtue of Lemma \ref{appendis_L1}, for any given stepsize $\Delta\in (0, 0.04]$ and initial value $(x_0, \ell)=(25, 1)$, one observes that $\{|X^{\Delta}_{k}|\}_{k\geq 1}$ will tend to
infinity super-exponentially with a positive probability, see Figs. \ref{exp1_1} and \ref{exp1_2}. Both simulations show that the classical EM method does not capture the dynamic properties of the underlying SDE \eqref{logistic1}, while the second simulation shows that the classical EM method can blow up very quickly.  This contrasts with the extinction of the underlying  SDS \eqref{logistic1}.

\begin{figure}[!htbp]
\begin{center}
\includegraphics[angle=0, height=6.5cm, width=14cm]{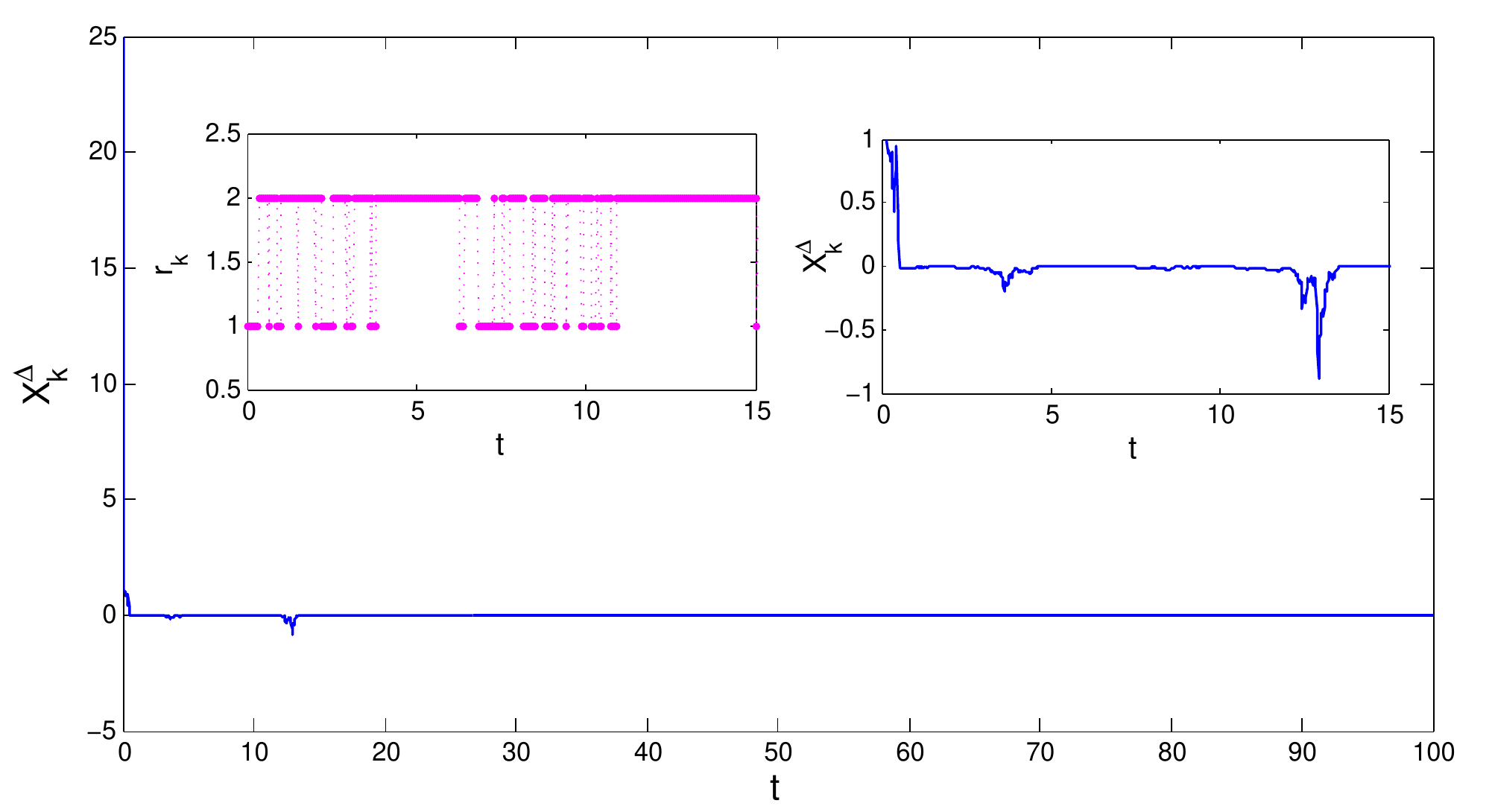}  
\caption{A sample path of  the classical EM solution $X_k^{\Delta}$ and corresponding state $r_k$. The pink trajectory represents Markov chain while the blue trajectory represents the numerical solution of Scheme \eqref{EM_sds1} with
 $\Delta=0.02$ and $t\in[0,100]$. (For interpretation of the colors in the figure(s), the reader is referred to the web version of this article.)}\label{exp1_1}
\end{center}  
\vspace{-1em}
\end{figure}
\begin{figure}[!htbp]
\begin{center}
\includegraphics[angle=0, height=6.5cm, width=14cm]{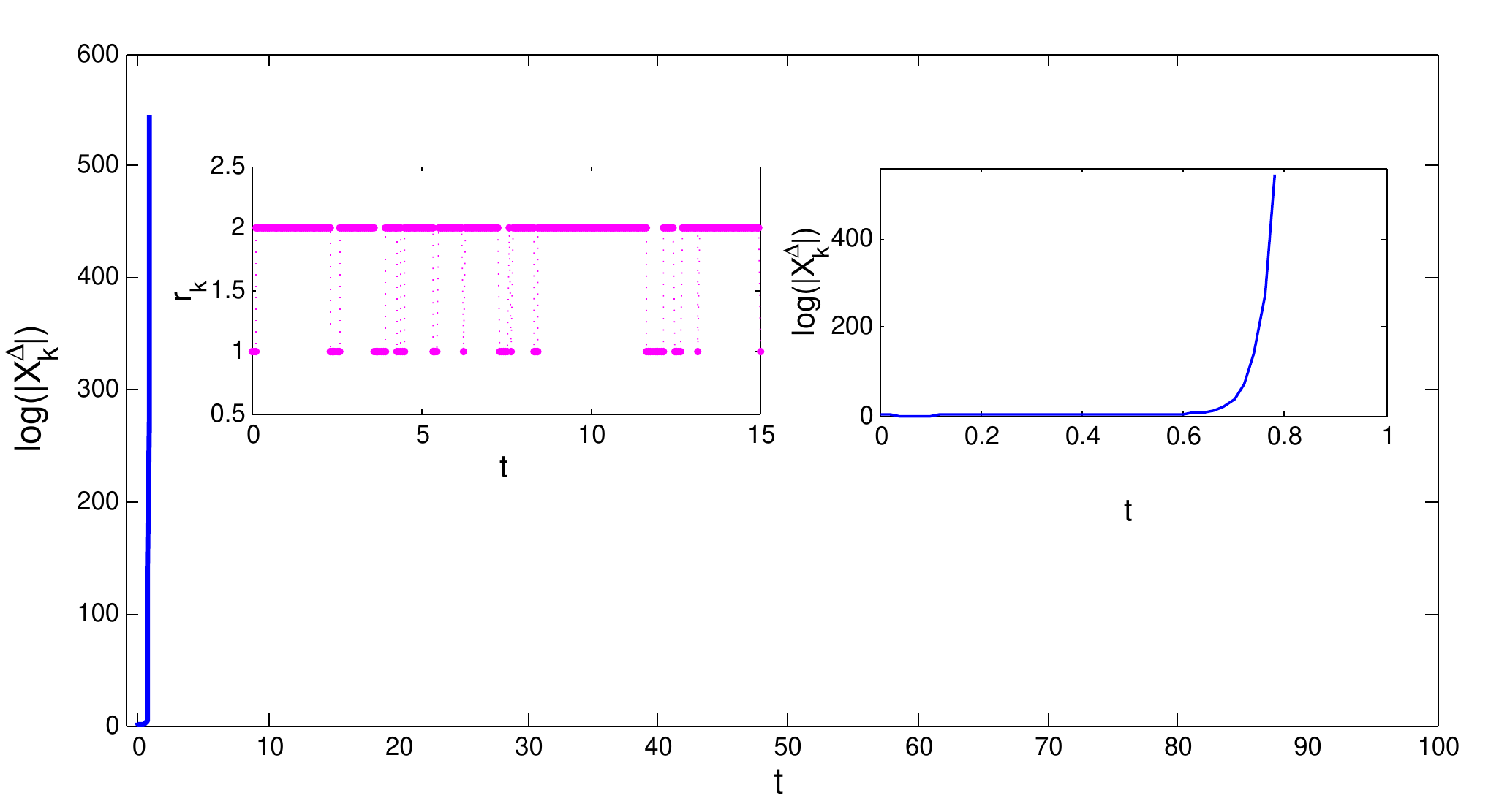}  
\caption{ A sample path of  the classical EM solution $\log(|X_k^{\Delta}|)$ and corresponding state $r_k$. The pink trajectory represents Markov chain while the blue trajectory represents the numerical solution of Scheme \eqref{EM_sds1} with
 $\Delta=0.02$ and $t\in[0,100]$.}\label{exp1_2}
\end{center}                               
\vspace{-1em}
\end{figure}

In order to represent the simulations by Scheme \eqref{TEM_1} and illustrate its effectiveness, we divide it into five steps.

\vspace*{4pt}\noindent{\bf Step~1.} {M{\scriptsize ATLAB}} code. Next we specify the {M{\scriptsize ATLAB}} code for simulating $X^E_k$ and $X_k$:
\vspace{-1.2em}
\begin{lstlisting}
%MATLAB code for simulating X^E_k and X_k
clear all;
X(1)=25; XE(1)=25; r=1; T=100; Y=log(XE(1)); Z=log(X(1));
b=[2,1]; a=[1.8,2.5]; cgm=[0.8,2]; beat=b-0.5*cgm.^2;
Gam=[-8 8;2 -2]; dt=0.02; dB=sqrt(dt)*randn(1,T/dt);
c=expm(Gam*dt); v=log(10*dt^(-2/5)); %Obviously, v>Z;
for n=1:T/dt
    Y=Y+(beat(r)-a(r)*exp(Y))*dt+cgm(r)*dB(n);
    Z=Z+(beat(r)-a(r)*exp(Z))*dt+cgm(r)*dB(n);
    if rand<c(r,1)
        r=1;
    else
        r=2;
    end
    XE(n+1)=exp(Y);
    if Z>v
        Z=v;
    end
    X(n+1)=exp(Z);
end
\end{lstlisting}
\vspace*{4pt}\noindent{\bf Step~2.} Approximating the error
$\mathbb{E}|x(T)-X_{\Delta}(T)|^{p}$.
To compute the approximation error, we run $M$ independent trajectories where $x^{(j)}(t)$ and $X^{(j)}_{\Delta}(t)$ represent the $j$th trajectories of exact solution $x(t)$ and the numerical solution $X_{\Delta}(t)$ respectively. Thus
$$
\mathbb{E}|x(T)-X_{\Delta}(T)|^{p}=\frac{1}{M}\sum_{j=1}^{M}|x^{(j)}(T)-X_{\Delta}^{(j)}(T)|^p.
$$
\vspace*{4pt}\noindent{\bf Step~3.} The log-log error plot with $M=2000$. The simulation procedure is carried out by steps 1 and 2.
 The red dashed line depicts log-log error while the blue solid line is a reference line of slope $1/2$ in Fig. \ref{exp_1fig1}. Fig. \ref{exp_1fig1} depicts the approximation error $
\mathbb{E}|x(32)-X^E_{\Delta}(32)|
$ between the exact solution of the SDS \eqref{logistic1} and the numerical solution by Scheme \eqref{TEM2}, and the error $
\mathbb{E}|x(32)-X_{\Delta}(32)|
$ between the exact solution and that by  Scheme \eqref{TEM_1}  with $\theta=0.5$,  as the function of stepsize $\Delta\in\{2^{-1},2^{-2},\ldots,2^{-14}\}$.
\begin{figure}[!htp]
\begin{center}
\includegraphics[angle=0, height=7cm, width=16cm]{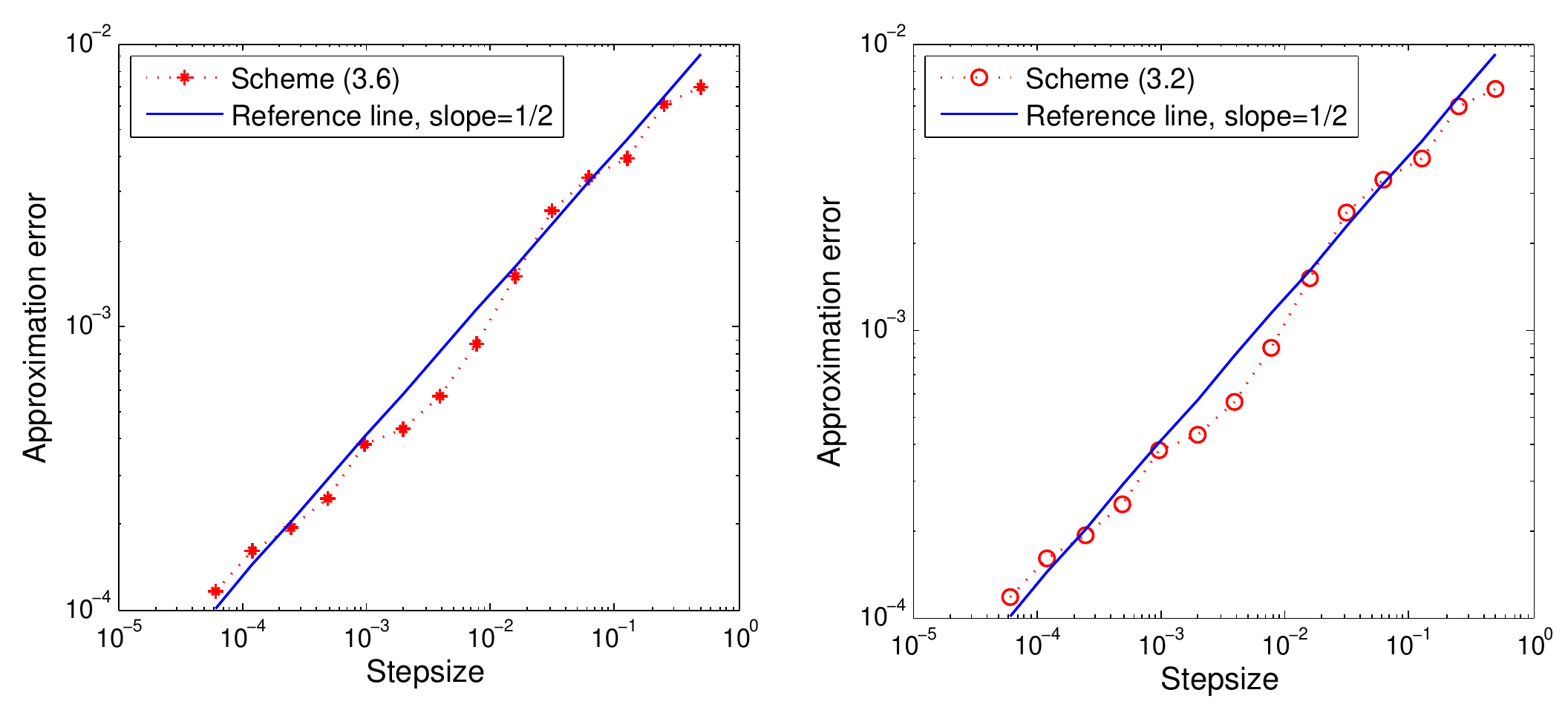}  
\caption{The red asterisk trajectory depicts the approximation error $\mathbb{E}|x(32)-X^E_{\Delta}(32)|$ of the exact solution of SDS \eqref{logistic1} and the numerical solution by Scheme \eqref{TEM2} while the red circle trajectory depicts the  approximation error  $\mathbb{E}|x(32)-X_{\Delta}(32)|$ of the exact solution of SDS \eqref{logistic1} and the numerical solution by Scheme \eqref{TEM_1} as the functions of    stepsize
$\Delta\in\{2^{-1},2^{-2},\ldots,2^{-14}\}$.}
\label{exp_1fig1}                                 
\end{center}                               
\vspace{-3em}
\end{figure}

\vspace*{4pt}\noindent{\bf Step~4.}
To compare to the simulations of the classical EM method.  The simulation procedure is carried out by
step  1, and all parameters are  same as the classical EM method. The two simulations shown in Figs. \ref{exp1_3} and \ref{exp1_4} are based on Scheme \eqref{TEM_1}. Both figures show clearly that Scheme \eqref{TEM_1}  reproduces the dynamic properties of the underlying SDS \eqref{logistic1}.

\begin{figure}[!htbp]
\begin{center}
\includegraphics[angle=0, height=6cm, width=12cm]{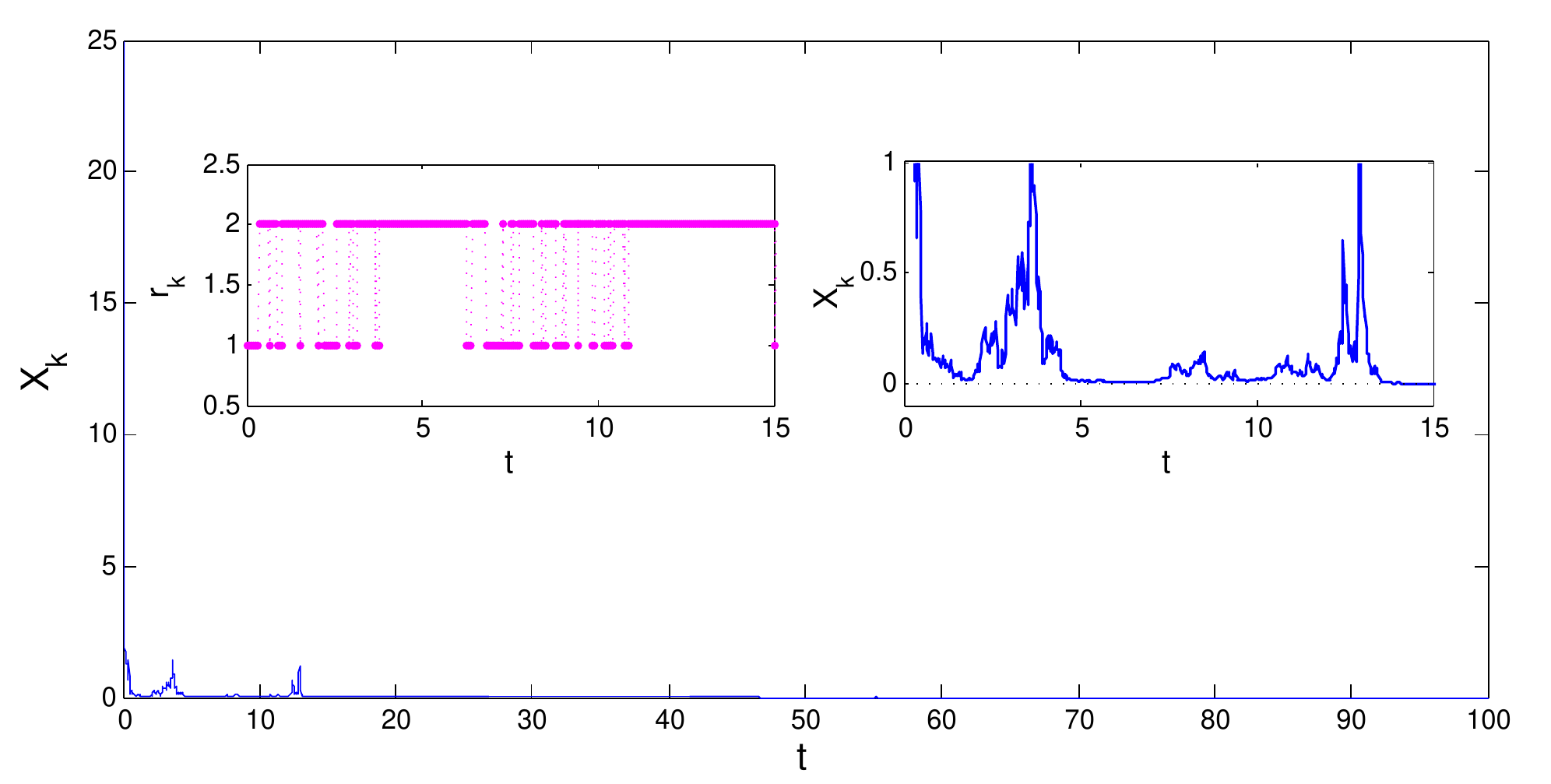}  
\caption{A sample path of  the numerical solutions $X_k$ and corresponding state $r_k$. The pink trajectory represents Markov chain while the blue trajectory represents the numerical solution of Scheme \eqref{TEM_1} with
 $\Delta=0.02$, $\theta=0.4$ and $t\in[0,100]$.}\label{exp1_3}
\end{center}                               
\vspace{-1em}
\end{figure}
\begin{figure}[tb]
\begin{center}
\includegraphics[angle=0, height=6cm, width=12cm]{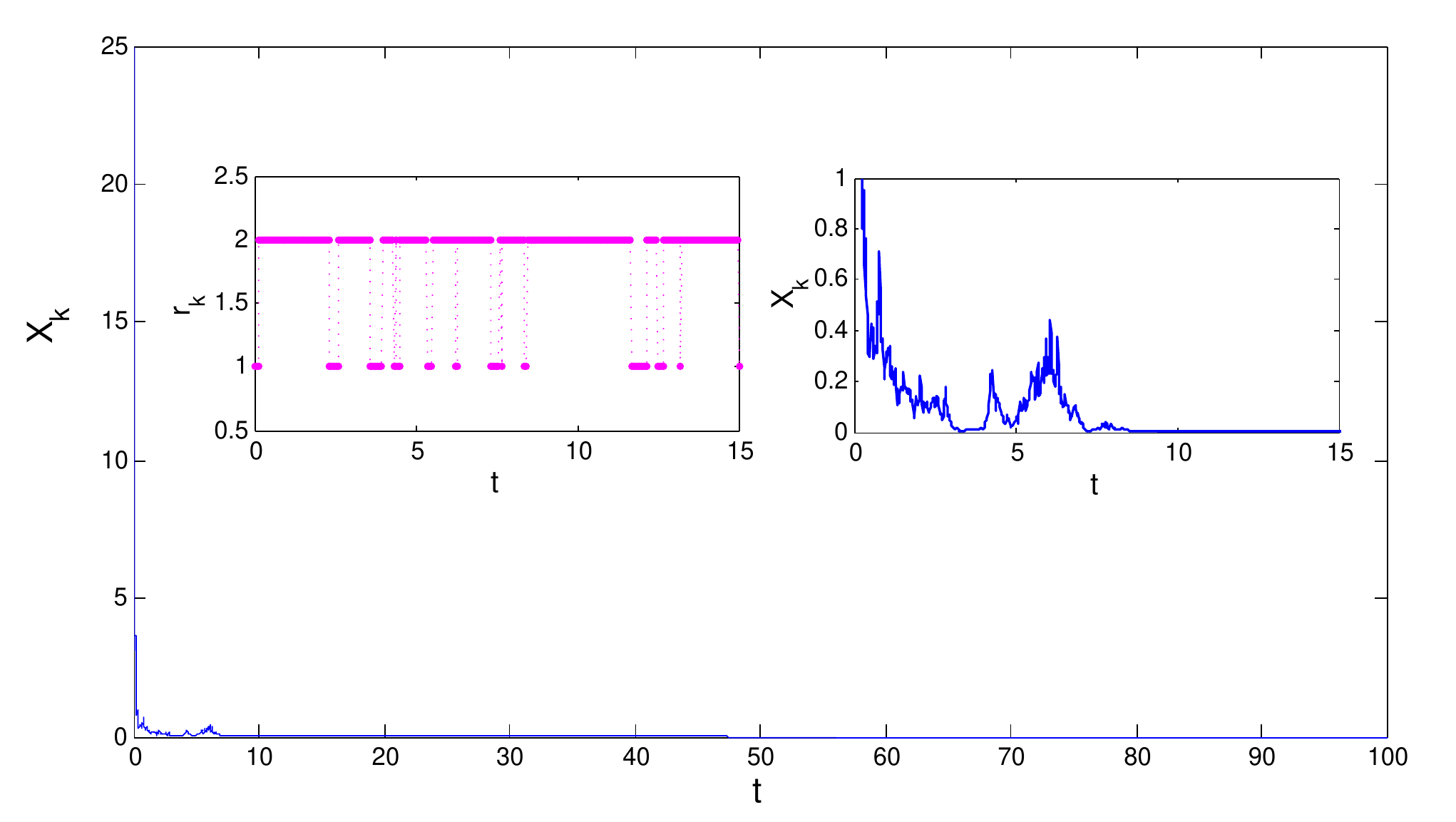}  
\caption{A sample path of  the numerical solutions $X_k$ and corresponding state $r_k$. The pink trajectory represents Markov chain while the blue trajectory represents the numerical solution of Scheme \eqref{TEM_1} with
 $\Delta=0.02$, $\theta=0.4$ and $t\in[0,100]$.}\label{exp1_4}
\end{center}                              
 \vspace{-1em}
\end{figure}

\noindent{\bf Step~5.}  Further  show that Scheme \eqref{TEM_1}  can reproduce this extinction very well.
Fig. \ref{exp1_5} depicts $500$  sample paths of the numerical solution of Scheme \eqref{TEM_1}.  This figures show clearly that Scheme \eqref{TEM_1}  reproduces positivity and extinction of the underlying SDS \eqref{logistic1}  (see the enlargement in Fig. \ref{exp1_5}).

\begin{figure}[!thbp]
\begin{center}
\includegraphics[angle=0, height=7cm, width=12cm]{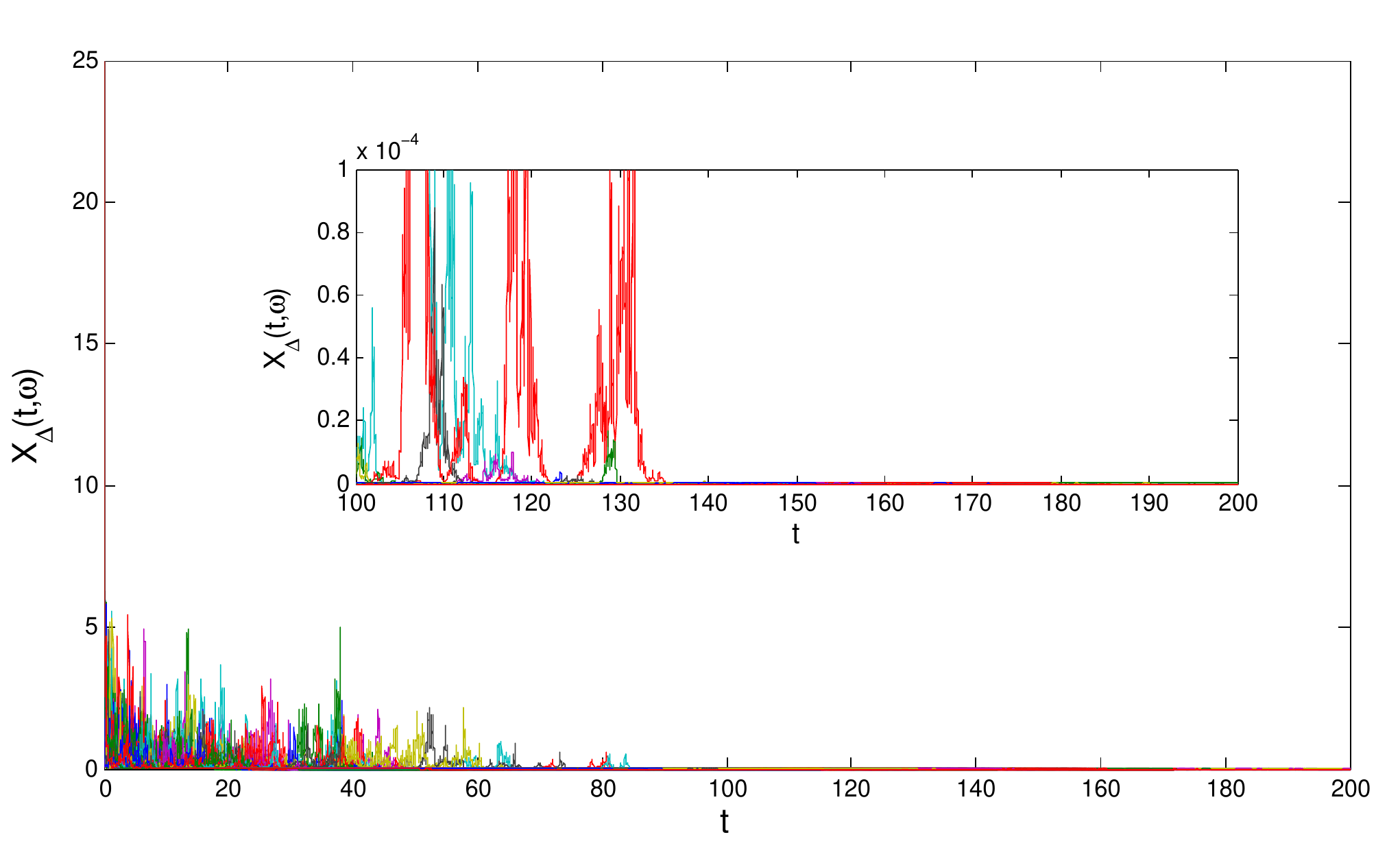}  
\caption{500  trajectories of the numerical solution of Scheme \eqref{TEM_1} with
 $\Delta=0.02$, $\theta=0.4$ and $t\in[0,200]$.}\label{exp1_5}
\end{center}                              
\vspace{-1em}
\end{figure}
}
 \end{expl}

Before closing this section we carry out some simulations to illustrate the efficiency of Scheme \eqref{TEM_1} in the approximation of invariant measures.

\begin{expl}\la{exp2}
{\rm
In this example   we consider SDS \eqref{logistic1} with the Markov chain  $r(t)$ is on the state space $\mathbb{S}=\{1,2,3\}$ with the generator
$$
\Gamma=\left(
  \begin{array}{ccc}
    -10 & 0 & 10\\
    2 & -2 & 0\\
    0 & 1 & -1
  \end{array}
\right),
$$
and the coefficients in each state
are given in Table \ref{T2}.  By solving the linear equation \eqref{eq:a1.2} we obtain the unique stationary (probability) distribution
$\pi=(\pi_1,\pi_2,\pi_3)=(\frac{1}{16},\frac{5}{16},\frac{10}{16})$.
\renewcommand\arraystretch{0.8}
\begin{table}[!htbp]
\centering
\begin{tabular}{|c|c|c|c|c|}
\hline
\diagbox[width=7em,trim=r]{\small{States}}{\small{Coefficients}}
&~~~~~~\small{$b(i)$}~~~~~~&~~~~~~\small{$a(i)$}~~~~~~&~~~~~~\small{$\sigma(i)$}
~~~~~~&\small{$\beta(i)=b(i)-0.5\sigma^2(i)$}\\
\hline
$i=1$&0.7&0.3&$\sqrt{3}$&-0.8\\
\hline
$i=2$&0.4&0.8&0.06&0.3982\\
\hline
$i=3$&1&0.5&0.04&0.9992\\
\hline
\end{tabular}
  \caption{Values of the coefficients in Example \ref{exp2}}\la{T2}\vspace{-0.5em}
\end{table}

Compute
$$
\pi a=0.5813>0,~~~~\pi\beta=0.6989>0.
$$
Therefore, by  Theorems \ref{log_Th*5.4} and \ref{log_th*6.1}, SDS \eqref{logistic1} is stochastically permanent and  asymptotically stable in distribution, namely the probability measure $\mathbf{P}_{t}(x_0,  \ell ;\cdot\times\cdot)$ of the solution
$x(t)$ tends to an invariant measure $\mu(\cdot\times\cdot)$ asymptotically as $t\rightarrow \infty$. On the other hand, by virtue of Theorem \ref{log:nu_permanence}, the numerical solutions $X_{k}$ are  stochastically permanent. Meanwhile,   by Theorems \ref{log_th*6.2}
and \ref{log_th*6.3}, the probability measure $\mathbf{P}^{\Delta}_{k}(x_0,  \ell ;\cdot\times\cdot)$ of the solution using Scheme \eqref{TEM_1} with any initial value
$(x_0,\ell)\in\mathbb{R}_+\times\mathbb{S}$ tends to a unique numerical invariant measure $\mu^{\Delta}(\cdot\times\cdot)$ asymptotically as $k\rightarrow \infty$, and  $\mu^{\Delta}(\cdot\times\cdot)\rightarrow\mu(\cdot\times\cdot)$ as $\Delta\rightarrow 0$.
\begin{figure}[!htbp]
\begin{center}
\includegraphics[angle=0, height=5.8cm, width=16cm]{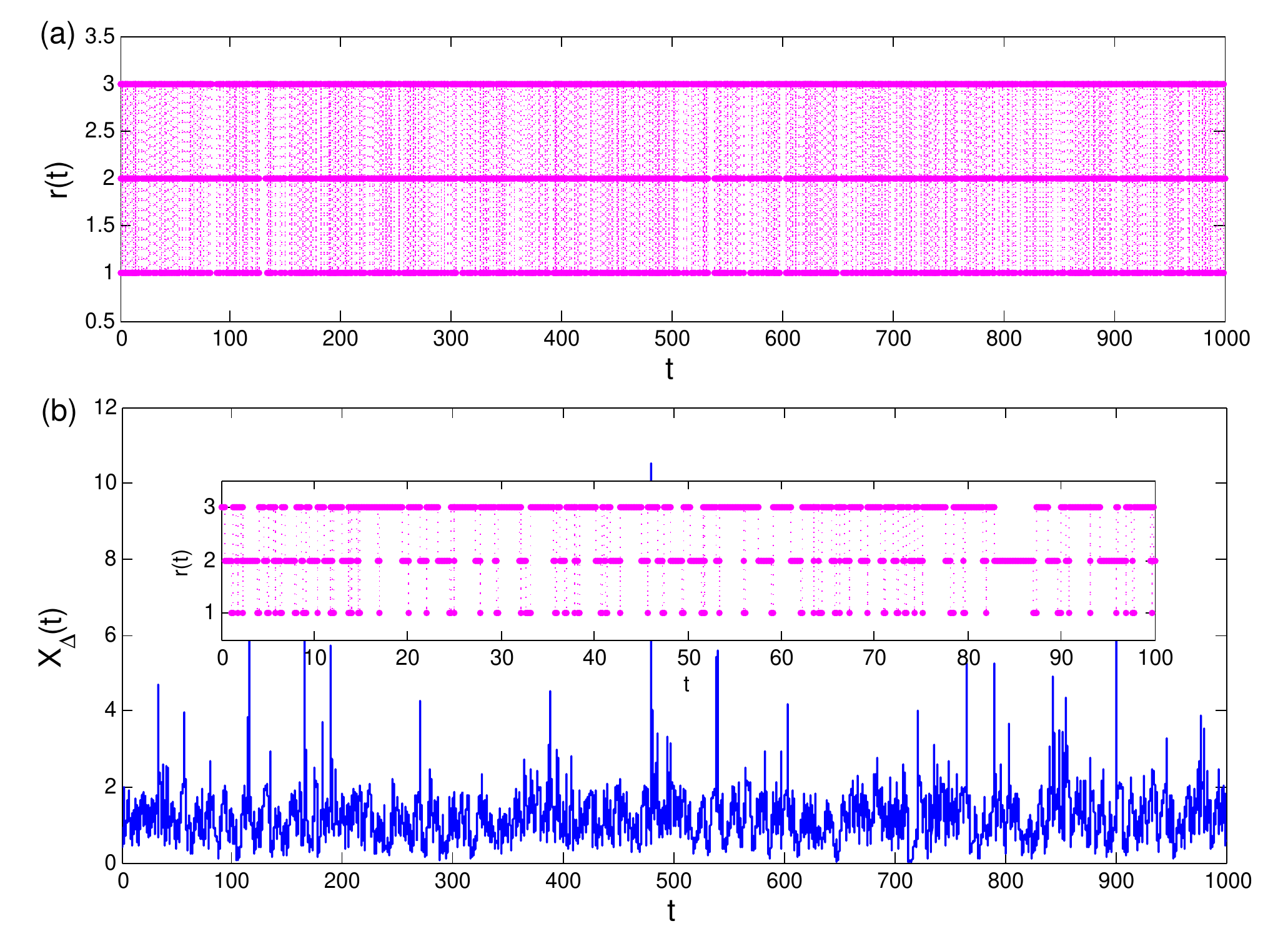}  
\caption{(a) Computer simulation of a sample path of Markov chain $r(t)$. (b) A sample path  of  numerical solution of Scheme \eqref{TEM_1} with
 $\Delta=10^{-2}$ and $\theta=0.4$ (the blue solid line).}\label{permanent_markov3}
\end{center}                              
\vspace{-1.5em}\end{figure}
\begin{figure}[!htbp]
\begin{center}
\includegraphics[angle=0, height=6cm, width=12cm]{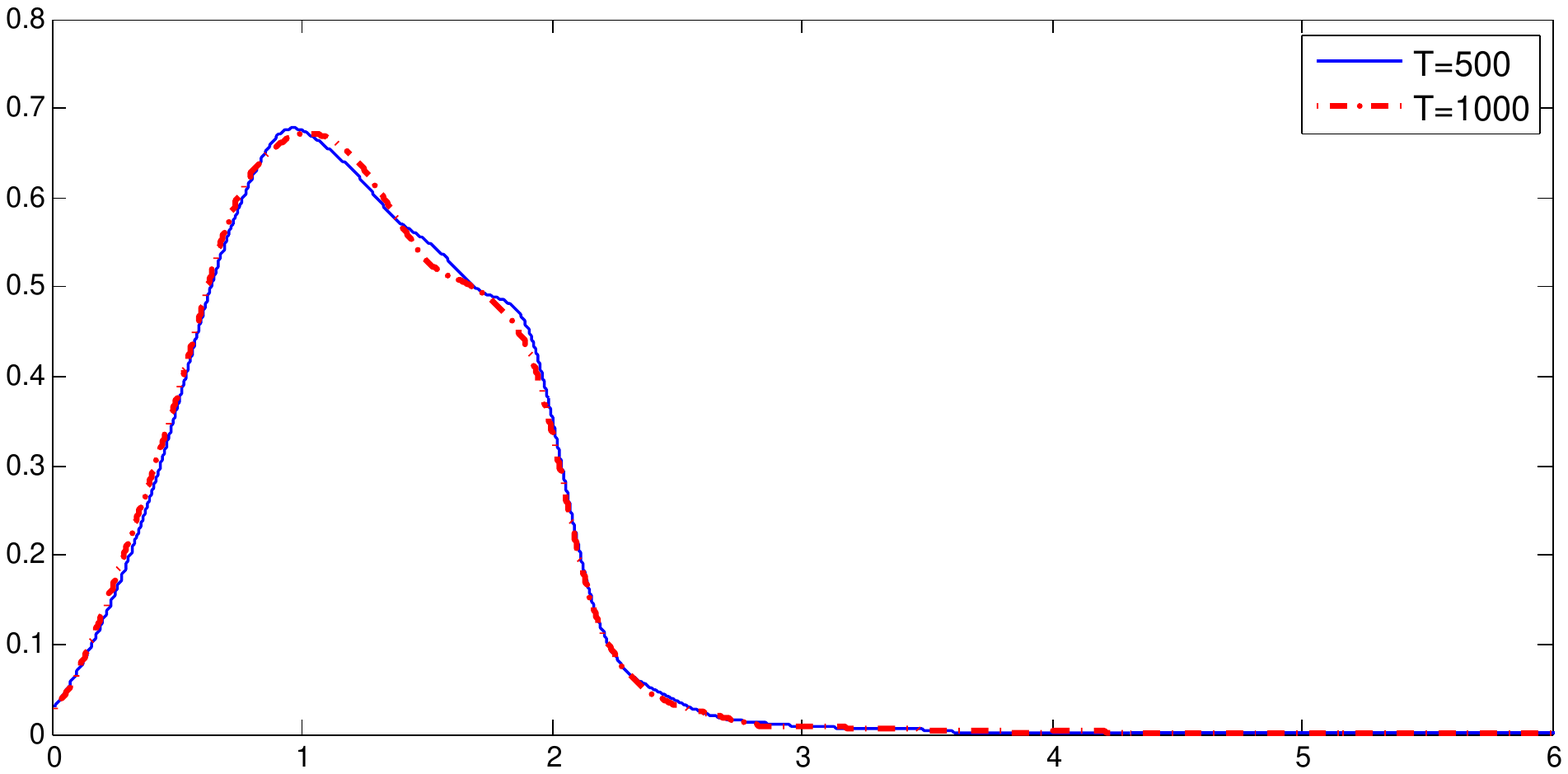}  
\caption{The empirical density of $\mu^{\Delta}$ for $10000$ sample points with
$T=500$ and $T=1000$, respectively.}\label{density}
\end{center}                              
\vspace{-1em}\end{figure}

Next, in order to test the efficiency of the scheme, we carry out numerical experiments by implementing Scheme \eqref{TEM_1} using {M{\scriptsize ATLAB}}.
 Let $(x_0,\ell)=(0.5,3)$ and take $\Delta=10^{-2}$, $K=10$ and $\theta=0.4$.
Fig. \ref{permanent_markov3} (a) depicts the path of the Markov chain while Fig. \ref{permanent_markov3} (b) further compares the path of   the numerical solution $X_{\Delta}(t)$.
 Fig.  \ref{density} depicts the empirical density of $\mu^{\Delta}$, which predicts the stationary distribution.

}
 \end{expl}

\section{First order strong convergence}\la{section_7}
In this section  we focus on the numerical approximation of  the stochastic logistic population system without regime switching (i.e.  $m=1$), and we will show that construct the explicit schemes with strongly converges with rate one. For $m=1$, we may consider without loss of generality that $\mathbb{S} =\{i\}$, $b(i)\equiv b, a(i)\equiv a$ and $\sigma(i)\equiv \sigma$, as a special case, the  stochastic logistic population system
\be\label{SDE_logistic1}
\mathrm{d}x(t)=x(t)\Big[\big(b-ax(t)\big)\mathrm{d}t + \sigma\mathrm{d}B(t)\Big]
\ee
called the subsystem of the SDS \eqref{logistic1} is permanent or not depends only on its parameters.
\begin{remark}
For the special case $\mathbb{S} =\{i\}$, the  all results on numerical solutions in Sections  \ref{section3}, \ref{strong2} and \ref{sections5} still hold.
  \end{remark}
For the subsystem \eqref{SDE_logistic1},
the EM method   \eqref{log_eq5} degenerate into  the following form
\begin{align}\la{SDE_log_eq5}
\left\{
\begin{array}{ll}
Y_0=\log x_0,&\\
Y_{k+1}=Y_{k}+\big(\beta -a  \mathrm{e}^{Y_k}\big)\Delta +\sigma  \Delta B_k,&
\end{array}
\right.
\end{align}
for any integer $k\geq 0$  and any  $\Delta\in (0, 1]$.
To study the
rate of convergence of numerical solutions $\{X^E_k\}_{k\geq 0}$ and $\{X_k\}_{k\geq 0}$ of the subsystem \eqref{SDE_logistic1}, we need the following lemma, the proof of which can be found in Appendix D.
\begin{lemma}\la{SDE_log_th2}
For any $q>0$ there exists a constant $C_T>0$  such that
\begin{align*}
\mathbb{E}\bigg[\sup_{k=0,\ldots, \lfloor T/\Delta\rfloor} |Y_{k}-y(t_{k})|^{q}\bigg]
\leq& C_T\Delta^{q}
\end{align*}
for any $\Delta\in(0, 1]$ and $T>0$.
\end{lemma}

\begin{theorem}\la{SDE_log_th3}
For any $p>0$ there exists a constant $C_T>0$   such that
\begin{align*}
 \sup_{k=0,\ldots, \lfloor T/\Delta\rfloor}\mathbb{E}\bigg[ |x(t_k)-X^E_k|^{p}\bigg]
\leq& C_T\Delta^{p}
\end{align*}
for any $\Delta\in(0, 1]$ and $T>0$.
\end{theorem}
{\bf  Proof.}~~   Using H\"{o}lder's inequality,  for any $\delta>1$, we have
\begin{align*}
 \mathbb{E}\bigg[  |x(t_k)-X^E_k|^{p}\bigg]
\leq&\bigg[\mathbb{E} \big(\mathrm{e}^{py(t_k)}+ \mathrm{e}^{pY_k}\big)^{(\delta-1)/\delta} \bigg]^{\delta/(\delta-1)} \bigg[\mathbb{E}\sup_{k=0,\ldots, \lfloor T/\Delta\rfloor} |y(t_k)-Y_k|^{p\delta}\bigg]^{1/\delta}.
\end{align*}
Thus, by applying Lemma \ref{SDE_log_th2}, Corollary \ref{SDE_log_le4}   and \eqref{log_eq3}, we infer that
\begin{align*}
\mathbb{E}\bigg[  |x(t_k)-\mathrm{e}^{Y_k}|^{p}\bigg] \leq C_T\Delta^{p}
\end{align*}
for any $\Delta\in(0, 1]$. The proof is complete. \eproof

\begin{theorem}\la{SDE_log_th4}
If $a>0$, for any $p>0$ there exists a constant $C_T>0$  such that
\begin{align*}
\mathbb{E}\bigg[\sup_{k=0,\ldots, \lfloor T/\Delta\rfloor} |x(t_k)-X^E_k|^{p}\bigg]
\leq& C_T\Delta^{p}
\end{align*}
for any $\Delta\in (0, \Delta^*)$ and $T>0$.
\end{theorem}
{\bf  Proof.}~~Now the mean value theorem implies
\begin{align*}
 \mathbb{E}\bigg[\sup_{k=0,\ldots, \lfloor T/\Delta\rfloor} |X(t_k)-\mathrm{e}^{Y_k}|^{p}\bigg]
\leq& \mathbb{E}\bigg[\sup_{k=0,\ldots, \lfloor T/\Delta\rfloor}\big(\mathrm{e}^{py(t_k)}+ \mathrm{e}^{pY_k}\big) |y(t_k)-Y_k|^{p}\bigg].
\end{align*}
 Using H\"{o}lder's inequality,  for any   $\frac{ p-2 }{p}\vee 0<\frac{1}{\delta}<1$, we have
\begin{align*}
&\mathbb{E}\bigg[\sup_{k=0,\ldots, \lfloor T/\Delta\rfloor} |x(t_k)-\mathrm{e}^{Y_k}|^{p}\bigg]\nn\\
\leq&\bigg[\mathbb{E}\sup_{k=0,\ldots, \lfloor T/\Delta\rfloor}\big(\mathrm{e}^{py(t_k)}+ \mathrm{e}^{pY_k}\big)^{(\delta-1)/\delta} \bigg]^{\delta/(\delta-1)} \bigg[\mathbb{E}\sup_{k=0,\ldots, \lfloor T/\Delta\rfloor} |y(t_k)-Y_k|^{p\delta}\bigg]^{1/\delta}.
\end{align*}
Thus, by applying Lemmas \ref{SDE_log_th2}, \ref{log_le*7.5}   and \ref{log_le*6.12} as well as \eqref{log_eq3}, we infer that
\begin{align*}
\mathbb{E}\bigg[\sup_{k=0,\ldots, \lfloor T/\Delta\rfloor} |x(t_k)-\mathrm{e}^{Y_k}|^{p}\bigg] \leq C_T\Delta^{p}.
\end{align*}
The proof is complete. \eproof

\begin{remark}\la{re6.2}
Using the same method as employed in the proofs of Section   \ref{section3}, we can easily also obtain the first order strong convergence rate  for the approximation of the original  SDE \eqref{SDE_logistic1} by $X_k:=\mathrm{e}^{Z_k}$, and hence is omitted to avoid repetition.
  \end{remark}

\begin{remark}\la{re6.3}
By  virtue of  Theorem \ref{log_th*6.1},  we know that  SDE  \eqref{SDE_logistic1} is asymptotically stable in distribution. On the other hand, under the condition  of   Theorem \ref{log_th*6.1},
 by solving the Fokker-Planck equation (see details in \cite{Pasqual2011}),
the process $x(t)$ has a unique stationary distribution $\mu(\cdot)$, and obeys the Gamma distribution with parameter
$$
\alpha_1=\frac{2b}{\sigma^2}-1,~~~~~\alpha_2=\frac{2a}{\sigma^2},
$$
with a notation abuse slightly, we write $x\sim  Ga(\alpha_1, \alpha_2)$, with  density
$$
p(x)=\frac{(\alpha_2)^{\alpha_1}}{\bar{\Gamma}(\alpha_1)}x^{\alpha_1-1}
\mathrm{e}^{-\alpha_2 x},~~~~~~x>0,
$$
where $\bar{\Gamma}(\cdot)$ is the Gamma function. By the strong law of
large numbers we deduce that
$$
\lim_{t\rightarrow \infty}\frac{1}{t}\int_{0}^{t}x(s)\mathrm{d}s=\int_{0}^{\infty}x p(x)\mathrm{d}x:=\frac{\alpha_1}{\alpha_2}=\frac{2b-\sigma^2}{2a}~~~a.s.
$$
\end{remark}

Let us discuss an example and present some simulations to illustrate our theory before closing this section to highlight the advantages of our new results on the convergence rates.

\begin{expl}\la{exp2}
{\rm
 In this example we consider autonomous stochastic logistic model   \eqref{SDE_logistic1}  of the form
\be\label{exp_logistic1}
\mathrm{d}x(t)=x(t) \big(0.5-0.8x(t)\big)\mathrm{d}t + 0.3x(t)\mathrm{d}B(t)
\ee
with an initial value $x_0=50$. Then we have
$$
a>0,~~~~~\beta=b-0.5\sigma^2=0.455>0.
$$
Theorems \ref{log_permanence} and \ref{log_th*6.1} tell us that  $x(t)$ is stochastically permanent (see the red solid line of Fig. \ref{exp_3permanence}) and has a unique stationary distribution. Meanwhile, Remark \ref{re6.3} also shows that the distribution of  $x(t)$
weakly converges to the unique invariant probability measure $\mu(\cdot)$, the Gamma distribution with $\alpha_1=91/9$ and $\alpha_2=160/9$.

By virtue of Theorem \ref{SDE_log_th4} and Remark \ref{re6.2}, the numerical solutions $X^{E}_{\Delta}(t)$, $X_{\Delta}(t)$   approximates the exact solution in the mean square sense with error estimate $\Delta$, respectively. It follows from Theorem  \ref{log:nu_permanence} that given $a>0$ and  $\beta>0$,   the numerical solution $X_{\Delta}(t)$  is stochastically permanent (see the  blue dashed line of Fig. \ref{exp_3permanence}).  Moreover, by Theorems \ref{log_th*6.2} and \ref{log_th*6.3}, the probability measure of the solution using Scheme \eqref{TEM_1} with any initial value
$x_0\in \mathbb{R}_+$ tends to a unique numerical invariant measure $\mu^{\Delta}(\cdot)$ asymptotically as $k\rightarrow \infty$, and  $\mu^{\Delta}(\cdot)\rightarrow \mu(\cdot)$ as $\Delta\rightarrow 0$.

\begin{figure}[!htp]
\begin{center}
\includegraphics[angle=0, height=5.4cm, width=14cm]{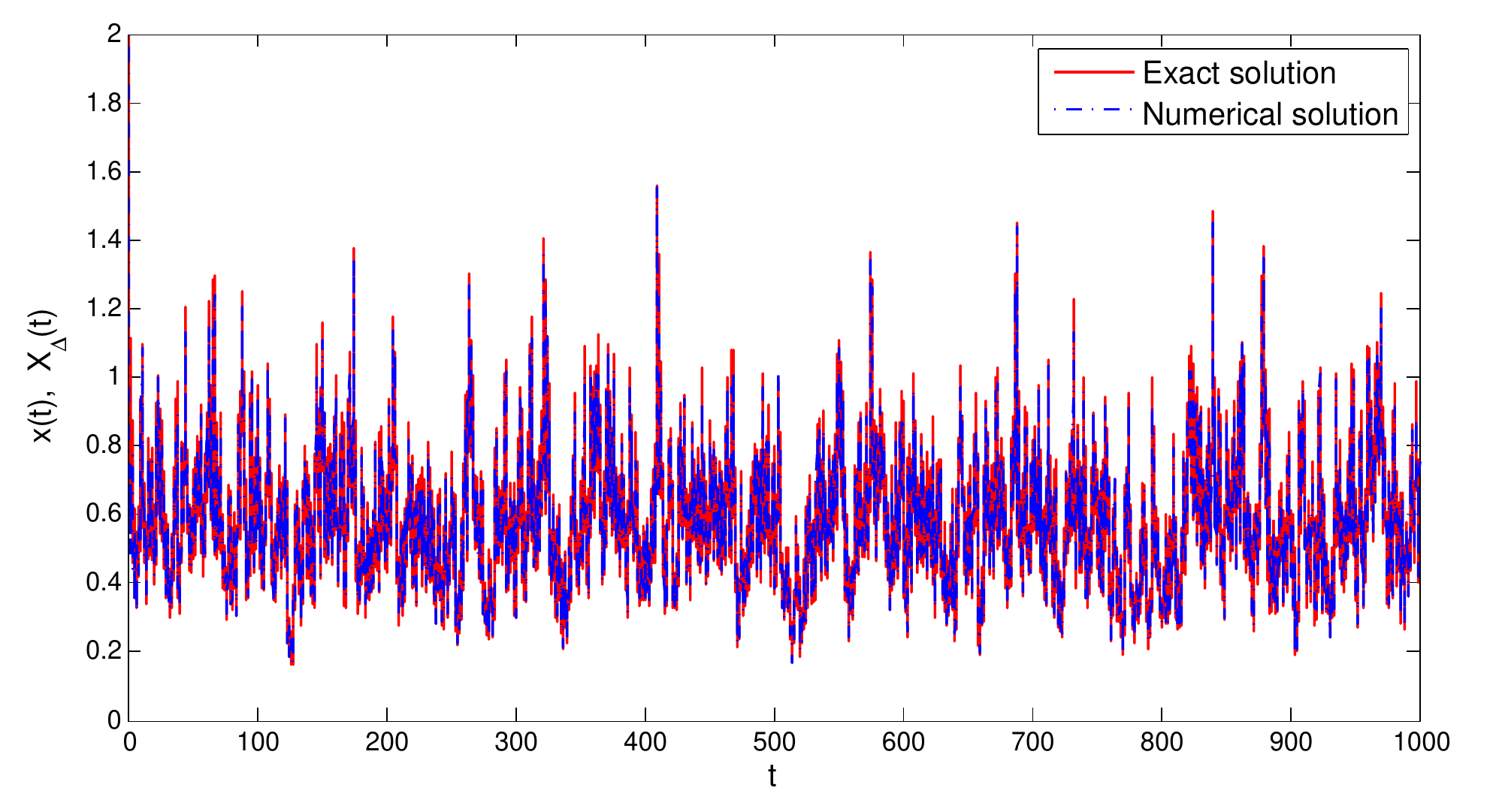}  
\caption{Sample paths of the exact solution (the red solid line) and numerical solution of Scheme \eqref{TEM_1} with
 $\Delta=0.1$, $K=20$ and $\theta=0.4$ (the blue dashed line). (For interpretation of the colors in the figure(s), the reader is referred to the web version of this article.)}                                \label{exp_3permanence}
\end{center}                               
\vspace{-2em}
\end{figure}
\begin{figure}[!htp]
\begin{center}
\includegraphics[angle=0, height=5.4cm, width=16cm]{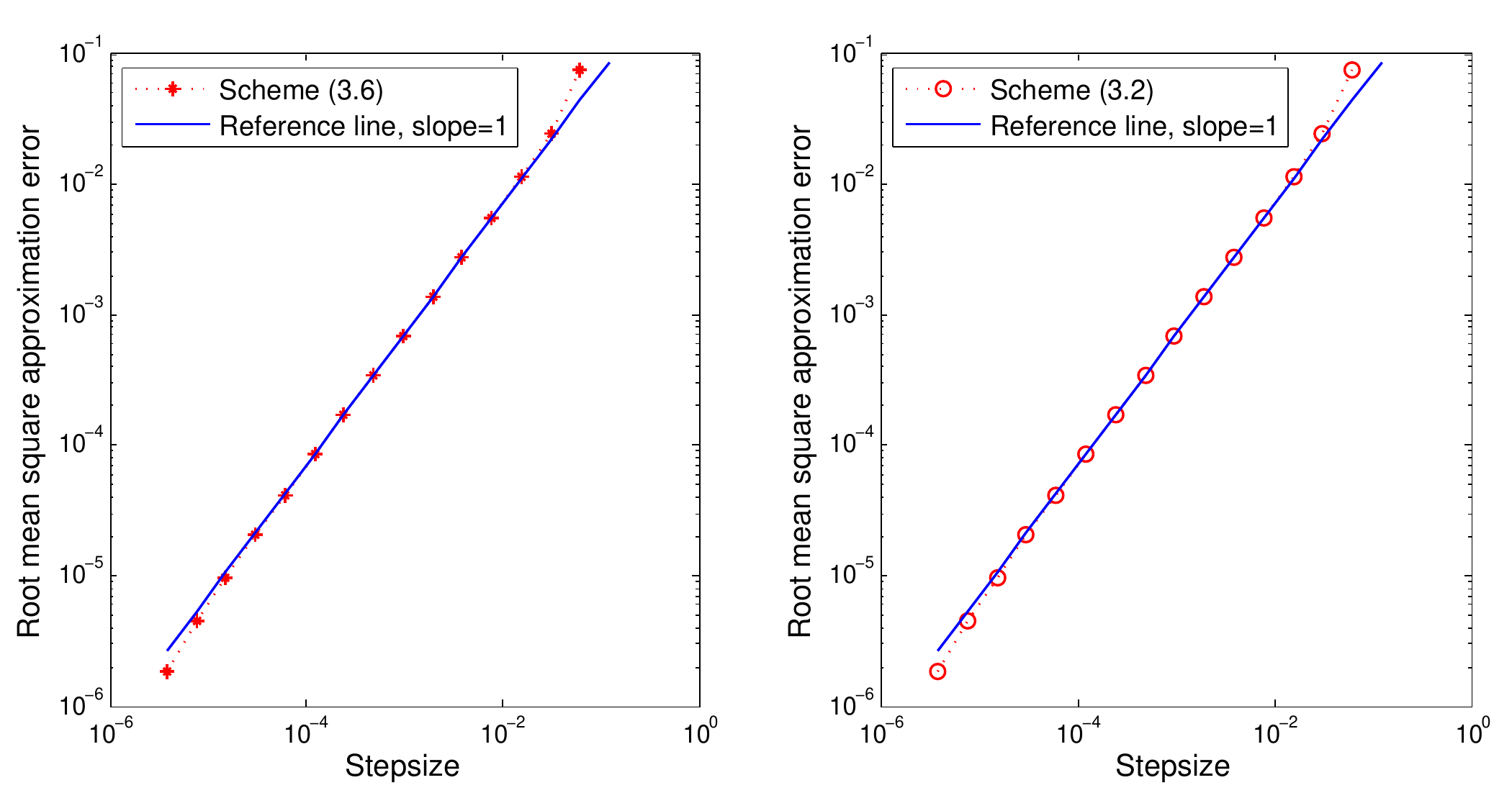}  
\caption{The red asterisk trajectory depicts the root mean square approximation error $\big(\mathbb{E}|x(2)-X^E_{\Delta}(2)|^2\big)^{1/2}$ of the exact solution of SDE \eqref{exp_logistic1} and the numerical solution by Scheme \eqref{TEM2} while the red circle trajectory depicts the root mean square approximation error $\big(\mathbb{E}|x(2)-X_{\Delta}(2)|^2\big)^{1/2}$ of the exact solution of SDE \eqref{exp_logistic1} and the numerical solution by Scheme \eqref{TEM2} as the functions
of  stepsize  $\Delta\in \{2^{-4}, 2^{-5}, \ldots, 2^{-18}\}$.}
\label{exp_2fig1}                                 
\end{center}                               
\vspace{-2em}
\end{figure}

To test the efficiency of the scheme we carry out numerical experiments by implementing Schemes \eqref{TEM_1} and \eqref{TEM2} using {M{\scriptsize ATLAB}}.
 Fig. \ref{exp_2fig1} plots the root mean square approximation error $
\big(\mathbb{E}|x(2)-X^E_{\Delta}(2)|^2\big)^{1/2}
$ between the exact solution of SDE \eqref{exp_logistic1} and the numerical solution by Scheme \eqref{TEM2}, and the
error $
\big(\mathbb{E}|x(2)-X_{\Delta}(2)|^2\big)^{1/2}
$ between the exact solution and that of Scheme \eqref{TEM_1}, as the functions
of  stepsize  $\Delta\in \{2^{-4}, 2^{-5}, \ldots, 2^{-18}\}$, for $10^4$ sample points.

One observes that the schemes proposed in \cite{Liu13,li2018jde,Li2019IMA,mao15a} are not preserve positivity  and therefore are not well defined when directly applied to SDS \eqref{logistic1}, which don't work for the above stochastic logistic models. However, the performance of Scheme \eqref{TEM_1} is very nice for this case,
 see Figs. \ref{exp_3permanence} and \ref{exp_3EM_S}.
\begin{figure}[!htp]
\begin{center}
\includegraphics[angle=0, height=5.5cm, width=16.5cm]{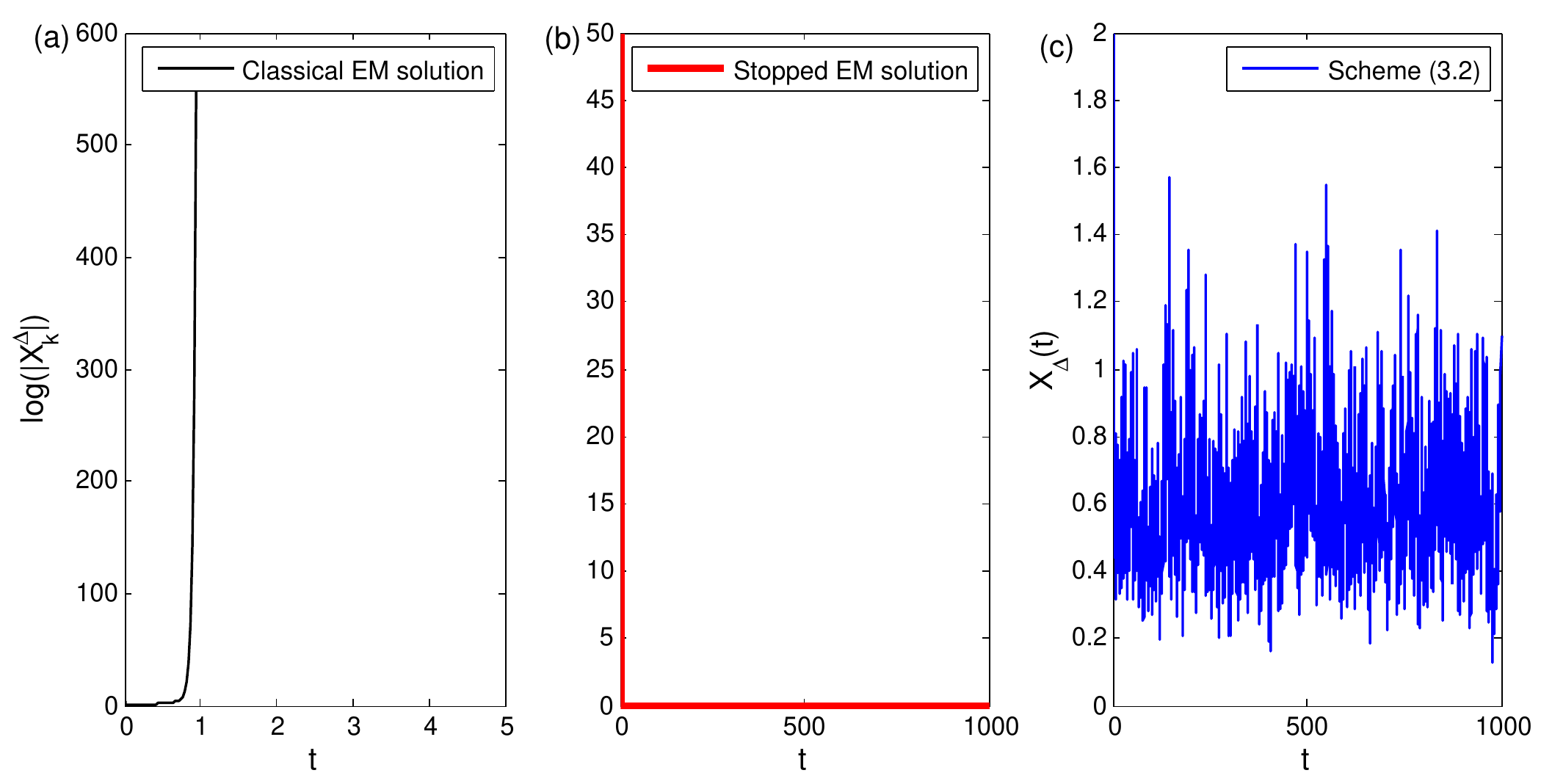}  
\caption{(a) A sample path of the classical EM solution. (b) A sample path of the stopped EM solution. (c) A sample path of the numerical solution of Scheme \eqref{TEM_1}  with the same stepsize $\Delta=0.025$ and $t\in [0,1000]$.}                                \label{exp_3EM_S}
\end{center}                               
\vspace{-2em}
\end{figure}
\begin{figure}[!htp]
\begin{center}
\includegraphics[angle=0, height=6.5cm, width=14cm]{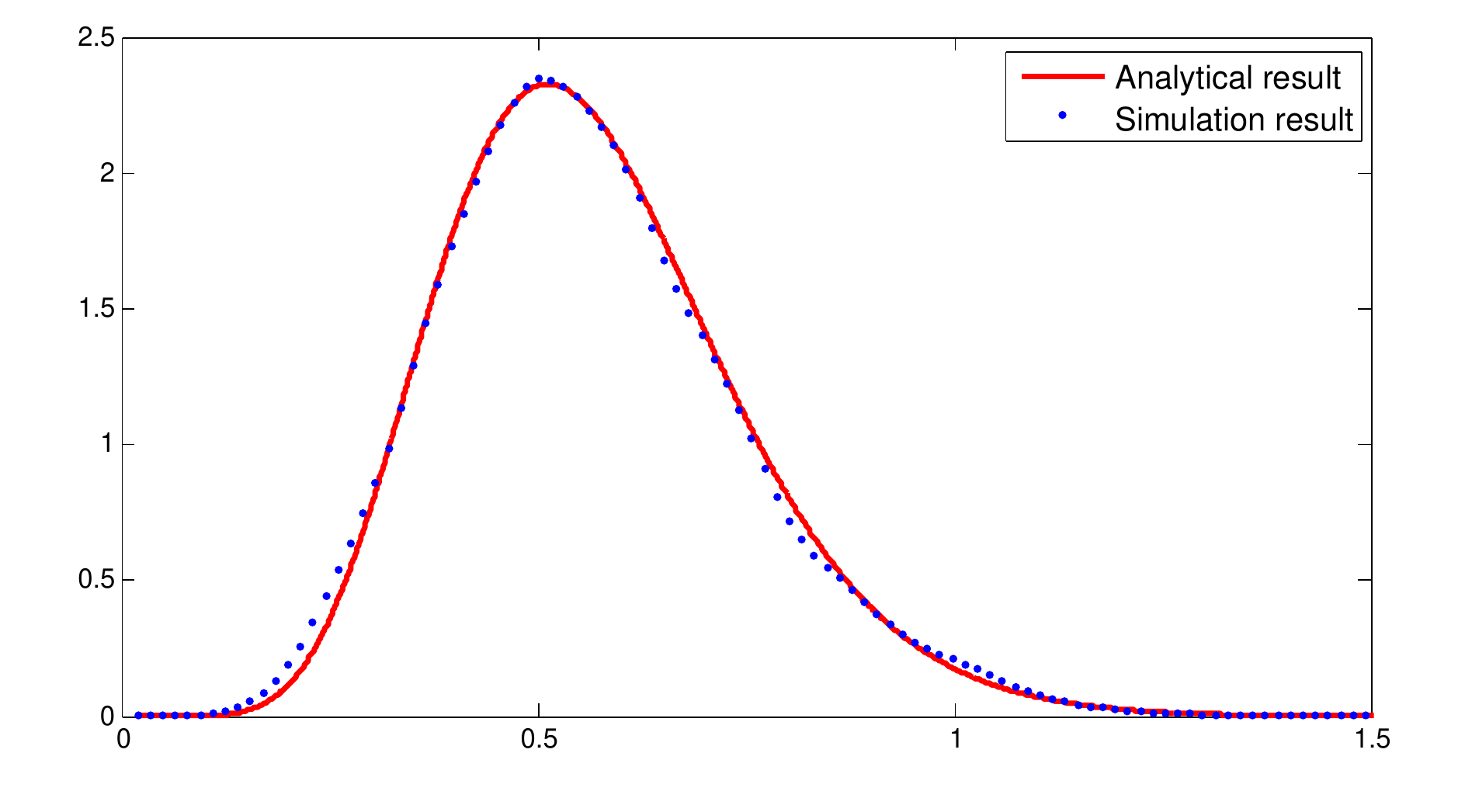}  
\caption{The red solid line indicates the density function of the Gamma distribution  $Ga(91/9, 160/9)$, the solid blue dots indicates the empirical density function of $\mu^{\Delta}$.}                                \label{exp_3density}
\end{center}                               
\vspace{-2em}
\end{figure}

To further illustrate the result of Theorems \ref{log_Th*5.4} and  \ref{log_th*6.3}. First, we generate  sample paths of the exact solution and numerical solution of
Scheme \eqref{TEM_1} in interval $[0, 1000]$, see Fig. \ref{exp_3permanence}. It is evident to see that these two sample paths overlap with each other. Secondly, to measure the similarity quantitatively, we use
  the Kolmogorov-Smirnov  test with a significance level of $0.05$ to check if the stationary distribution of $X_{\Delta}(t)$ is the Gamma distribution. At this level of significance, by {M{\scriptsize ATLAB}} we do confirm that the stationary distribution of $X_{\Delta}(t)$ is the Gamma distribution. So the numerical invariant measure approximates the underlying exact invariant measure very well.
   Finally, to more intuitively illustrate the result of Theorem \ref{log_th*6.3}, we plot the empirical density function of $\mu^{\Delta}(\cdot)$ and the density function of the Gamma distribution $Ga(91/9, 160/9)$ in Fig. \ref{exp_3density}. One observes obviously from the Fig. \ref{exp_3density} that  the computer simulation results obtained with our method approaches the analytical result which can be obtained by the Fokker-Plank equation.
 Furthermore, the similarity between the paths as well as the distributions is significant.  Thus, this example illustrates the significance of the results of Theorems \ref{log_Th*5.4} and  \ref{log_th*6.3}.

}
\end{expl}

\subsection*{Appendix A.}
Using the technique  in the proof Lemma \ref{log_le3.1} yield the following lemma, and hence is omitted to avoid repetition.

\begin{thm}\la{log_le4}
For any $p>0$, the EM scheme defined by \eqref{log_eq5} has the property that
\begin{align*}
\sup_{\Delta\in (0, 1]}\sup_{0\leq k\leq \lfloor T/\Delta\rfloor}\mathbb{E}\big[\mathrm{e}^{pY_{k}}\big]
\leq  C_T,~~~~\forall~T>0,
\end{align*}
 where $\lfloor T/\Delta\rfloor$  represents the integer part of $T/\Delta$.
\end{thm}

{\bf  Proof of Lemma \ref{logeq2.14}.}~~
By \eqref{log_eq4}, we have
\begin{align}\la{log_eq11}
y(t_{k+1})=& y(t_{k})+\int_{t_k}^{t_{k+1}}\big(\beta(r(s))-a(r(s)) \mathrm{e}^{y(s)}\big)\mathrm{d}s+\int_{t_k}^{t_{k+1}}\sigma(r(s)) \mathrm{d}B(s).\tag{A.1}
\end{align}
Using
\eqref{log_eq5} and \eqref{log_eq11} we have
\begin{align*}
Y_{k+1}-y(t_{k+1})
=:&Y_{k}-y(t_{k}) -a(r_k) \big(\mathrm{e}^{Y_k} - \mathrm{e}^{y(t_k)}\big)\Delta+J_k,
\end{align*}
where $J_k=J^{(1)}_{k}+J^{(2)}_{k}+J^{(3)}_{k}
+J^{(4)}_{k}+J^{(5)}_{k}$,   $J^{(1)}_{k}
 = -\int_{t_k}^{t_{k+1}}\big(a(r_k)-a(r(s))\big)x(t_k)  \mathrm{d}s,$
\begin{align*}
\!J^{(2)}_{k}
 \!=& \int_{t_k}^{t_{k+1}}\!\!a(r(s)) \int_{t_k}^{s}x(u) \big(b(r(u))\!-a(r(u))x(u)\big)\mathrm{d}u\mathrm{d}s,
 ~~  J^{(4)}_{k}
 = \! \int_{t_k}^{t_{k+1}}\big[\beta(r_k)\!-\beta(r(s))\big]\mathrm{d}s,\\
 J^{(3)}_{k}
 =& \int_{t_k}^{t_{k+1}}a(r(s)) \int_{t_k}^{s}\sigma(r(u))x(u) \mathrm{d}B(u)\mathrm{d}s,
~~   J^{(5)}_{k}
 =  \int_{t_k}^{t_{k+1}}\big[\sigma(r_k)-\sigma(r(s))\big] \mathrm{d}B(s).
\end{align*}
Define $u_k=Y_{k}-y(t_{k})$. Note that  $u_k\big(\mathrm{e}^{Y_k} -\mathrm{e}^{y(t_k)}\big)\geq 0$, we get
\begin{align*}
u^2_{k+1}
=& u^2_{k}+a^2(r_k) \big(\mathrm{e}^{Y_k} -\mathrm{e}^{y(t_k)}\big)^2\Delta^2+ (J_{k})^{2}+2u_k J_{k} \nn\\
 &
-2a(r_k)u_k\big(\mathrm{e}^{Y_k} -\mathrm{e}^{y(t_k)}\big)\Delta -2a(r_k)J_{k}\big(\mathrm{e}^{Y_k} -\mathrm{e}^{y(t_k)}\big)\Delta \nn\\
\leq& u^2_{k}
 +2a^2(r_k) \big(\mathrm{e}^{Y_k} -\mathrm{e}^{y(t_k)}\big)^2\Delta^2 +2 (J_{k})^2+2u_k J_{k}.
\end{align*}
One further observes that
 \begin{align}\la{logeq:3.14}
\mathbb{E}u^2_{k+1} \leq&  2 \check{a}^2 \Delta^2 \sum_{i=0}^{k}\mathbb{E}\Big[\big(\mathrm{e}^{Y_i} -\mathrm{e}^{y(t_i)}\big)^2\Big] +2\sum_{i=0}^{k}\mathbb{E}\big(u_{i} J_{i}\big)
 +10\sum_{i=0}^{k}\sum_{j=1}^{5}\mathbb{E} (J^{(j)}_{i} )^2.\tag{A.2}
\end{align}
Now the mean value theorem implies
$$
(\mathrm{e}^{x}-\mathrm{e}^{y})^2
\leq(\mathrm{e}^{x}+\mathrm{e}^{y})|\mathrm{e}^{x}-\mathrm{e}^{y}|\leq (\mathrm{e}^{x}+\mathrm{e}^{y})^2|x-y|~~~~\forall~x,y\in \mathbb{R}.
$$
The above inequality together with  Lemma \ref{log_le4} as well as H\"{o}lder's inequality implies
\begin{align*}
 2 \check{a}^2\Delta^{2} \mathbb{E}\Big[ \sum_{i=0}^{k}\big(\mathrm{e}^{Y_i}-\mathrm{e}^{y(t_i)}\big)^{2}\Big]
\leq&2 \check{a}^2 \Delta^{2}   \sum_{i=0}^{k} \mathbb{E}\Big[\big(\mathrm{e}^{Y_i}+\mathrm{e}^{y(t_i)}\big)^{2}
\big|Y_i-y(t_i)\big|\Big]\nn\\
\leq&\Delta   \bigg[ 2^{3}\check{a}^4\Delta^{2}\sum_{i=0}^{k} \big(\mathbb{E}\mathrm{e}^{4Y_i}
+\mathbb{E}\mathrm{e}^{4y(t_i)}\big)  +
 \sum_{i=0}^{k}\mathbb{E}|u_{i}|^{2} \bigg]\nn\\
\leq& CT \Delta^{2}+ \Delta\sum_{i=0}^{k}\mathbb{E}|u_{i}|^{2}.
\end{align*}
This together with \eqref{logeq:3.14} implies
 \begin{align}\la{logeq:3.16}
\mathbb{E}u^2_{k+1} \leq& \Delta\sum_{i=0}^{k}\mathbb{E}|u_{i}|^{2} +CT \Delta^{2} +2\sum_{i=0}^{k}\mathbb{E}\big(u_{i} J_{i} \big) +10\sum_{i=0}^{k}\sum_{j=1}^{5}\mathbb{E} (J^{(j)}_{i} )^2.\tag{A.3}
\end{align}
Then, by the Markov property ((4.16) in
\cite[p.116]{Mao06}) and \eqref{log_eq3}, we derive that
 \begin{align}\la{logeq:3.17}
\!\!\! \mathbb{E} \big[ \big( J^{(1)}_{k}\big)^2 \big]
\!  \leq& \Delta \mathbb{E} \Big[\int_{t_k}^{t_{k+1}} \big(a(r(s))-a(r_k)\big)^2x^2(t_k)\mathrm{d}s\Big]\nn\\
\!\!\!\leq&  4\check{a}^2\Delta\!\int_{t_k}^{t_{k+1}}\!\! \mathbb{E} \Big[x^2(t_k)\mathbb{E} \big(I_{\{r(s)\neq r_k\}}\big|{\cal{F}}_{t_k}\big)\Big]\mathrm{d}s
\!\leq \!   C\Delta^2\int_{t_k}^{t_{k+1}}\!\! \mathbb{E} \big[x^2(t_k) \big]\mathrm{d}s
\!\leq \!  C\Delta^3,\tag{A.4}
\end{align}
and
\begin{align}\la{logeq:3.18}
 \mathbb{E} \big[\big(J^{(2)}_{k}\big)^2 \big]
\leq&  \Delta \mathbb{E} \Big[  \int_{t_k}^{t_{k+1}} a^2(r(s))  \Big(\int_{t_k}^{s}x(u) \big(b(r(u))-a(r(u))x(u)\big)\mathrm{d}u\Big)^2\mathrm{d}s\Big]\nn\\
\leq&  \check{a}^2\Delta^2\mathbb{E} \Big[  \int_{t_k}^{t_{k+1}}\int_{t_k}^{s} x^2(u) \big(b(r(u))-a(r(u))x(u)\big)^2\mathrm{d}u\mathrm{d}s\Big]\nn\\
\leq&  2\check{a}^2\Delta^2  \int_{t_k}^{t_{k+1}}\int_{t_k}^{s}\Big(\check{b}^2 \mathbb{E} \big[x^2(u)\big] +\check{a}^2 \mathbb{E} \big[x^4(u)\big]\Big)\mathrm{d}u\mathrm{d}s
\leq  C\Delta^4.\tag{A.5}
\end{align}
By the It\^{o} isometry and \eqref{log_eq3},  we have
  \begin{align}\la{logeq:3.19}
  \mathbb{E} \big[ \big( J^{(3)}_{k}\big)^2  \big]
\leq& \Delta\mathbb{E} \Big[\int_{t_k}^{t_{k+1}}a^2(r(s))\Big( \int_{t_k}^{s}\sigma(r(u))x(u) \mathrm{d}B(u)\Big)^2\mathrm{d}s\Big]\nn\\
\leq& \check{a}^2 \Delta\mathbb{E} \Big[\int_{t_k}^{t_{k+1}}\Big( \int_{t_k}^{s}\sigma(r(u))x(u) \mathrm{d}B(u)\Big)^2\mathrm{d}s\Big]\nn\\
\leq& |\breve{\sigma}|^2 \check{a}^2 \Delta\int_{t_k}^{t_{k+1}}\int_{t_k}^{s}\mathbb{E} \big[x^2(u)\big] \mathrm{d}u  \mathrm{d}s\leq C\Delta^3.\tag{A.6}
\end{align}
Furthermore, using Lemma 6.10 in \cite[p.251]{Mao06}, we yield that
 \begin{align}\la{logeq:3.20}
 \mathbb{E} \big[ \big(J^{(4)}_{k}\big)^2  \big]
\leq& 4\check{\beta}^2\Delta\int_{t_k}^{t_{k+1}}   \mathbb{P} \big(r(s)\neq r_k\big)\mathrm{d}s
\leq C\Delta^3,\tag{A.7}
\end{align}
and
 \begin{align}\la{logeq:3.21}
  \mathbb{E} \big[\big(J^{(5)}_{k}\big)^2  \big]
\leq& 4|\breve{\sigma}|^2 \int_{t_k}^{t_{k+1}} \mathbb{P} \big(r(s)\neq r_k\big)\mathrm{d}s
\leq C\Delta^2.\tag{A.8}
\end{align}
Combining \eqref{logeq:3.16}-\eqref{logeq:3.21}, we obtain that
\begin{align}\la{logeq:3.22}
\mathbb{E}u^2_{k+1} \leq& \Delta\sum_{i=0}^{k}\mathbb{E}|u_{i}|^{2} +CT \Delta  +2\sum_{i=0}^{k}\mathbb{E}\Big[u_{i}\big(J^{(1)}_{i}+J^{(2)}_{i}
+J^{(3)}_{i}
+J^{(4)}_{i}+J^{(5)}_{i}\big)\Big].\tag{A.9}
\end{align}
Using the H\"{o}lder inequality and \eqref{logeq:3.17}, we obtain that
 \begin{align*}
2\sum_{i=0}^{k}\mathbb{E} \big[u_{i} J^{(1)}_{i} \big]
\leq& 2\sum_{i=0}^{k} \big(\mathbb{E}u^2_{i}\big)^{\frac{1}{2}}\Big(\mathbb{E} \Big[ \big(J^{(1)}_{i}\big)^2\Big]\Big)^{\frac{1}{2}}
 \leq 2\sum_{i=0}^{k} \big(\mathbb{E}u^2_{i}\big)^{\frac{1}{2}}\Big( C\Delta^3 \Big)^{\frac{1}{2}}\nn\\
 \leq&    2 \Delta^{\frac{1}{2}}\Big(\sum_{i=0}^{k} \mathbb{E}u^2_{i}\Big)^{\frac{1}{2}}\Big( CT\Delta \Big)^{\frac{1}{2}}
 \leq    \sum_{i=0}^{k} \mathbb{E}u^2_{i}\Delta + CT\Delta.
\end{align*}
By \eqref{log_eq3}  we derive that
 \begin{align*}
2\mathbb{E} \big[u_{k} J^{(2)}_{k} \big]
\leq& \mathbb{E}u^2_{k}\Delta+ \mathbb{E} \Big[  \int_{t_k}^{t_{k+1}} a^2(r(s))  \Big(\int_{t_k}^{s}x(u) \big(b(r(u))-a(r(u))x(u)\big)\mathrm{d}u\Big)^2\mathrm{d}s\Big]\nn\\
\leq&\mathbb{E}u^2_{k}\Delta+ \check{a}^2\Delta\mathbb{E} \Big[  \int_{t_k}^{t_{k+1}}\int_{t_k}^{s} x^2(u) \big(b(r(u))-a(r(u))x(u)\big)^2\mathrm{d}u\mathrm{d}s\Big]\nn\\
\leq&\mathbb{E}u^2_{k}\Delta+ 2\check{a}^2\Delta  \int_{t_k}^{t_{k+1}}\int_{t_k}^{s}\Big(\check{b}^2 \mathbb{E} \big[x^2(u)\big] +\check{a}^2 \mathbb{E} \big[x^4(u)\big]\Big)\mathrm{d}u\mathrm{d}s
\leq  \mathbb{E}u^2_{k}\Delta+ C\Delta^3.
\end{align*}
Since
$$
\mathbb{E}\Big[J^{(3)}_{k}|{\cal{F}}_{t_k}\Big]
=\mathbb{E}\bigg[\int_{t_k}^{t_{k+1}}\int_{t_k}^{s}
a(r(s))\sigma(r(u))x(u)\mathrm{d}B(u) \mathrm{d}s\Big|{\cal{F}}_{t_k}\bigg]=0,
$$
 we have that $
 \mathbb{E} \big[u_{k} \big(J^{(3)}_{k}+J^{(5)}_{k}\big) \big]
 =0.
$
By the Markov property, one observes
  \begin{align*}
2\mathbb{E} \big[u_{k}  J^{(4)}_{k}  \big]
\leq& \mathbb{E}u^2_{k}\Delta+2\check{\beta}^2\int_{t_k}^{t_{k+1}} \mathbb{E} \Big[ \mathbb{E} \big(I_{\{r(s)\neq r_k\}}\big|{\cal{F}}_{t_k}\big)\Big]\mathrm{d}s
\leq
 \mathbb{E}u^2_{k}\Delta+C\Delta^2.
\end{align*}
By inserting these four estimates in \eqref{logeq:3.22} we end up with
 \begin{align*}
\mathbb{E}u^2_{k+1} \leq& 4\Delta\sum_{i=0}^{k}\mathbb{E}|u_{i}|^{2}+CT \Delta
\end{align*}
and Gronwall's Lemma completes the proof of the assertion.
\eproof

\subsection*{Appendix B.}
{\bf  Proof of Lemma \ref{log_le*6.10}.}~~
Using the well-known  Taylor formula we get
\begin{align*}
\mathrm{e}^{\bar{Z}_{k+1}}= \mathrm{e}^{Z_{k}}+\mathrm{e}^{Z_{k}}(\bar{Z}_{k+1}-Z_{k})
+\frac{1}{2}\mathrm{e}^{\bar{\xi}_k}(\bar{Z}_{k+1}-Z_{k})^2,
\end{align*}
where $\bar{\xi}_k\in (\bar{Z}_{k+1}\wedge Z_{k}, \bar{Z}_{k+1}\vee Z_{k})$. Clear,
 \begin{align*}
&  \frac{1}{2}\mathrm{e}^{\bar{\xi}_k}(\bar{Z}_{k+1}-Z_{k})^2
\leq \frac{1}{2}\mathrm{e}^{Z_{k}}\mathrm{e}^{|\bar{Z}_{k+1}-Z_{k}|}(\bar{Z}_{k+1}
-Z_{k})^2\nn\\
=&\frac{1}{2}\mathrm{e}^{Z_{k}}\exp{\Big(\big|\big( \beta(r_k)-a(r_k) \mathrm{e}^{Z_k}\big)\Delta +\sigma(r_k) \Delta B_k\big|\Big)}\Big|\big(\beta(r_k)-a(r_k) \mathrm{e}^{Z_k}\big)\Delta +\sigma(r_k) \Delta B_k\Big|^2\nn\\
\leq&\mathrm{e}^{Z_{k}}\tilde{\mathcal{U}}_k,
\end{align*}
where
$$
\tilde{\mathcal{U}}_k:=  C\exp{\Big( |\breve{\sigma}| |\Delta B_k|\Big)}\Big(  \Delta^{2(1-\theta)} +  |\Delta B_k|^2\Big).
$$
Therefore, we derive from \eqref{log:eq:5.8} and \eqref{T_w2} that for any integer $k\geq 0$,
\begin{align*}
 \mathrm{e}^{\bar{Z}_{k+1}}\leq& \mathrm{e}^{Z_{k}}+ \mathrm{e}^{Z_{k}}\Big[\big(\beta(r_k)-a(r_k) \mathrm{e}^{Z_k}\big)\Delta +\sigma(r_k) \Delta B_k+\tilde{\mathcal{U}}_k\Big].
\end{align*}
Then using the above inequality,   we have
\begin{align}\la{2}
 (1+\mathrm{e}^{\bar{Z}_{k+1}})^{p}\leq& (1+\mathrm{e}^{Z_{k}})^p\big(1+ \bar{\varsigma}_k\big)^p.\tag{B.1}
\end{align}
where
$$
 \bar{\varsigma}_k:= (1+\mathrm{e}^{Z_{k}})^{-1}\mathrm{e}^{Z_{k}}\Big[\big(\beta(r_k)-a(r_k) \mathrm{e}^{Z_k}\big)\Delta +\sigma(r_k) \Delta B_k+\tilde{\mathcal{U}}_k\Big],
$$
and we can see that $\bar{\varsigma}_k>-1$. By  virtue of \cite[Lemma 3.3]{li2018jde},
without loss the generality we prove \eqref{1} only for $0<p\leq 1$. It follows from \eqref{2} that
\begin{align}
&\E \Big[\big(1+\mathrm{e}^{-\bar{Z}_{k+1}}\big)^{p}  \big|\mathcal{F}_{t_k}\Big]\nn\\
\leq&\big(1+\mathrm{e}^{-Z_{k}}\big)^{p} \bigg\{1+ p \E\big[\bar{\varsigma}_k \big|\mathcal{F}_{t_k}\big] + \frac{p(p-1)}{2} \E\big[\bar{\varsigma}_k^2\big|\mathcal{F}_{t_k}\big]+ \frac{p(p-1)(p-2)}{6} \E\big[\bar{\varsigma}_k^3\big|\mathcal{F}_{t_k}\big]\bigg\}.
\tag{B.2}
\end{align}
  By \eqref{logeq:5.10} and Lemma \ref{leB}, we derive that
\begin{align*}
\E\big[\tilde{\mathcal{U}}_k\big|\mathcal{F}_{t_k}\big]
\leq&C \E \bigg[   \Delta^{2(1-\theta)}\mathrm{e}^{\frac{|\Delta B_k|^2}{4\Delta} } + \mathrm{e}^{\frac{|\Delta B_k|^2}{4\Delta} } |\Delta B_k|^2 \bigg]  \leq C\Delta,
\end{align*}
and
\begin{align*}
 \E\big[\bar{\varsigma}_k \big|\mathcal{F}_{t_k}\big]  \leq & (1+\mathrm{e}^{Z_{k}})^{-1}\mathrm{e}^{Z_{k}}\Big[\big(\beta(r_k)-a(r_k) \mathrm{e}^{Z_k}\big)\Delta  +C\t\Big]\nn\\
  \leq & -a(r_k)(1+\mathrm{e}^{Z_{k}})^{-1}\mathrm{e}^{2Z_{k}} \Delta  +C(1+\mathrm{e}^{Z_{k}})^{-1}\mathrm{e}^{Z_{k}} \t \nn\\
=& -a(r_k)(1+\mathrm{e}^{Z_{k}})\Delta+2 a(r_k) \Delta-a(r_k)(1+\mathrm{e}^{Z_{k}})^{-1} \Delta \nn\\ &+C\t-C(1+\mathrm{e}^{Z_{k}})^{-1}  \t\nn\\
\leq& -a(r_k)(1+\mathrm{e}^{Z_{k}})\Delta+C\t-C(1+\mathrm{e}^{Z_{k}})^{-1}  \t.
\end{align*}
Similarly, we can also prove that
\begin{align*}
 \E\big[\bar{\varsigma}_k^2 \big|\mathcal{F}_{t_k}\big]  \geq -C \t,~~~~~ \E\big[\bar{\varsigma}_k^3 \big|\mathcal{F}_{t_k}\big]\leq  C \t^{\frac{3}{2}}.
\end{align*}
Thus, we obtain that
\begin{align*}
&\E \Big[\big(1+\mathrm{e}^{-\bar{Z}_{k+1}}\big)^{p} \Big]\nn\\
\leq&\E \Big[\big(1+\mathrm{e}^{-Z_{k}}\big)^{p} \Big(1 -pa(r_k)(1+\mathrm{e}^{Z_{k}})\Delta+C\t-C(1+\mathrm{e}^{Z_{k}})^{-1}  \t\Big)\Big]\nn\\
\leq&-p \hat{a}\E \Big[\big(1+\mathrm{e}^{-Z_{k}}\big)^{p+1}\Big]\Delta
+(1+C\t)\E \Big[\big(1+\mathrm{e}^{-Z_{k}}\big)^{p}\Big] -C\E \Big[(1+\mathrm{e}^{Z_{k}})^{p-1}\Big]  \t
\leq  C.
\end{align*}
The proof is therefore complete.\eproof

{\bf  Proof of Lemma \ref{log:Le5.6}.}~~Since $\pi a =0$, we have
\begin{align*}
Z_{k+1}=\bar{Z}_{k+1}
=&Z_{k}+ \beta(r_k) \Delta +\sigma(r_k) \Delta B_k .
\end{align*}
It is easy to see that
\begin{align*}
\mathrm{e}^{Z_{k+1}}\leq \mathrm{e}^{Z_{k}}+\mathrm{e}^{Z_{k}}(Z_{k+1}-Z_{k})+\frac{1}{2}\mathrm{e}^{Z_{k}}(Z_{k+1}-Z_{k})^2
+\mathrm{e}^{Z_{k}}\bar{\mathcal{U}}_k,
\end{align*}
implies that
 \begin{align}\la{logeq:5.27}
 \mathrm{e}^{Z_{k+1}}
\leq   \mathrm{e}^{Z_{k}} \big(1+\varsigma_k\big),\tag{B.3}
\end{align}
where
 \begin{align*}
\varsigma_k=& \beta(r_k) \Delta +\sigma(r_k) \Delta B_k+\frac{\sigma^2(r_k)}{2} \big(\Delta B_k\big)^2 +\sigma(r_k) \Delta B_k \beta(r_k)  \Delta
+\frac{1}{2} \beta^2(r_k)\Delta^2+\bar{\mathcal{U}}_k,
\end{align*}
and we can see that $\varsigma_k>-1$. Note that
 \begin{eqnarray*}
\pi[\beta-(\pi \beta)\I_m]=0,\ \ \ \ \ \sum_{i=1}^{m}\pi_i=1.
\end{eqnarray*}
It follows from Lemma \ref{Le:Yin} (1) that the equation
 \begin{eqnarray*}
\dis \Gamma v=\beta-(\pi \beta)\I_m
\end{eqnarray*}
has a solution $v=(v_1,   \cdots, v_m)^T\in \mathbb{R}^m$. Thus we have
 \begin{align} \la{eq_4.27}
-\beta(i)+\sum_{j=1}^{m}\gamma_{i j}v_j=-\pi \beta>0,\ \ \ i\in \mathbb{S}.\tag{B.4}
\end{align}
Choose a constant $0<\rho_0\leq 1$ such that for each $0<\rho\leq \rho_0$,   $$\xi^{\rho,v}_i:=1-v_i \rho>0,  \ i=1,\cdots,m.$$
Using the techniques in the proof of Lemma \ref{L:3},  it follows from \eqref{logeq:5.27} that
\begin{align*}
 \E \Big[ \mathrm{e}^{\rho Z_{k+1}} \xi^{\rho,v}_{r_{k+1}} \big|\mathcal{F}_{t_k}\Big]
\leq& \mathrm{e}^{\rho Z_{k}}  \bigg\{\E\big[\xi^{\rho,v}_{r_{k+1}}\big|\mathcal{F}_{t_k}\big]+\rho \E\big[\varsigma_k\xi^{\rho,v}_{r_{k+1}}\big|\mathcal{F}_{t_k}\big]  + \frac{\rho(\rho-1)}{2} \E\big[\varsigma_k^2\xi^{\rho,v}_{r_{k+1}}\big|\mathcal{F}_{t_k}\big]\nn\\ &~~~~~~~~~~~~~~~~~~~~ +\frac{\rho(\rho-1)(\rho-2)}{6} \E\big[\varsigma_k^3\xi^{\rho,v}_{r_{k+1}}\big|\mathcal{F}_{t_k}\big]\bigg\}.
\end{align*}
  By \eqref{logeq:5.9}, \eqref{logeq:5.10} and Lemma \ref{leB}, we derive that
\begin{align*}
\E\big[\varsigma_k\xi^{\rho,v}_{r_{k+1}}\big|\mathcal{F}_{t_k}\big]
 =&   \E\big[\xi^{\rho,v}_{r_{k+1}}\E\big(\varsigma_k\big|\mathcal{G}_{t_k}\big)\big|\mathcal{F}_{t_k}\big]
\leq   \Big[ \beta(r_k)\Delta+\frac{\sigma^2(r_k)}{2}\Delta
 +C\Delta^{\frac{3}{2}}\Big]
\E\big[\xi^{\rho,v}_{r_{k+1}}\big|\mathcal{F}_{t_k}\big]\nn\\
\leq&  \big(\beta(r_k) +\frac{\sigma^2(r_k)}{2}\big)\xi^{\rho,v}_{r_{k}}\Delta  +C\Delta^{\frac{3}{2}}
,
\end{align*}
and
\begin{align*}
\E\big[\varsigma^2_k \big|\mathcal{G}_{t_k}\big]=&\E\bigg\{   \Big[\beta(r_k) \Delta +\sigma(r_k) \Delta B_k+\frac{\sigma^2(r_k)}{2} \big(\Delta B_k\big)^2 \nn\\
&~~~~~~~~~~~ +\sigma(r_k) \Delta B_k \beta(r_k)  \Delta
+\frac{1}{2} \beta^2(r_k)\Delta^2+\bar{\mathcal{U}}_k\Big]^2\Big|\mathcal{G}_{t_k}\bigg\}\nn\\
\geq&\E\bigg\{ \sigma^2(r_k) \big(\Delta B_k\big)^2 +2\sigma(r_k) \Delta B_k\Big[ \beta(r_k) \Delta +\frac{\sigma^2(r_k)}{2} \big(\Delta B_k\big)^2\nn\\
& ~~~~~~~~~~~~~~~~~~~~~~~~~~~~~~~~~~~~~~+\sigma(r_k) \Delta B_k \beta(r_k) \Delta
+\frac{1}{2} \beta^2(r_k) \Delta^2+\bar{\mathcal{U}}_k\Big]  \Big|\mathcal{G}_{t_k}\bigg\}\nn\\
=&   \sigma^2(r_k) \Delta  +2\sigma^2(r_k) \beta(r_k)  \Delta^2
 +2\sigma(r_k) \E\big[\Delta B_k\bar{\mathcal{U}}_k\big|\mathcal{G}_{t_k}\big]
\geq \sigma^2(r_k) \Delta -C \Delta^{2},
\end{align*}
implies that
\begin{align*}
\E\big[\varsigma_k^2\xi^{\rho,v}_{r_{k+1}}\big|\mathcal{F}_{t_k}\big]
 =&   \E\big[\xi^{\rho,v}_{r_{k+1}}\E\big(\varsigma^2_k\big|\mathcal{G}_{t_k}\big)\big|\mathcal{F}_{t_k}\big]
\geq \sigma^2(r_k)\xi^{\rho,v}_{r_{k}} \Delta -C \Delta^{2}.
\end{align*}
Similarly, we can also prove that $ \E\big[\varsigma_k^{3}\xi^{\rho,v}_{r_{k+1}}|\F_{t_k}\big]\leq C\t^{\frac{3}{2}}.$
 Thus, we obtain that
\begin{align}\la{eq_4.28}
&\E \Big[ \mathrm{e}^{\rho Z_{k+1}}   \xi^{\rho,v}_{r_{k+1}} \big|\mathcal{F}_{t_k}\Big]\nn\\
\leq& \mathrm{e}^{\rho Z_{k}}  \bigg\{\xi^{\rho,v}_{r_{k}}+\sum_{j\in \mathbb{S}}\xi^{\rho,v}_{j}\left(\gamma_{r_k j}\t+o(\t)\right) +\rho \big(\beta(r_k) +\frac{\sigma^2(r_k)}{2}\big)\xi^{\rho,v}_{r_{k}}\Delta  +C\Delta^{\frac{3}{2}}\nn\\
 &~~~~~~~~~~~~~~~~~~~~~~ ~~~~~+ \frac{\rho(\rho-1)}{2} \sigma^2(r_k)\xi^{\rho,v}_{r_{k}}\Delta +C \Delta^{2}
  \bigg\}\nn\\
\leq& \mathrm{e}^{\rho Z_{k}}\xi^{\rho,v}_{r_{k}}  \bigg\{1+\frac{1}{\xi^{\rho,v}_{r_{k}}}\sum_{j\in \mathbb{S}}\xi^{\rho,v}_{j} \gamma_{r_k j}\t +\rho \beta(r_k) \Delta + \frac{\rho^2}{2} \sigma^2(r_k) \Delta +o(\t)
  \bigg\}.\tag{B.5}
\end{align}
By the properties of the generator, we have
\begin{align}\la{eq_4.29}
\frac{1}{\xi^{\rho,v}_{i}}\sum_{j\in \mathbb{S}}\xi^{\rho,v}_{j} \gamma_{i j}=&\frac{1}{1-v_i \rho}\sum_{j=1}^m(1-v_j \rho) \gamma_{i j}\nn\\
=&-\rho\bigg(\sum_{j=1}^m \gamma_{i j} v_j+\frac{v_i \rho}{1-v_i \rho}\sum_{j=1}^m \gamma_{i j} v_j\bigg).\tag{B.6}
\end{align}
It follows from \eqref{eq_4.27}, \eqref{eq_4.28} and \eqref{eq_4.29} that
\begin{align*}
&\E \Big[ \mathrm{e}^{\rho Z_{k+1}}   \xi^{\rho,v}_{r_{k+1}} \big|\mathcal{F}_{t_k}\Big]\nn\\
\leq& \mathrm{e}^{\rho Z_{k}}\xi^{\rho,v}_{r_{k}}  \bigg\{1-\rho\bigg(\sum_{j=1}^m \gamma_{{r_k} j} v_j+\frac{v_{r_k} \rho}{1-v_{r_k} \rho}\sum_{j=1}^m \gamma_{{r_k} j} v_j\bigg)\t +\rho \beta(r_k)  \Delta + \frac{\rho^2}{2} \sigma^2(r_k) \Delta +o(\t)
  \bigg\}\nn\\
=& \mathrm{e}^{\rho Z_{k}}\xi^{\rho,v}_{r_{k}}  \bigg\{1-\rho\bigg(\beta(r_k)-\pi \beta+\frac{v_{r_k} \rho}{1-v_{r_k} \rho}\sum_{j=1}^m \gamma_{{r_k} j} v_j- \frac{\rho}{2} \sigma^2(r_k) \Delta\bigg)\t +\rho \beta(r_k)  \Delta +o(\t)
  \bigg\}\nn\\
  =& \mathrm{e}^{\rho Z_{k}}\xi^{\rho,v}_{r_{k}}  \bigg\{1-\rho\bigg(-\pi \beta+\frac{v_{r_k} \rho}{1-v_{r_k} \rho}\sum_{j=1}^m \gamma_{{r_k} j} v_j- \frac{\rho}{2} \sigma^2(r_k) \Delta\bigg)\t +o(\t)
  \bigg\}.
\end{align*}
Choose a constant $0<\rho_1\leq \rho_0$ such that for any $0<\rho\leq \rho_1$,
\begin{eqnarray*}
\lambda_i:=-\pi \beta+\rho\bigg(\frac{v_i}{1-v_i \rho}\sum\limits_{j=1}^{m}\gamma_{ij}v_j-\frac{\rho}{2} \sigma^2(i)\bigg)>0.
\end{eqnarray*}
 Choose    $\Delta_2^{*}\in(0,1)$ sufficiently small
such that
$\Delta_2^{*}< 2/\rho\hat{\lambda}$,   $o\left(\Delta_2^{*}\right)\leq \rho\hat{\lambda} \Delta_2^{*}/2$. Then, for any $\Delta\in(0, \Delta_2^{*}]$  yields
\begin{align*}
 \E \Big[ \mathrm{e}^{\rho Z_{k+1}}   \xi^{\rho,v}_{r_{k+1}} \big|\mathcal{F}_{t_k}\Big] \leq \big(1-\frac{\rho\hat{\lambda}}{2} \t
  \big)\mathrm{e}^{\rho Z_{k}}\xi^{\rho,v}_{r_{k}}.
\end{align*}
Since the proof method is as that of Lemma \ref{L:3} we omit the details. \eproof

{\bf  Proof of Lemma \ref{log_le*6.11}.}~~Obviously,
\begin{align*}
\bar{Z}^{y_0,\ell}_{k+1}
= Z^{y_0,\ell}_k+\big(\beta(r^{\ell}_k)-a(r^{\ell}_k)\mathrm{e}^{Z^{y_0,\ell}_{k}}\big)\Delta+\sigma(r^{\ell}_k) \Delta B_k,
\end{align*}
and
\begin{align*}
\bar{Z}^{\bar{y}_0,\ell}_{k+1}
= Z^{\bar{y}_0,\ell}_k+\big(\beta(r^{\ell}_k)-a(r^{\ell}_k)\mathrm{e}^{Z^{\bar{y}_0,\ell}_{k}}\big)\Delta+\sigma(r^{\ell}_k) \Delta B_k.
\end{align*}
It is easy to see that
\begin{align*}
\!\bar{Z}^{y_0,\ell}_{k+1}\!-\bar{Z}^{\bar{y}_0,\ell}_{k+1}
=\! Z^{y_0,\ell}_{k}\!- Z^{\bar{y}_0,\ell}_k \!-a(r^{\ell}_k)\big(\mathrm{e}^{Z^{y_0,\ell}_{k}}
\!-\mathrm{e}^{Z^{\bar{y}_0,\ell}_{k}}\big)\Delta=\!  Z^{y_0,\ell}_{k}- Z^{\bar{y}_0,\ell}_k  \!-a(r^{\ell}_k)\big(X^{x_0,\ell}_{k}
\!-X^{\bar{x}_0,\ell}_{k} \big)\Delta,
\end{align*}
which implies that
\begin{align}\la{log_eq*6.20}
\big|\bar{Z}^{y_0,\ell}_{k+1}-\bar{Z}^{\bar{y}_0,\ell}_{k+1}\big|
 = \big| Z^{y_0,\ell}_{k}- Z^{\bar{y}_0,\ell}_k \big| -a(r^{\ell}_k)\big|X^{x_0,\ell}_{k}
-X^{\bar{x}_0,\ell}_{k} \big|\Delta.\tag{B.7}
\end{align}
Note that
\begin{align}\la{log_eq*6.21}
\big| Z^{y_0,\ell}_{k}- Z^{\bar{y}_0,\ell}_{k} \big|\leq&\big|\bar{Z}^{y_0,\ell}_{k}-\bar{Z}^{\bar{y}_0,\ell}_{k}\big|.\tag{B.8}
\end{align}
In fact, if $\bar{Z}^{y_0,\ell}_{k}\vee\bar{Z}^{\bar{y}_0,\ell}_{k}\leq \log(K\Delta^{-\theta } )$,  \eqref{log_eq*6.21}  holds obviously. If $\bar{Z}^{y_0,\ell}_{k}\wedge\bar{Z}^{\bar{y}_0,\ell}_{k}> \log(K\Delta^{-\theta } )$, \eqref{log_eq*6.21}  holds obviously. If $\bar{Z}^{y_0,\ell}_{k}\leq  \log(K\Delta^{-\theta } )<\bar{Z}^{\bar{y}_0,\ell}_{k}$, we have $$Z^{y_0,\ell}_{k}=\bar{Z}^{y_0,\ell}_{k}\leq Z^{\bar{y}_0,\ell}_{k}=\log(K\Delta^{-\theta } )<\bar{Z}^{\bar{y}_0,\ell}_{k},$$
and
\begin{align*}
 \big| Z^{y_0,\ell}_{k}- Z^{\bar{y}_0,\ell}_{k} \big|^2-\big|\bar{Z}^{y_0,\ell}_{k}-\bar{Z}^{\bar{y}_0,\ell}_{k}\big|^2
=& \big(Z^{\bar{y}_0,\ell}_{k}+\bar{Z}^{\bar{y}_0,\ell}_{k}\big)\big(Z^{\bar{y}_0,\ell}_{k}-\bar{Z}^{\bar{y}_0,\ell}_{k}\big)- 2Z^{y_0,\ell}_{k}\big(Z^{\bar{y}_0,\ell}_{k}
-\bar{Z}^{\bar{y}_0,\ell}_{k}\big)  \nn\\
=& \big(Z^{\bar{y}_0,\ell}_{k}+\bar{Z}^{\bar{y}_0,\ell}_{k}- 2Z^{y_0,\ell}_{k}\big)
\big(Z^{\bar{y}_0,\ell}_{k}-\bar{Z}^{\bar{y}_0,\ell}_{k}\big)\leq 0.
\end{align*}
Then \eqref{log_eq*6.21} follows immediately. If $\bar{Z}^{\bar{y}_0,\ell}_{k}\leq  \log(K\Delta^{-\theta } )<\bar{Z}^{y_0,\ell}_{k}$, \eqref{log_eq*6.21} holds also by
symmetry on $\bar{Z}^{\bar{y}_0,\ell}_{k}$ and $\bar{Z}^{y_0,\ell}_{k}$.   Thus, the desired inequality \eqref{log_eq*6.21} holds for all cases. It follows from \eqref{log_eq*6.20} and \eqref{log_eq*6.21} that for any integer $k\geq 0$,
\begin{align*}
\big| Z^{y_0,\ell}_{k+1}- Z^{\bar{y}_0,\ell}_{k+1} \big|\leq& \big| Z^{y_0,\ell}_{k}- Z^{\bar{y}_0,\ell}_k \big| -a(r^{\ell}_k)\big|X^{x_0,\ell}_{k}
-X^{\bar{x}_0,\ell}_{k} \big|\Delta\nn\\
\leq& | y_0 - \bar{y}_0  | -\sum_{i=0}^{k}a(r^{\ell}_i)\big|X^{x_0,\ell}_{i}
-X^{\bar{x}_0,\ell}_{i} \big|\Delta.
\end{align*}
Then we have
$
 \mathbb{E}\big| Z^{y_0,\ell}_{k+1}- Z^{\bar{y}_0,\ell}_{k+1} \big| \leq  | y_0 - \bar{y}_0  | -\hat{a}\sum_{i=0}^{k} \mathbb{E}\big|X^{x_0,\ell}_{i}
-X^{\bar{x}_0,\ell}_{i} \big|\Delta.
$
Due to $\hat{a}>0$,
\begin{align*}
 \sum_{i=0}^{\infty} \mathbb{E}\big|X^{x_0,\ell}_{i}
-X^{\bar{x}_0,\ell}_{i} \big|\Delta \leq  \frac{| y_0 - \bar{y}_0  |}{\hat{a}}<\infty.
\end{align*}
The proof is therefore complete.\eproof

{\bf  Proof of Lemma \ref{log_le*6.12}.}~~By \eqref{log:eq:5.9}, we see that
\begin{align*}
 \mathrm{e}^{Z_{k+1}}
\leq& x_{0}+ \sum_{i=0}^{k}\mathrm{e}^{Z_{i}}\bigg[\big(\beta(r_i)-a(r_i) \mathrm{e}^{Z_i}\big)\Delta +\sigma(r_i) \Delta B_i+\sigma(r_i) \Delta B_i\big(\beta(r_i)-a(r_i) \mathrm{e}^{Z_i}\big) \Delta\nn\\
&~~~~~~~~~~~~ +\frac{\sigma^2(r_i)}{2} \big(\Delta B_i\big)^2
+\frac{1}{2}\big(\beta(r_i)-a(r_i) \mathrm{e}^{Z_i}\big)^2\Delta^2+\bar{C} \mathrm{e}^{\frac{|\Delta B_i|^2}{8\Delta}}\Big(  \Delta^3 +  |\Delta B_i|^3\Big)\bigg],
\end{align*}
where $\bar{C}= {2}/{3}\exp{\big((\check{\beta}+2|\breve{\sigma}|^2)\Delta\big)}
\big(|\check{\beta}|\vee|\breve{\sigma}|\big)^3$.
Then we have
\begin{align*}
 \mathrm{e}^{2Z_{k+1}}
\leq\! 8\bigg\{& x^2_{0}\! + \Big[\sum_{i=0}^{k}\mathrm{e}^{Z_{i}}\big(\beta(r_i)-a(r_i) \mathrm{e}^{Z_i}\big)\Delta\Big]^2 \!+\Big[\sum_{i=0}^{k}\mathrm{e}^{Z_{i}}\sigma(r_i) \Delta B_i\Big]^2
\!+\frac{|\breve{\sigma}|^4}{4}\Big[\sum_{i=0}^{k}\mathrm{e}^{Z_{i}} \big(\Delta B_i\big)^2\Big]^2\nn\\
&+\Delta^2\Big[\sum_{i=0}^{k}\mathrm{e}^{Z_{i}}\big(\beta(r_i)-a(r_i) \mathrm{e}^{Z_i}\big) \sigma(r_i) \Delta B_i\Big]^2
+\frac{\Delta^4}{4}\Big[\sum_{i=0}^{k}\mathrm{e}^{Z_{i}}\big(\beta(r_i)-a(r_i) \mathrm{e}^{Z_i}\big)^2\Big]^2\nn\\
&
+\bar{C}^2\Delta^6\Big[\sum_{i=0}^{k}\mathrm{e}^{Z_{i}} \exp{\Big(\frac{|\Delta B_i|^2}{8\Delta}\Big)} \Big]^2
 +\bar{C}^2\Big[\sum_{i=0}^{k}\mathrm{e}^{Z_{i}} \exp{\Big(\frac{|\Delta B_i|^2}{8\Delta}\Big)} |\Delta B_i|^3 \Big]^2\bigg\},
\end{align*}
which implies that
 \begin{align*}
&\mathbb{E}\Big[\sup_{0\leq k\Delta\leq T} X^{2}_{k+1}\Big]\nn\\
\leq& 8\bigg\{x^2_{0}+\mathbb{E} \sup_{0\leq k\Delta\leq T} \Big[\sum_{i=0}^{k} X_{i} \big(\beta(r_i)-a(r_i) X_{i}\big)\Delta\Big]^2
+\mathbb{E} \sup_{0\leq k\Delta\leq T}\Big[\sum_{i=0}^{k}\sigma(r_i) X_{i} \Delta B_i\Big]^2\nn\\
&+\frac{|\breve{\sigma}|^4}{4}\mathbb{E} \sup_{0\leq k\Delta\leq T}\Big[\sum_{i=0}^{k}X_{i}  \big(\Delta B_i\big)^2\Big]^2
+\Delta^2\mathbb{E} \sup_{0\leq k\Delta\leq T}\Big[\sum_{i=0}^{k}X_{i} \big(\beta(r_i)-a(r_i)X_{i}\big) \sigma(r_i) \Delta B_i\Big]^2\nn\\
&
+\frac{\Delta^4}{4}\mathbb{E} \sup_{0\leq k\Delta\leq T}\Big[\sum_{i=0}^{k}X_{i} \big(\beta(r_i)-a(r_i) X_{i}\big)^2\Big]^2 +\bar{C}^2\Delta^6\mathbb{E} \sup_{0\leq k\Delta\leq T}\Big[\sum_{i=0}^{k}X_{i}  \exp{\Big(\frac{|\Delta B_i|^2}{8\Delta}\Big)} \Big]^2\nn\\
&+\bar{C}^2\mathbb{E} \sup_{0\leq k\Delta\leq T}\Big[\sum_{i=0}^{k}X_{i}  \exp{\Big(\frac{|\Delta B_i|^2}{8\Delta}\Big)} |\Delta B_i|^3 \Big]^2\bigg\}\nn\\
\leq 8\bigg\{& x^{2}_{0}+2T |\breve{\beta}|^2\Delta   \sum_{i=0}^{\lfloor T/\Delta\rfloor} \mathbb{E}X^2_{i} +2T\check{a}^2\Delta    \sum_{i=0}^{\lfloor T/\Delta\rfloor} \mathbb{E}  X^4_{i} +\mathbb{E} \sum_{i=0}^{\lfloor T/\Delta\rfloor}\sigma^2(r_i) X^2_{i} \Delta
+2|\breve{\beta}|^4T\Delta^3 \sum_{i=0}^{\lfloor T/\Delta\rfloor} \mathbb{E}X^2_{i}\nn\\
& +\frac{|\breve{\sigma}|^4}{4}\lfloor T/\Delta\rfloor\mathbb{E} \Big[\sum_{i=0}^{\lfloor T/\Delta\rfloor}X^2_{i} \mathbb{E}\Big( \big(\Delta B_i\big)^4\big|\mathcal{F}_{t_i}\Big)\Big]
+\Delta^2\mathbb{E}  \Big[\sum_{i=0}^{\lfloor T/\Delta\rfloor}X^2_{i} \big(\beta(r_i)-a(r_i)X_{i}\big)^2 \sigma^2(r_i) \Delta \Big] \nn\\
&
 +2\check{a}^4 T\Delta^3  \sum_{i=0}^{\lfloor T/\Delta\rfloor}  \mathbb{E}X^6_{i}
 +\sqrt{2}\bar{C}^2T\Delta^5\sum_{i=0}^{\lfloor T/\Delta\rfloor}\mathbb{E} X^2_{i}   +C\lfloor T/\Delta\rfloor  \sum_{i=0}^{\lfloor T/\Delta\rfloor}\mathbb{E}X^2_{i}  \Delta^3  \bigg\}\nn\\
\leq C\Big(& x^{2}_{0}+ \Delta   \sum_{i=0}^{\lfloor T/\Delta\rfloor} \mathbb{E}X^2_{i} + \Delta    \sum_{i=0}^{\lfloor T/\Delta\rfloor} \mathbb{E}  X^4_{i}       + \Delta^3  \sum_{i=0}^{\lfloor T/\Delta\rfloor}  \mathbb{E}X^6_{i}
  \Bigg).
\end{align*}
Using Lemma \ref{log_le*6.10} we obtain
$
\mathbb{E}\Big[\sup\limits_{0\leq k\Delta\leq T} X^{2}_{k}\Big]
\leq  C.
$
The proof is therefore complete.\eproof

{\bf  Proof of Lemma \ref{log_le*6.13}.}~~Due to $\vartheta<1$, by \eqref{logeq:5.45} we have
\begin{align}\la{log_eq*6.27}
 \big(1+\mathrm{e}^{-\bar{Z}_{k+1}}\big)^{\frac{\vartheta}{2}}
\leq&\big(1+\mathrm{e}^{-Z_{k}}\big)^{\frac{\vartheta}{2}} \bigg[1+ \frac{\vartheta}{2} \varsigma_k  +\frac{\vartheta(\vartheta-2)(\vartheta-4)}{48} \varsigma_k^3 \bigg].\tag{B.9}
\end{align}
 Obviously,
\begin{align}\la{log_eq*6.28}
 \!\! \big(1+\mathrm{e}^{-Z_{k}}\big)^{\frac{\vartheta}{2}}
 = \! \big(1+\mathrm{e}^{-Z_{k}}\big)^{\frac{\vartheta}{2}}I_{\Omega^c_k}+\big(1+\mathrm{e}^{-Z_{k}}\big)^{\frac{\vartheta}{2}}I_{\Omega_k}
\leq\! \big(1+\mathrm{e}^{-\bar{Z}_{k}}\big)^{\frac{\vartheta}{2}} \!+\big(1+K^{-1}\Delta^{\theta }\big)^{\frac{\vartheta}{2}}I_{\Omega_k},\!
\tag{B.10}
\end{align}
where $\Omega_k$ is defined by \eqref{Omega}.
Using \eqref{log_eq*6.27} and \eqref{log_eq*6.28} yields
\begin{align*}
 \big(1+\mathrm{e}^{-Z_{k+1}}\big)^{\frac{\vartheta}{2}}
\!\leq&\big(1+\mathrm{e}^{-Z_{k}}\big)^{\frac{\vartheta}{2}} \bigg[1+ \frac{\vartheta}{2}\varsigma_k   + \frac{\vartheta(\vartheta-2)(\vartheta-4)}{48} \varsigma_k^3 \bigg] +2I_{\Omega_k} \nn\\
\leq&\big(1 +x^{-1}_{0} \big)^{\frac{\vartheta}{2}} +\frac{\vartheta}{2}\sum^{k}_{i=0} \big(1 +\mathrm{e}^{-Z_{i}}\big)^{\frac{\vartheta}{2}} \Big(\varsigma_k    +  \frac{(\vartheta-2)(\vartheta -4)}{24} \varsigma_k^3 \Big)
 +2\sum^{k}_{i=0}I_{\Omega_k},
\end{align*}
which implies that
\begin{align*}
 \big(1+\mathrm{e}^{-Z_{k+1}}\big)^{\vartheta}
\leq&4\big(1+x^{-1}_{0} \big)^{\vartheta}+\frac{\vartheta^2(\vartheta-2)^2(\vartheta-4)^2}{24^2} \Big[\sum^{k}_{i=0}\big(1+\mathrm{e}^{-Z_{i}}\big)^{\frac{\vartheta}{2}} \varsigma_k^3 \Big]^2\nn\\
&+ \vartheta^2 \Big[\sum^{k}_{i=0}\big(1+\mathrm{e}^{-Z_{i}}\big)^{\frac{\vartheta}{2}} \varsigma_k  \Big]^2 +16\Big[\sum^{k}_{i=0}I_{\Omega_k}\Big]^2.
\end{align*}
Then we have
\begin{align}\la{log_eq*6.29}
\mathbb{E}\sup_{0\leq k\Delta\leq T} \big(1+\mathrm{e}^{-Z_{k+1}}\big)^{\bar{p}}
\leq&4 \mathbb{E}\sup_{0\leq k\Delta\leq T}\Big[\sum^{k}_{i=0}\big(1+\mathrm{e}^{-Z_{i}}\big)^{\frac{\bar{p}}{2}} \varsigma_k  \Big]^2 +\mathbb{E}\sup_{0\leq k\Delta\leq T}\Big[\sum^{k}_{i=0}\big(1+\mathrm{e}^{-Z_{i}}\big)^{\frac{\bar{p}}{2}} \varsigma_k^3 \Big]^2\nn\\
&+4\big(1+x^{-1}_{0} \big)^{\bar{p}}
+16\mathbb{E}\sup_{0\leq k\Delta\leq T}\Big[\sum^{k}_{i=0}I_{\Omega_k}\Big]^2,\tag{B.11}
\end{align}
we deduce that
\begin{align}\la{log_eq*6.31}
  &\mathbb{E}\sup_{0\leq k\Delta\leq T}\Big[\sum^{k}_{i=0}\big(1+\mathrm{e}^{-Z_{i}}\big)^{\frac{\vartheta}{2}} \varsigma_k  \Big]^2\nn\\
\leq&7\check{a}^2\Delta^2\mathbb{E}\sup_{0\leq k\Delta\leq T}\bigg(\sum^{k}_{i=0}\big(1+\mathrm{e}^{-Z_{i}}\big)^{\frac{\vartheta}{2}-1} \bigg)^2+7|\breve{\beta}|^2\Delta^2\mathbb{E}\sup_{0\leq k\Delta\leq T}\bigg(\sum^{k}_{i=0}\big(1+\mathrm{e}^{-Z_{i}}\big)^{\frac{\vartheta}{2}}  \bigg)^2\nn\\
&+C \Delta^{4(1-\theta )}\mathbb{E}\sup_{0\leq k\Delta\leq T}\bigg(\sum^{k}_{i=0}\big(1+\mathrm{e}^{-Z_{i}}\big)^{\frac{\vartheta}{2}}  \bigg)^2
+4|\breve{\sigma}|^2 \mathbb{E}\sup_{0\leq k\Delta\leq T}\bigg(\sum^{k}_{i=0}\big(1+\mathrm{e}^{-Z_{i}}\big)^{\frac{\vartheta}{2}}  (\Delta B_i)^2\bigg)^2 \nn\\
& +7 \mathbb{E}\sup_{0\leq k\Delta\leq T}\Big|\sum^{k}_{i=0}\sigma(r_i)\big(1+\mathrm{e}^{-Z_{i}}\big)^{\frac{\vartheta}{2}-1}\mathrm{e}^{-Z_{k}} \Delta B_i \Big|^2+7 |\breve{\sigma}|^2  \mathbb{E}\sup_{0\leq k\Delta\leq T}\bigg(\sum^{k}_{i=0} \big(1+\mathrm{e}^{-Z_{i}}\big)^{\frac{\vartheta}{2}} \mathcal{U}_k \bigg)^2\nn\\
&+7\Delta^2 \mathbb{E}\sup_{0\leq k\Delta\leq T}\Big|\sum^{k}_{i=0}\sigma(r_i)\big(1+\mathrm{e}^{-Z_{i}}\big)^{\frac{\vartheta}{2}-1}\big(a(r_i) -\beta(r_i)\mathrm{e}^{-Z_{i}}\big)\Delta B_i\Big|^2\nn\\
%
\leq&7T^2\check{a}^2  +7T|\breve{\beta}|^2\Delta \sum^{k}_{i=0} \mathbb{E} \big(1+\mathrm{e}^{-Z_{i}}\big)^{\vartheta}+CT \Delta^{3-4\theta }  \sum^{\lfloor T/\Delta\rfloor}_{i=0}\mathbb{E}\big(1+\mathrm{e}^{-Z_{i}}\big)^{ \vartheta } \nn\\
&
+12T|\breve{\sigma}|^2 \Delta \sum^{\lfloor T/\Delta\rfloor}_{i=0}\mathbb{E}\big(1+\mathrm{e}^{-Z_{i}}\big)^{\vartheta} +7|\breve{\sigma}|^2\Delta \sum^{\lfloor T/\Delta\rfloor}_{i=0}\mathbb{E} \big(1+\mathrm{e}^{-Z_{i}}\big)^{\vartheta}+14T\check{a}^2|\breve{\sigma}|^2\Delta^2   \nn\\
& +14|\breve{\sigma}|^2 |\breve{\beta}|^2\Delta^3 \sum^{\lfloor T/\Delta\rfloor}_{i=0}\mathbb{E}\big(1+\mathrm{e}^{-Z_{i}}\big)^{\vartheta}    +C  T \Delta^2  \sum^{\lfloor T/\Delta\rfloor}_{i=0}\mathbb{E}\big(1+\mathrm{e}^{-Z_{i}}\big)^{\vartheta},\tag{B.12}
\end{align}
and
\begin{align}\la{log_eq*6.32}
\mathbb{E}\sup_{0\leq k\Delta\leq T}\Big[\sum^{k}_{i=0}\big(1+\mathrm{e}^{-Z_{i}}\big)^{\frac{\vartheta}{2}} \varsigma^3_k  \Big]^2
   \leq &\lfloor T/\Delta\rfloor\mathbb{E} \Big[\sum^{\lfloor T/\Delta\rfloor}_{i=0}\big(1+\mathrm{e}^{-Z_{i}}\big)^{ \vartheta } \mathbb{E}\big(\varsigma^6_k|\mathcal{G}_{t_k}\big)  \Big]\nn\\
   \leq &CT\Delta \sum^{\lfloor T/\Delta\rfloor}_{i=0}\mathbb{E} \big(1+\mathrm{e}^{-Z_{i}}\big)^{ \vartheta }.\tag{B.13}
\end{align}
By Chebyshev's inequality and   Lemma \ref{log_le*6.10},
\begin{align}\la{log_eq*6.30}
 \mathbb{E}\sup_{0\leq k\Delta\leq T}\Big[\sum^{k}_{i=0}I_{\Omega_k}\Big]^2
\leq& \lfloor T/\Delta\rfloor\sum^{\lfloor T/\Delta\rfloor}_{i=0}\mathbb{P}\Big\{\bar{Z}_{i+1}> \log(K\Delta^{-\theta })\Big\}\nn\\
\leq& \lfloor T/\Delta\rfloor\sum^{\lfloor T/\Delta\rfloor}_{i=0}\frac{\mathbb{E}\exp{\big(2\theta ^{-1}\bar{Z}_{i+1}\big)}}{K^{2\theta ^{-1}}\Delta^{-2}} \leq \frac{C T^2}{K^{2\theta ^{-1}} },\tag{B.14}
\end{align}
Inserting \eqref{log_eq*6.31}, \eqref{log_eq*6.32} and  \eqref{log_eq*6.30} into \eqref{log_eq*6.29}, and using Lemma \ref{log:Le5.9} we obtain
\begin{align*}
\mathbb{E}\sup_{0\leq k\Delta\leq T} \big(1+\mathrm{e}^{-Z_{k}}\big)^{\vartheta}
\leq C_T.
\end{align*}
The proof is therefore complete.\eproof

\newpage

\subsection*{Appendix C.}
{\bf  Proof of Lemma \ref{appendis_L1}.}~~We first note that
\begin{align*}
|X^{\Delta}_{1}|\geq& x_0\min_{i\in \mathbb{S}}\{|\sigma(i)|\}|\Delta B_0|-x_0\Big(1+ \check{b} \Delta+\check{a}  x_0 \Delta \Big)\geq \frac{\mathrm{e}}{\Delta}
\end{align*}
if
$$
|\Delta B_0|\geq \frac{\mathrm{e}/\Delta+x_0\big(1+ \check{b} \Delta+\check{a}  x_0 \Delta \big)}{x_0\min_{i\in \mathbb{S}}\{|\sigma(i)|\}}
$$
In other words, we have
\begin{align}\la{app_eq0}
\mathbb{P}\Big(|X^{\Delta}_{1}|\geq\frac{\mathrm{e}}{\Delta}\Big)
\geq \mathbb{P}\Big(|\Delta B_0|\geq \frac{\mathrm{e}/\Delta+x_0\big(1+ \check{b} \Delta+\check{a}  x_0 \Delta \big)}{x_0\min_{i\in \mathbb{S}}\{|\sigma(i)|\}}\Big)>0  \tag{C.1}
\end{align}
due to  $\Delta B_0\sim\mathcal{N}(0,\Delta)$. Let $M_i= \big(a(i)-1.4\big)/|\sigma(i)|$, we observe that, for $k\geq 1$, if $|X^{\Delta}_{k}|\geq  {\exp(2^{k-1})}/{\Delta}$ and $|\Delta B_k|\leq  \hat{M}\exp(2^{k-1})$ hold, then
$$
|X^{\Delta}_{k+1}|\geq \frac{\exp(2^{k})}{\Delta}.
$$
In fact,
\begin{align*}
|X^{\Delta}_{k+1}|=&|X^{\Delta}_{k}|\big|1+b(r_k)\Delta
-a(r_k)X^{\Delta}_{k}\Delta+\sigma(r_k)\Delta B_k\big|\nn\\
\geq&|X^{\Delta}_{k}|\Big[a(r_k)|X^{\Delta}_{k}|\Delta-1-b(r_k)\Delta-|\sigma(r_k)||\Delta B_k|
\Big]\nn\\
\geq&\frac{\exp(2^{k-1})}{\Delta}\Big[a(r_k)\exp(2^{k-1})-1-b(r_k)\bar{\Delta}
-|\sigma(r_k)|\hat{M}\exp(2^{k-1})
\Big]\nn\\
=&\frac{\exp(2^{k})}{\Delta}\Big(a(r_k)-\exp(-2^{k-1})-b(r_k)\exp(-2^{k-1})\bar{\Delta}
-|\sigma(r_k)|\hat{M}
\Big)\nn\\
\geq&\frac{\exp(2^{k})}{\Delta}\Big(a(r_k)-\mathrm{e}^{-1}-b(r_k)\mathrm{e}^{-1}\bar{\Delta}
-|\sigma(r_k)|\hat{M}
\Big)\geq \frac{\exp(2^{k})}{\Delta}.
\end{align*}
We therefore have
$$
\Big\{|X^{\Delta}_{1}|\geq  \frac{\mathrm{e}}{\Delta}~\mathrm{and}~
|\Delta B_k|\leq  \hat{M}\exp(2^{k-1}) ,~\forall~k\geq 1\Big\}\subset \Big\{|X^{\Delta}_{k}|\geq  \frac{\exp(2^{k-1})}{\Delta},~\forall~k\geq 1\Big\}.
$$
Since $X^{\Delta}_{1}$ and $\Delta B_k$ for $k\geq 1$ are all independent,
\begin{align*}
\mathbb{P}\bigg(|X^{\Delta}_{k}|\geq  \frac{\exp(2^{k-1})}{\Delta},~\forall~k\geq 1\bigg)
\geq& \mathbb{P}\bigg(|X^{\Delta}_{1}|\geq  \frac{\mathrm{e}}{\Delta}~\mathrm{and}~
|\Delta B_k|\leq \hat{M}\exp(2^{k-1}),~\forall~k\geq 1\bigg)\nn\\
=&\mathbb{P}\Big(|X^{\Delta}_{1}|\geq  \frac{\mathrm{e}}{\Delta}\Big)
\mathbb{P}\Big(|\Delta B_k|\leq  \hat{M}\exp(2^{k-1}),~\forall~k\geq 1\Big).
\end{align*}
By the conditional probability formula implies
\begin{align}\la{app-eq1}
\mathbb{P}\bigg(|X^{\Delta}_{k+1}|\geq  \frac{\exp(2^{k})}{\Delta},~\forall~k\geq 1\bigg||X^{\Delta}_{1}|\geq  \frac{\mathrm{e}}{\Delta}\bigg)
=&\frac{\mathbb{P}\Big(|X^{\Delta}_{k}|\geq  \frac{\exp(2^{k-1})}{\Delta},~\forall~k\geq 1\Big)}{\mathbb{P}\Big(|X^{\Delta}_{1}|\geq  \frac{\mathrm{e}}{\Delta}\Big)}\nn\\
\geq&
\mathbb{P}\Big(|\Delta B_k|\leq  \hat{M}\exp(2^{k-1}),~\forall~k\geq 1\Big)\nn\\
=&\prod_{k=1}^{\infty}
\mathbb{P}\Big(|\Delta B_k|\leq  \hat{M}\exp(2^{k-1}) \Big).\tag{C.2}
\end{align}
Now, because $\Delta B_k\sim\mathcal{N}(0,\Delta)$, we have
\begin{align*}
\mathbb{P}\Big(|\Delta B_k|> \hat{M}\exp(2^{k-1}) \Big)
=&\mathbb{P}\Big(\frac{|\Delta B_k|}{\sqrt{\Delta}}> \frac{\hat{M}\exp(2^{k-1})}{\sqrt{\Delta}}\Big)
\leq \mathbb{P}\Big(\frac{|\Delta B_k|}{\sqrt{\Delta}}> \hat{M}\exp(2^{k-1})\Big)\nn\\
=&\frac{2}{\sqrt{2\pi}}\int_{\hat{M}\exp(2^{k-1})}^{\infty}
\mathrm{e}^{-\frac{x^2}{2}}\mathrm{d}x
\leq \int_{\hat{M}\exp(2^{k-1})}^{\infty}
\frac{x\mathrm{e}^{-\frac{x^2}{2}}}{\hat{M}\exp(2^{k-1})}
\mathrm{d}x\nn\\
\leq&\frac{2}{\exp(2^{k-1})}\int_{\hat{M}\exp(2^{k-1})}^{\infty}
 x\mathrm{e}^{-\frac{x^2}{2}}
\mathrm{d}x\nn\\
=& 2
\exp\Big(-2^{k-1}-\hat{M}^2\exp(2^{k})/2 \Big)\nn\\
\leq& 2
\exp\Big(-2^{k-1}-2^{-3}\exp(2^{k}) \Big).
\end{align*}
But, by the elementary inequality $\log(1-u)\geq-2u$ for $0\leq u < 0.5$, we derive
\begin{align*}
&\log\bigg(\prod^{\infty}_{k=1}\bigg[1-2
\exp\Big(-2^{k-1}-2^{-3}\exp(2^{k}) \Big)\bigg]\bigg)\nn\\
=&\sum^{\infty}_{k=1}\log\bigg( 1-
2
\exp\big(-2^{k-1}-2^{-3}\exp(2^{k}) \big) \bigg)
\geq -4\sum^{\infty}_{k=1}
\exp\big(-2^{k-1}-2^{-3}\exp(2^{k}) \big).
\end{align*}
Noting that $\mathrm{e}^{x}\geq1+x$ and $2^{k-1}\geq2(k-1)$, we then get
\begin{align*}
&-4\sum^{\infty}_{k=1}
\exp\big(-2^{k-1}-2^{-3}\exp(2^{k}) \big)\nn\\
\geq&-4\sum^{\infty}_{k=1}\exp\big(-2^{k-1}-2^{-3}- 2^{k-3} \big)
\geq -4\sum^{\infty}_{k=1}\exp\big(-2(k-1)-2^{-3}- 2(k-3) \big)\nn\\
=& -4\exp\big(8-2^{-3}  \big)\sum^{\infty}_{k=1}\mathrm{e}^{-4k}
= -4\exp\big(8-2^{-3}  \big)
\lim_{k\rightarrow \infty}
\frac{\mathrm{e}^{-4}(1-\mathrm{e}^{-4k})}{1-\mathrm{e}^{-4}}
=-
\frac{4\mathrm{e}^{4-2^{-3}}}{1-\mathrm{e}^{-4}}.
\end{align*}
Hence, in \eqref{app-eq1},
\begin{align*}
\log\bigg(\mathbb{P}\bigg(|X^{\Delta}_{k+1}|\geq  \frac{\exp(2^{k})}{\Delta},~\forall~k\geq 1\bigg||X^{\Delta}_{1}|\geq  \frac{\mathrm{e}}{\Delta}\bigg)\bigg)\geq-
\frac{4\mathrm{e}^{4-2^{-3}}}{1-\mathrm{e}^{-4}}
\end{align*}
and the result follows.\eproof
\subsection*{Appendix D.}
We can easily obtain from Lemma \ref{log_le4} the following corollary.
\begin{corollary}\la{SDE_log_le4}
For any $p>0$, the EM scheme defined by \eqref{SDE_log_eq5} has the property that
\begin{align*}
\sup_{\Delta\in (0, 1]}\sup_{0\leq k\leq \lfloor T/\Delta\rfloor}\mathbb{E}\big[\mathrm{e}^{pY_{k}}\big]
\leq  C_T,~~~~~\forall~T>0,
\end{align*}
  where $\lfloor T/\Delta\rfloor$  represents the integer part of $T/\Delta$.
\end{corollary}

{\bf  Proof of Lemma \ref{SDE_log_th2}.}~~By \eqref{log_eq4}, we have
\begin{align}\la{SDE_log_eq11}
y(t_{k+1})
 =& y(t_{k})+\beta\Delta-a\int_{t_k}^{t_{k+1}} \mathrm{e}^{y(s)} \mathrm{d}s+ \sigma  \Delta B_k.\tag{D.1}
\end{align}
Using
\eqref{SDE_log_eq5} and \eqref{SDE_log_eq11} we have
\begin{align*}
Y_{k+1}-y(t_{k+1})
=&Y_{k}-y(t_{k}) -a\big(\mathrm{e}^{Y_k} -\mathrm{e}^{y(t_k)}\big)\Delta+a \int_{t_k}^{t_{k+1}} \int_{t_k}^{s}\mathrm{d}x(u) \mathrm{d}s\nn\\
=:&Y_{k}-y(t_{k}) -a\big(\mathrm{e}^{Y_k} -\mathrm{e}^{y(t_k)}\big)\Delta+\Xi_{k},
\end{align*}
where
\begin{align}\la{SDE_log_eq13}
 \Xi_{k}
 \! =a\int_{t_k}^{t_{k+1}}\int_{t_k}^{s}x(u) \big(b-ax(u)\big)\mathrm{d}u\mathrm{d}s
  + a\sigma\int_{t_k}^{t_{k+1}}\int_{t_k}^{s}x(u)\mathrm{d}B(u) \mathrm{d}s
  =:\Xi^{(1)}_{k}+\Xi^{(2)}_{k}.\tag{D.2}
\end{align}
Let us define $u_k=Y_{k}-y(t_{k})$. Note that  $u_k\big(\mathrm{e}^{Y_k} -\mathrm{e}^{y(t_k)}\big)\geq 0$, we get
\begin{align}\la{SDE_log_eq14}
u^2_{k+1}
\leq& u^2_{k}
-2au_k\big(\mathrm{e}^{Y_k} -\mathrm{e}^{y(t_k)}\big)\Delta+2a^2 \big(\mathrm{e}^{Y_k} -\mathrm{e}^{y(t_k)}\big)^2\Delta^2 +2\Xi^2_{k }+2u_k\Xi_{k }\nn\\
\leq& u^2_{k}+2a^2 \big(\mathrm{e}^{Y_k} -\mathrm{e}^{y(t_k)}\big)^2\Delta^2 +2\Xi^2_{k }+2u_k\Xi_{k }\nn\\
\leq&  2a^2\Delta^2 \sum_{i=0}^{k}\big(\mathrm{e}^{Y_i} -\mathrm{e}^{y(t_i)}\big)^2 +2\sum_{i=0}^{k}\Xi^2_{i }+2\sum_{i=0}^{k}u_{i}\Xi_{i }.\tag{D.3}
\end{align}
Let $\mathfrak{M}_0=0,$ and
$
\mathfrak{M}_k=\sum\limits_{i=0}^{k-1}u_{i}\Xi^{(2)}_{i }
$  for any $k\geq 1$,
since
$$
\mathbb{E}\Big[\Xi^{(2)}_{k }|{\cal{F}}_{t_k}\Big]
=a\sigma\mathbb{E}\bigg[\int_{t_k}^{t_{k+1}}\int_{t_k}^{s}x(u)\mathrm{d}B(u) \mathrm{d}s\Big|{\cal{F}}_{t_k}\bigg]=0.
$$
It is then easy to show that
$$
\mathbb{E}\Big[\mathfrak{M}_{k+1}|{\cal{F}}_{t_k}\Big]
= \mathbb{E}\Big[\mathfrak{M}_{k}+u_k\Xi^{(2)}_{k }|{\cal{F}}_{t_k}\Big]
= \mathfrak{M}_{k}+u_k\mathbb{E}\Big[\Xi^{(2)}_{k }|{\cal{F}}_{t_k}\Big]
=\mathfrak{M}_{k},
$$
This implies immediately that $\mathfrak{M}_{k}$ is a martingale and the Burkholder-Davis-Gundy inequality implies that
\begin{align*}
 \mathbb{E}\Big[\sup_{k=0,\ldots, l}\big|\mathfrak{M}_{k}\big|^q\Big]\leq C_T \mathbb{E}\bigg[\Big|\sum_{i=0}^{l-1}\big(u_{i}\Xi^{(2)}_{i }\big)^2\Big|^{\frac{q}{2}}\bigg]
 \leq C_T \mathbb{E}\bigg[\bigg(\sum_{i=0}^{l-1}u^2_{i}\big|\Xi^{(2)}_{i }\big|^2\bigg)^{\frac{q}{2}}\bigg]
\end{align*}
for any $q\geq 2$ and $l=0,\ldots, \lfloor T/\Delta\rfloor$. Using this and Jensen's inequality in \eqref{SDE_log_eq14} we now arrive at
\begin{align}\la{SDE_log_eq16}
&\mathbb{E}\bigg[\sup_{k=0,\ldots, l} |u_{k+1}|^{2q}\bigg]
\leq   2^q\mathbb{E}\bigg[\sup_{k=0,\ldots, l}\Big|a^2\Delta^2 \sum_{i=0}^{k}\big(\mathrm{e}^{Y_i} -\mathrm{e}^{y(t_i)}\big)^2 + \sum_{i=0}^{k}\Xi^2_{i }+ \sum_{i=0}^{k}u_{i}\Xi_{i }\Big|^q\bigg]\nn\\
\leq& 6^q \mathbb{E}\bigg[(a\Delta)^{2q}\bigg[\sum_{i=0}^{l}\big(\mathrm{e}^{Y_i} -\mathrm{e}^{y(t_i)}\big)^2\bigg]^q +\bigg(\sum_{i=0}^{l}\Xi^2_{i }\bigg)^q+\sup_{k=0,\ldots, l}\Big|\sum_{i=0}^{k}u_{i}\Xi_{i }\Big|^q\bigg]\nn\\
\leq& 6^q \mathbb{E}\bigg[(a\Delta)^{2q}\big(\lfloor T/\Delta\rfloor\big)^{q-1} \sum_{i=0}^{l}\big(\mathrm{e}^{Y_i} -\mathrm{e}^{y(t_i)}\big)^{2q}+\big(\lfloor T/\Delta\rfloor\big)^{q-1} \sum_{i=0}^{l}\Xi^{2q}_{i }\nn\\
&~~~~~~~~~~~~+2^{q}\big(\lfloor T/\Delta\rfloor\big)^{q-1}
\sum_{i=0}^{l}|u_{i}|^q\big|\Xi^{(1)}_{i }\big|^q +2^{q} \sup_{k=0,\ldots, l}\Big|\mathfrak{M}_{k+1}\Big|^q\bigg]\nn\\
\leq& C_T\mathbb{E}\bigg[ \Delta^{2q}\big(\lfloor T/\Delta\rfloor\big)^{q-1} \sum_{i=0}^{l}\big(\mathrm{e}^{Y_i}-\mathrm{e}^{y(t_i)}\big)^{2q}+\big(\lfloor T/\Delta\rfloor\big)^{q-1} \sum_{i=0}^{l}\big(\Xi^{(1)}_{i }\big)^{2q}\nn\\
&~~~~~ +\big(\lfloor T/\Delta\rfloor\big)^{q-1} \sum_{i=0}^{l}\big(\Xi^{(2)}_{i }\big)^{2q}+ \big(\lfloor T/\Delta\rfloor\big)^{q-1}
\sum_{i=0}^{l}|u_{i}|^q\big|\Xi^{(1)}_{i }\big|^q + \bigg(\sum_{i=0}^{l}u^2_{i}\big|\Xi^{(2)}_{i }\big|^2\bigg)^{\frac{q}{2}}\bigg]\nn\\
\leq& C_T\mathbb{E}\bigg[ \Delta^{2q}\big(\lfloor T/\Delta\rfloor\big)^{q-1} \sum_{i=0}^{l}\big(\mathrm{e}^{Y_i}\!-\mathrm{e}^{y(t_i)}\big)^{2q}\!+\big(\lfloor T/\Delta\rfloor\big)^{q-1} \sum_{i=0}^{l}\big|\Xi^{(1)}_{i }\big|^{2q}\!+\big(\lfloor T/\Delta\rfloor\big)^{q-1} \sum_{i=0}^{l}\big|\Xi^{(2)}_{i }\big|^{2q}
\nn\\
&~~~~~ + \big(\lfloor T/\Delta\rfloor\big)^{q-1}
\sum_{i=0}^{l}|u_{i}|^q\big|\Xi^{(1)}_{i }\big|^q +\big(\lfloor T/\Delta\rfloor\big)^{q/2-1} \sum_{i=0}^{l}|u_{i}|^{q}\big|\Xi^{(2)}_{i }\big|^{q} \bigg]\tag{D.4}
\end{align}
for any $q\geq 2$ and $l=0,\ldots, \lfloor T/\Delta\rfloor$. It is easy to see that
\begin{align*}
 \mathbb{E}\Big[ \Delta^{2q}\big(\lfloor T/\Delta\rfloor\big)^{q-1} \sum_{i=0}^{l}\big(\mathrm{e}^{Y_i}-\mathrm{e}^{y(t_i)}\big)^{2q}\Big]
\leq&\Delta^{q+1} T^{q-1} \sum_{i=0}^{l} \mathbb{E}\Big[\big(\mathrm{e}^{Y_i}
+\mathrm{e}^{y(t_i)}\big)^{q}\big|\mathrm{e}^{Y_i}-\mathrm{e}^{y(t_i)}\big|^{q}\Big]
.
\end{align*}
On the other hand, by applying Corollary \ref{SDE_log_le4} and \eqref{log_eq3}, we infer that
\begin{align}\la{SDE_log_eq17}
\sup_{0\leq k\leq \lfloor T/\Delta\rfloor}\mathbb{E}\Big[\mathrm{e}^{2qY_k}\Big]\leq C_T,~~~~\sup_{0\leq t\leq   T }\mathbb{E}\Big[\mathrm{e}^{2qy(t)}\Big]=\sup_{0\leq t\leq  T} \mathbb{E}\Big[x^{2q}(t)\Big]\leq C_T.\tag{D.5}
\end{align}
Now note that
\begin{align}\la{SDE_log_eq18}
 \mathbb{E}\Big[ \Delta^{2q}\big(\lfloor T/\Delta\rfloor\big)^{q-1} \sum_{i=0}^{l}\big(\mathrm{e}^{Y_i}-\mathrm{e}^{y(t_i)}\big)^{2q}\Big]
\leq&\Delta^{q+1}T^{q-1}  \sum_{i=0}^{l} \mathbb{E}\Big[\big(\mathrm{e}^{Y_i}+\mathrm{e}^{y(t_i)}\big)^{2q}\big|Y_i-y(t_i)\big|^{q}\Big]\nn\\
\leq& C_T\Delta^{2q}+C_T\Delta\sum_{i=0}^{l}\mathbb{E}|u_{i}|^{2q}.\tag{D.6}
\end{align}
By \eqref{SDE_log_eq13} and \eqref{SDE_log_eq17}, for any $q\geq 2$,
\begin{align}\la{SDE_log_eq19}
\mathbb{E}\big|\Xi^{(1)}_{i }\big|^{2q}
\leq& a^{2q} \mathbb{E}\bigg[\bigg(\int_{t_k}^{t_{k+1}}\int_{t_k}^{s}x(u) \big(b+ax(u)\big)\mathrm{d}u\mathrm{d}s\bigg)^{2q}\bigg]\nn\\
\leq& a^{2q} \mathbb{E}\bigg[\bigg(\int_{t_k}^{t_{k+1}}1^{\frac{2q}{2q-1}}
ds\bigg)^{2q-1}\bigg(\int_{t_k}^{t_{k+1}}\bigg(\int_{t_k}^{s}x(u) \big(b+ax(u)\big)\mathrm{d}u\bigg)^{2q}\mathrm{d}s\bigg)\bigg]\nn\\
\leq&(a\vee b)^{2q} (2a)^{2q} \Delta^{4q-2} \int_{t_k}^{t_{k+1}}
\int_{t_k}^{s}\Big( \mathbb{E}\big[x^{2q}(u)\big]  +\mathbb{E}\big[x^{4q}(u)\big] \Big)\mathrm{d}u\mathrm{d}s
\leq C_T\Delta^{4q},\tag{D.7}
\end{align}
and
\begin{align}\la{SDE_log_eq20}
\mathbb{E}\big|\Xi^{(2)}_{i }\big|^{2q}
\leq& (a \sigma)^{2q} \Delta^{2q-1} \int_{t_k}^{t_{k+1}}\mathbb{E}\bigg[
\bigg|\int_{t_k}^{s}x(u)\mathrm{d}B(u)\bigg|^{2q}\bigg]\mathrm{d}s\nn\\
 \leq& C_T \Delta^{3q-2} \int_{t_k}^{t_{k+1}}
 \int_{t_k}^{s}\mathbb{E}\big[x^{2q}(u)\big]\mathrm{d}u\mathrm{d}s
  \leq  C_T \Delta^{3q}.\tag{D.8}
\end{align}
Thus, the Canchy-Schwarz inequality give that
\begin{align}\la{SDE_log_eq21}
\mathbb{E}\big[|u_{i}|^q\big|\Xi^{(1)}_{i }\big|^q\big]\leq \Big(\mathbb{E}\big[|u_{i}|^{2q}\big]\Big)^{1/2}
\Big(\mathbb{E}\big[\big|\Xi^{(1)}_{i }\big|^{2q}\big]\Big)^{1/2}\leq C_T \Big(\mathbb{E}\big[|u_{i}|^{2q}\big]\Big)^{1/2}\Delta^{2q}.\tag{D.9}
\end{align}
Similar we also obtain
\begin{align}\la{SDE_log_eq22}
\mathbb{E}\big[|u_{i}|^q\big|\Xi^{(2)}_{i }\big|^q\big]\leq \Big(\mathbb{E}\big[|u_{i}|^{2q}\big]\Big)^{1/2}
\Big(\mathbb{E}\big[\big|\Xi^{(2)}_{i }\big|^{2q}\big]\Big)^{1/2}\leq C_T \Big(\mathbb{E}\big[|u_{i}|^{2q}\big]\Big)^{1/2}\Delta^{3q/2}.\tag{D.10}
\end{align}
Thus, for any integer $k\geq 0$,
substituting \eqref{SDE_log_eq18}-\eqref{SDE_log_eq22} into \eqref{SDE_log_eq16}, we know that
\begin{align*}
&\mathbb{E}\bigg[\sup_{k=0,\ldots, l} |u_{k+1}|^{2q}\bigg]\nn\\
\leq& C_T\bigg\{T^{q}\Delta^{2q}\!+T^{q-1}\Delta\sum_{i=0}^{l}\mathbb{E}|u_i|^{2q} +\!\big(\lfloor T/\Delta\rfloor\big)^{q-1}
\sum_{i=0}^{l}\mathbb{E}\Big[ \big|\Xi^{(1)}_{i }\big|^{2q}\Big]\!+\!\big(\lfloor T/\Delta\rfloor\big)^{q-1} \sum_{i=0}^{l}\mathbb{E}\Big[ \big|\Xi^{(2)}_{i }\big|^{2q}\Big]\nn\\
&~~ + \big(\lfloor T/\Delta\rfloor\big)^{q-1}
\sum_{i=0}^{l}\mathbb{E}\Big[ |u_{i}|^q\big|\Xi^{(1)}_{i }\big|^q \Big]+\big(\lfloor T/\Delta\rfloor\big)^{q/2-1} \sum_{i=0}^{l}\mathbb{E}\Big[ |u_{i}|^{q}\big|\Xi^{(2)}_{i }\big|^{q}\Big] \bigg\}\nn\\
\leq& C_T\bigg\{T^{q}\Delta^{2q}+T^{q-1}\Delta\sum_{i=0}^{l}\mathbb{E}|u_{i}|^{2q} +\big(\lfloor T/\Delta\rfloor\big)^{q}\Delta^{4q}
+\big(\lfloor T/\Delta\rfloor\big)^{q}\Delta^{3q}\nn\\
&~~~~~+ \big(\lfloor T/\Delta\rfloor\big)^{q-1}
\sum_{i=0}^{l} \Big(\mathbb{E} |u_{i}|^{2q} \Big)^{1/2}\Delta^{2q} +\big(\lfloor T/\Delta\rfloor\big)^{q/2-1} \sum_{i=0}^{l}
\Big(\mathbb{E} |u_{i}|^{2q} \Big)^{1/2}\Delta^{3q/2} \bigg\}\nn\\
\leq& C\bigg\{2T^{q}\Delta^{2q}+T^{q-1}\Delta\sum_{i=0}^{l}
\mathbb{E}|u_{i}|^{2q}
 +\big(T^{q-1}+T^{q/2-1}\big)
\sum_{i=0}^{l} \Big(\mathbb{E} |u_{i}|^{2q} \Big)^{1/2}\Delta^{q+1}
  \bigg\}\nn\\
\leq& C_T\bigg\{2T^{q}\Delta^{2q}+T^{q-1}\Delta\sum_{i=0}^{l}\mathbb{E}|u_{i}|^{2q}
 +\big(T^{q-1}+T^{q/2-1}\big)\Delta
\sum_{i=0}^{l} \Big(\mathbb{E} |u_{i}|^{2q} +\Delta^{2q}\Big)  \bigg\}\nn\\
\leq& C_T\bigg(3T^{q}\Delta^{2q}+2T^{q-1}\Delta\sum_{i=0}^{l}\mathbb{E}|u_{i}|^{2q} +T^{q/2-1}\Delta \sum_{i=0}^{l}  \mathbb{E} |u_{i}|^{2q} +T^{q/2}\Delta^{2q}   \bigg)\nn\\
\leq& C_T\bigg(T^{q}\Delta^{2q}+T^{q-1}\Delta\sum_{i=0}^{l}\mathbb{E}|u_{i}|^{2q} \bigg)
\end{align*}
for any $q\geq 2$ and $l=0,\ldots, \lfloor T/\Delta\rfloor$. By Gronwall's Lemma
\begin{align*}
\mathbb{E}\bigg[\sup_{k=0,\ldots, l} |u_{k+1}|^{2q}\bigg]
\leq& C_T\Delta^{2q}\exp(C_T)
\end{align*}
for any $q\geq 2$ and $l=0,\ldots, \lfloor T/\Delta\rfloor$. This completes now the proof of the assertion for $q\geq 2$. The case $q\in (0, 2)$ follows now by Lyapunov's inequality.
 \eproof


\end{document}